\title{Local time on the exceptional set of dynamical percolation, and the Incipient Infinite Cluster}
\author{Alan Hammond \and G\'abor Pete \and Oded Schramm}
\date{}

\documentclass[11pt]{article}

\usepackage{amsmath}
\usepackage{amssymb}
\usepackage{amsthm}
\usepackage{amsfonts}
\usepackage{graphicx}
\usepackage{mathrsfs}

\usepackage{colordvi}
\usepackage[usenames,dvipsnames]{xcolor}

\newif\iffigures\figurestrue
\figuresfalse

\newif\ifhyper\IfFileExists{hyperref.sty}{\hypertrue}{\hyperfalse}

\ifhyper\usepackage{hyperref}
\def\hitem#1#2{\item[\hypertarget{#1}{#2}]\expandafter\gdef\csname LBL#1ITM\endcsname{#2}}
\def\iref#1{\hyperlink{#1}{\csname LBL#1ITM\endcsname}}
\else
\def\hitem#1#2{\item[{#2}]\expandafter\gdef\csname LBL#1ITM\endcsname{#2}}
\def\iref#1{{\csname LBL#1ITM\endcsname}}
\fi

\newif\ifdraft
\drafttrue
\long\def\note#1/{\ifdraft {\bf [#1]}\fi}
\long\def\comment#1{}
\long\def\old#1{}

\def\P{\mathbb{P}}
\def\E{\mathbb{E}}
\def\Bb#1#2{{\def\md{\bigm| }#1\bigl(#2\bigr)}}
\def\BB#1#2{{\def\md{\Bigm| }#1\Bigl(#2\Bigr)}}
\def\Bs#1#2{{\def\md{\mid}#1(#2)}}
\def\Pb{\Bb\P}
\def\Eb{\Bb\E}
\def\PB{\BB\P}
\def\EB{\BB\E}
\def\Ps{\Bs\P}
\def\Es{\Bs\E}
\def\Pso#1{\Bs{\P_{#1}}}
\def\Pbo#1{\Bb{\P_{#1}}}
\def\PBo#1{\BB{\P_{#1}}}
\def\Eso#1{\Bs{\E_{#1}}}
\def\Ebo#1{\Bb{\E_{#1}}}
\def\EBo#1{\BB{\E_{#1}}}

\newcommand{\p}{\partial}
\newcommand{\dist}{\mathrm{dist}}
\def\d{\mathrm{d}}

\def \proof {{ \medbreak \noindent {\bf Proof. }}}
\def\proofof#1{{ \medbreak \noindent {\bf Proof of #1.}}}
\def\proofcont#1{{ \medbreak \noindent {\bf Proof of #1, continued.}}}

\numberwithin{equation}{section}
\numberwithin{figure}{section}

\newtheorem{theorem}{Theorem}[section]
\newtheorem{corollary}[theorem]{Corollary}
\newtheorem{lemma}[theorem]{Lemma}
\newtheorem{proposition}[theorem]{Proposition}
\newtheorem{conjecture}[theorem]{Conjecture}
\newtheorem{question}[theorem]{Question}

\newtheorem{definition}[theorem]{Definition}
\theoremstyle{remark}

\newcommand{\supp}{\mathop{\mathrm{supp}}}

\newcommand{\eps}{\epsilon}
\newcommand{\R}{\mathbb R}
\newcommand{\N}{\mathbb N}
\newcommand{\Q}{\mathbb Q}

\newcommand{\Z}{\mathbb Z}

\newcommand{\cA} {{\mathcal A}}
\newcommand{\cB} {{\mathcal B}}

\newcommand{\cN} {{\mathcal N}}
\newcommand{\cY} {{\mathcal Y}}

\newcommand{\cR} {{\mathcal R}}
\newcommand{\cM} {{\mathcal M}}
\newcommand{\cJ} {{\mathcal J}}
\newcommand{\cD} {{\mathcal D}}
\newcommand{\cC} {{\mathcal C}}

\def\bl{\begin{lemma}}
\def\el{\end{lemma}}
\def\bth{\begin{theorem}}
\def\eth{\end{theorem}}
\def\bc{\begin{corollary}}
\def\ec{\end{corollary}}
\def\bcj{\begin{conjecture}}
\def\ecj{\end{conjecture}}
\def\bpr{\begin{proposition}}
\def\epr{\end{proposition}}
\def\bde{\begin{definition}}
\def\ede{\end{definition}}
\newcommand{\be}{\begin{eqnarray}}
\newcommand{\ee}{\end{eqnarray}}

\def\Cl{\mathscr{C}}
\def\MA{\overline{M}}
\def\MO{M}
\def\muA{\overline{\mu}}
\def\muO{\mu}
\def\FO{\mathscr{F}}
\def\FA{\overline{\mathscr{F}}}
\def\IIC{\iic}
\def\1{1\!\!1}
\def\llra{\longleftrightarrow}
\def\lra{\leftrightarrow}
\def\given{\,|\,}
\def\bgiven{\;\big|\;}
\def\Bgiven{\;\Big|\;}
\def\size{\textsc{size}}
\def\vol{\textsc{vol}}
\def\radius{\textsc{radius}}
\def\help{\textsc{help}}

\def\spec{\mathsf{Spec}}
\def\Piv{\mathsf{Piv}}

\def\what{\widehat}
\def\fet{\mathsf{FET}}
\def\fetic{\mathsf{FETIC}}
\def\FETIC{\fetic}
\def\FET{\fet}
\def\iic{\mathsf{IIC}}
\def\build#1_#2^#3{\mathrel{ \mathop{\kern 0pt#1}\limits_{#2}^{#3}}}
\def\exc{\mathcal{E}}
\def\Exc{\mathcal{E}}

\def\coup{\mathbf{Q}}

\def\intg{\mathsf{Int}}
\def\thin{\mathsf{Fine}}
\def\thinning{\mathsf{Thinning}}
\def\good{\mathsf{Good}}
\def\vgood{\mathsf{VeryGood}}

\def\pthin{\P_{\rm thin}}
\def\pnorm{\P_{\rm norm}}

\def\th{{\rm thin}}
\def\no{{\rm norm}}

\def\annsmall{{A_{\rho(r),\,2\rho(r)}}}

\def\horiz{H}
\def\boxx{B}
\def\upp{U}
\def\low{L}
\def\clup{\mathcal{C}^{U_n}}
\def\cllow{\mathcal{C}^{L_n}}
\def\smallup{\mathcal{S}^{U_n}}
\def\smalllow{\mathcal{S}^{L_n}}

\def\tmall{\mathcal{T}}

\def\arrival{\mathcal{A}}

\def\thinrv{T}

\def\piic{\P_{\iic}}

\def\hexlatt{\mathcal{H}}
\def\orighex{0}
\def\Tg{\mathbb{T}}

\input epsf.sty
\input labelfig.tex

\begin{document}
\maketitle

\begin{abstract}
In dynamical critical site percolation on the triangular lattice or bond percolation on $\Z^2$, we define and study a local time measure on the exceptional times at which the origin is in an infinite cluster. We show that at a typical time with respect to this measure, the percolation configuration has the law of Kesten's Incipient Infinite Cluster. In the most technical result of this paper, we show that, on the other hand, at the first exceptional time, the law of the configuration is different. We also study the collapse of the infinite cluster near typical exceptional times, and establish a relation between static and dynamic exponents, analogous to Kesten's near-critical relation. 
\end{abstract}

\tableofcontents

\section{Introduction}\label{s.intro}

{\bf Critical planar percolation} is a central object of probability theory and statistical mechanics; see \cite{Grimmett,WWperc} for background.  The best understood example is $\mathsf{Bernoulli}(1/2)$ site percolation on the triangular lattice $\Tg$, where conformal invariance and hence convergence of interfaces to SLE$_6$ is known \cite{Smi, SchSLE, SmiICM, CNconv}. Nevertheless, many results are known for critical bond percolation on $\Z^2$ and other nice lattices, as well. In particular, almost everything in the present paper will apply equally to site percolation on~$\Tg$ and bond percolation on $\Z^2$.

In {\bf dynamical percolation}, a model introduced independently by \cite{HaPeSt} and Itai Benjamini, the status of each bit (site or bond) is continuously and independently resampled from the $\mathsf{Bernoulli}(p)$ measure, at times given by independent Poisson clocks of rate one. We will always consider site percolation on $\Tg$ and bond percolation on $\Z^2$, at the critical value $p=p_c=1/2$.
One of the principal reasons that dynamical percolation is interesting is that it provides a natural coupling of an uncountable number of copies of the underlying percolation process, and there may exist some exceptional instances of these copies
that satisfy certain events that have zero probability in static percolation. The existence (or non-existence) of such {\bf exceptional times} is called dynamical sensitivity (or stability) of the event, and the key event in question is of course the existence of an infinite cluster. See \cite{Steif} for a survey, but here is a brief summary of the subject. It was proved in \cite{HaPeSt} that for $p\not=p_c$ on any graph, both the existence and non-existence of infinite clusters are dynamically stable; then, dynamical stability of non-existence also holds at $p=p_c$ on $\Z^d$ with $d\geq 19$ and on regular trees; and finally, there exist non-regular but spherically symmetric trees with no infinite clusters at $p_c$ in static percolation, but with exceptional times in dynamical percolation. See \cite{Davar,PerSchSt} for more recent results on trees. The first example of dynamical sensitivity at $p_c$ in a transitive graph was given by \cite{SchSt}, proving it for the triangular lattice $\Tg$. This paper 
used discrete Fourier analysis, a tool that was introduced by \cite{BKS} for the closely related problem of noise sensitivity of percolation. This technique was further developed in \cite{GPS1}, proving that the set of exceptional times almost surely has Hausdorff dimension 31/36, and showing dynamical sensitivity of critical percolation also for bond percolation on $\Z^2$. Further studies of dynamical sensitivity and stability include \cite{BSplanes, exclusion, Ahlberg} for percolation type processes, \cite{BrSt} and \cite[Section 5]{DCGP} for Ising and random cluster Glauber dynamics, and \cite{BHPS,Hoff,web} for some other processes.

The rare appearances of infinite structure at the exceptional times are reminiscent of the {\bf Incipient Infinite Cluster}: a term used by physicists to refer to the large-scale connected structure present in critical percolation, and defined mathematically by Kesten as follows.  

\begin{definition}\label{d.IIC}
 The incipient infinite cluster, denoted by  $\IIC$, is the weak limit of the probability measures $\P_{p_c}( \, \cdot \, \big\vert 0 \leftrightarrow n )$ as $n \to \infty$, provided that the limit exists. 
\end{definition}

Here, $\{ 0 \lra n \}$ denotes the event that the open cluster of the origin reaches to distance~$n$. (We will formulate a precise definition shortly.)  
The existence of the $\IIC$ for numerous lattices in two dimensions was proved by Kesten~\cite{KestenIIC}. In high dimensions, properties of $\iic$ and its scaling limits have been investigated in detail using the lace expansion \cite{HaraSladeOne,HaraSladeTwo}. In two dimensions, several other natural means of locating large structures at criticality --- such as using the above definition with the condition $0 \leftrightarrow n$ replaced by the requirement that the open cluster of the origin have size at least $n$, or the weak limit as $n \to \infty$ of the largest cluster in $[-n,n]^2$ viewed from a uniformly chosen vertex in the cluster --- have been shown to also be equal to~$\iic$~\cite{Jar}. These results support the view that, at least in dimension two, any natural means of selecting a limit of large scale critical structure is the $\IIC$. One may ask then how the $\IIC$ may be found in dynamical percolation --- and this question is central to the present paper.

\subsection{The first exceptional time}

There is one very natural means of selecting an exceptional time at which the cluster of the origin in dynamical percolation is infinite:

\begin{definition}
Let $\mathcal{E}$ denote the random set of times at which the cluster of the origin is infinite.  
We define the {\it first exceptional time} $\fet$ to be $\inf \mathcal{E} \cap (0,\infty)$. That $\fet < \infty$ almost surely follows from the principal result of \cite{SchSt} for $\Tg$ and from \cite{GPS1} for $\Z^2$. Note that $\fet$ is positive almost surely, since some positive time passes before there is a change in any bit (be it site or bond) in the boundary of the finite cluster of the origin in $\omega_0$. The law of $\omega_{\fet}$ will be denoted by $\fetic$, the {\it first exceptional time infinite cluster}.    
\end{definition}

Although it may be a natural candidate for the appearance of the incipient infinite cluster in dynamical percolation, $\fetic$ is not the right choice: 

\begin{theorem}\label{t.FETIC}
The laws $\FETIC$ and $\IIC$ are not equal.
\end{theorem}

Proving Theorem \ref{t.FETIC} is this paper's most complex task. Roughly speaking, we show that the cluster of the origin under $\fetic$ is somewhat thinner than under $\iic$. Indeed, as we will state more precisely in the next subsection, while the configuration at a ``typical'' exceptional time turns out to have the law of $\iic$, with many other exceptional times nearby, $\FET$ appears at the endpoint of a unit-order interval in which exceptional times are absent; in fact, finite approximations to $\fetic$ may be constructed by size-biasing according to the length of the interval lacking connection from $0$ to a high distance $R$ leading up to a moment of such a connection. As such, $\fetic$ assigns more mass to configurations which are liable to break apart easily under the perturbation provided by dynamical percolation. What makes the proof difficult is to detect this imbalance also in the limit $R\to\infty$.
We will explain these vague ideas in more detail when we start proving Theorem~\ref{t.FETIC} in Section~\ref{s.FETIC}.

It seems natural to suppose that the two measures differ to a greater degree:  

\begin{conjecture}\label{c.sing}
The measures $\fetic$ and $\iic$ are singular with respect to each other.
\end{conjecture}

The above intuitive explanation about how biasing by the length of the waiting time makes $\fetic$ thinner than $\iic$ might suggest that $\iic$ stochastically dominates $\fetic$. However, $\iic$ does not satisfy the FKG inequality (which we shortly review), and so it may be that such a general conclusion does not follow from the negative conditioning represented by longer waiting times.

\begin{question}\label{q.stochdom}
Does $\iic$ stochastically dominate $\fetic$?
\end{question}

The invasion percolation cluster $\mathsf{IPC}$ is an infinite cluster associated to the critical point which is built by self-organized criticality. It was shown in \cite{DaSaVa} that $\IIC$ and $\mathsf{IPC}$ are singular with respect to each other on $\Z^2$. On the other hand, although $\IIC$ dominates $\mathsf{IPC}$ on regular trees \cite{invperctree}, this is not so on $\Z^2$ \cite{SapDoesnot}.

It was pointed out to us by Alain-Sol Sznitman that, instead of considering the distribution at the first entry to a given subset of the state space in a Markov process, which is $\fetic$ in our case, it is often more convenient to study the so-called equilibrium measure on the subset. For dynamical percolation on the ball $B_R$ and the subset $\cA:=\big\{\zeta \in  \{0,1\}^{B_R} : 0\lra R\text{ in }\zeta\big\}$, this measure is proportional at $\zeta \in \cA$ to the probability that dynamical percolation started at $\zeta$ and stopped at an independent exponential time $T$ leaves the set $\cA$ at the first update and does not return to it before $T$. The virtue of considering this measure could be that is has closer connections to the potential theory of the Markov process (Green's functions, Dirichlet forms, etc.; see \cite[Section 1.3]{Sznitman}) than the first entry time, hence it might be easier to address the analogues of Theorem~\ref{t.FETIC}, Conjecture~\ref{c.sing} and Question~\ref{q.stochdom} for 
this measure.

\subsection{The local time measure and the $\IIC$}\label{s.twodef}

Our first effort to seek the $\IIC$ in dynamical percolation was hampered by biasing created by the procedure for selection. In light of this, it is natural to try again by considering the law of the configuration obtained by selecting an exceptional time at a ``uniform'' moment. However, this notion of uniformity requires more structure on the exceptional time set in order to make sense. For this reason, and because of its intrinsic interest, we construct a local time measure $\mu$ on the exceptional time set $\exc$ as a weak limit of certain measures $\muO_r$ on the set of connection times to a large distance $r  \in \N$.

The simplest construction would be to define an approximative local time $\muA_r$ for distance $r \in \N$ by setting 
\begin{equation}\label{e.muA}
\MA_r(\omega_s):=\frac{\1\{0\leftrightarrow r\}}{\Ps{0\leftrightarrow r}},\qquad \textrm{and} \qquad 
\muA_r[a,b]:=\int_a^b \MA_r(\omega_s)\,ds \, ,
\end{equation}
and then hope that these measures have a limit $\muA[a,b]$ in some sense, as $r\to\infty$. However, we have encountered some  technical difficulties in trying to prove this convergence, hence will rely on the following slightly more complicated, but still very natural definition, which turns out to be easier to handle.

%
%

A local time is supposed to measure how much time the dynamical percolation process $\omega_s$ spends near $\Exc$. For this, we need some notion of how close a percolation configuration $\omega$ is to satisfying $0\lra\infty$. The simplest such notion was proposed in~(\ref{e.muA}): the existence of a connection to a large distance $r$. But we seem to get a more canonical notion by looking at how much a finite piece of the percolation configuration actually helps in realizing a connection to infinity. Namely, for any finite set $H$ of bits, we let $\omega^H$ denote the restriction of $\omega$ to $H$, and define the random variable
\begin{equation}\label{e.MO}
\MO_H(\omega):=\lim_{R\to\infty} \frac{\Ps{0\leftrightarrow R \given \omega^H}}{\Ps{0\leftrightarrow R}}.
\end{equation}
Of course, it is not at all obvious that the limit over $R$ exists. However,
\begin{equation}\label{e.MOIIC}
 \frac{\Ps{0\leftrightarrow R \given \omega^H}}{\Ps{0\leftrightarrow R}} =
 \frac{\Ps{0\leftrightarrow R,\  \omega^H}}{\Ps{0\leftrightarrow R}\Ps{\omega^H}}
=  \frac{\Ps{ \omega^H \given 0\leftrightarrow R}}{\Ps{\omega^H}},
\end{equation}
whose right-hand side indeed has a weak limit in high $R$ --- this is nothing other than the $\IIC$, whose construction was carried out in dimension two by Kesten~\cite{KestenIIC}. Thus, 
the limit in (\ref{e.MO}) indeed exists, so that we may define
\begin{equation}\label{e.muO}
\MO_r(\omega_s):=\MO_{B_r}(\omega_s),\qquad
\muO_r[a,b]:=\int_a^b \MO_r(\omega_s)\,ds \, .
\end{equation}
Note that $\E\MO_H(\omega_s)=1$ for any $H$, hence $\E\muO_r[a,b]=b-a$, independently of $r$, and we may hope to get a non-degenerate random measure in the limit $r\to\infty$. Moreover, and this is the main advantage of $\MO_r$ over $\MA_r$, the sequence $\{\muO_r[a,b]\}_{r\in\N}$ is a martingale with respect to the full filtration $\FO_r[a,b]$ generated by $\{\omega^{B_r}_s : s\in[a,b]\}$ (see~(\ref{e.muOMG}) in Section~\ref{s.localtime} for the proof). Thus, martingale convergence results can be used to prove the following:

\begin{theorem}\label{t.basic}
The limit $\muO[a,b]=\lim_{r\to\infty}\muO_r[a,b]$ of (\ref{e.muO}) exists almost surely and in $L^2$.

Assuming that the limit $\muA[a,b]=\lim_{r\to\infty}\muA_r[a,b]$ of (\ref{e.muA}) exists in $L^2$, the two local time measures almost surely coincide: $\muA[a,b]=\muO[a,b]$ for all intervals $[a,b]$ simultaneously.
\end{theorem}

So, we now have a measure from which we wish to sample uniformly to obtain a candidate for a law coinciding with $\iic$. However, $\muO$ is a $\sigma$-finite measure on $\R$ so that further work is needed to make valid the notion of sampling a uniform point in the measure. The next two theorems give constructions of such a point and show that indeed the law of the configuration at the selected time is $\iic$.

\begin{theorem}[Quenched sampling]
\label{t.IICq}
For almost every realization of the dynamical percolation process $\{\omega_s : s\in [0,\infty)\}$, and the corresponding local time measure $\muO$, there exists some $T_0<\infty$ such that for all $T>T_0$ we have $\muO[0,T]>0$. For such $T$, let $\chi_T$ be a random point from $[0,T]$ with law $\muO/\muO[0,T]$. Then, for almost all $\{\omega_s : s\in [0,\infty)\}$, the configuration $\omega(\chi_T)$ converges in law to $\IIC$, as $T\to\infty$.
\end{theorem}

\begin{theorem}[Annealed sampling]\label{t.IICa}\ 
\begin{itemize}
\item[{\bf (a)}] For any fixed $T>0$, let $\{\omega^*_s : s\in [0,T]\}$ be dynamical percolation reweighted (size-biased) by $\mu[0,T]$. Let $\chi^*_T$ be a random time from $[0,T]$ with law $\mu/\mu[0,T]$ for $\mu=\mu(\omega^*)$. Then, the configuration $\omega^*(\chi^*_T)$ has the distribution of the $\IIC$.
\item[{\bf (b)}] Given a sample of $\mu=\mu(\omega)$ on $\R$, let $\Pi_\mu$ be the Poisson point process with intensity $\mu$. 
One can make sense of conditioning $(\omega,\Pi_{\mu(\omega)})$ on $0\in \Pi_{\mu(\omega)}$; this is called $(\omega^*,\Pi_\mu^*)$, the Palm version of $(\omega,\Pi_\mu)$. Then $\omega^*_0$ has the law of the $\IIC$. 
\end{itemize}
\end{theorem} 

A concrete means of realizing the Palm version  of $(\omega,\Pi_\mu)$ from dynamical percolation $\omega$ is {Liggett's extra head construction}, which we will describe in Section~\ref{s.IIC}; see Figure~\ref{f.liggett}.
\medskip

Another application of the local time $\mu$ could be to run the dynamical percolation process $\omega$ according to $\mu(\omega)$. It should be possible to consider this time-changed dynamical percolation as a  Markov process on configurations satisfying $0\lra\infty$, with stationary measure $\IIC$; however, even the definition of the right state-space is unclear, especially if one wants $\IIC$ to be the unique stationary measure. We will not study these questions here.

\subsection{Structure of the paper}

In the rest of this Introduction, we summarize the necessary background in static and dynamical critical percolation. In Section~\ref{s.localtime}, we prove Theorem~\ref{t.basic}, and collect some 
properties of the finite and the limiting local time measures $\muA_r$, $\muO_r$, $\muO$. We then locate the $\IIC$ using the local time, proving Theorems~\ref{t.IICq} and~\ref{t.IICa} in Section~\ref{s.IIC}. The more substantial Section~\ref{s.FETIC} is devoted to telling apart $\fetic$ and $\iic$, with a thinning procedure on bounded configurations being introduced and analysed in order to prove Theorem~\ref{t.FETIC}. The proof of Theorem~\ref{t.FETIC} in fact exploits our identification of the $\iic$ in dynamical percolation, because the proof considers a uniform right-hand endpoint of a period of connection $0 \lra R$ and examines how long it takes for this connection to be reestablished as time advances; in finding an answer, we will exploit the fact that the law of the configuration in $B_R$ at this endpoint time is a close relative of critical percolation given $0 \lra R$ (and thus also of $\iic$). Section~\ref{s.collapse} contains Theorem~\ref{t.collapse}, a result addressing  the question of how instances of the $\iic$ embedded within dynamical percolation typically collapse as the time parameter is tuned at short distances to the moment at which the $\iic$ appears.

As mentioned above, all our results apply equally to critical site percolation on the triangular lattice $\Tg$ and critical bond percolation on $\Z^2$, except for the existence and values of some critical exponents, of course, but we will formulate our results without using these exponents. For the sake of definiteness, we will work with critical site percolation on $\Tg$, or rather, with critical percolation on the faces of the dual hexagonal lattice.

\subsection{Notation and percolation background}\label{ss.background}

Let $e_1$ and $e_2$ denote the Euclidean unit vectors. The lattice in $\R^2$ with generators $e_1$ and $\tfrac{1}{2} e_1 + \tfrac{\sqrt{3}}{2}e_2$
induces a Voronoi tiling of the plane whose faces are hexagons. We refer to the set of these hexagons, with the adjacency relation given by two hexagons sharing a common edge, as the {\bf hexagonal lattice} $\hexlatt$. The hexagon centred at the origin will be denoted by $\orighex$. 
Note that the set of hexagons intersecting the $x$-axis forms a bi-infinite simple path. Define $d: \hexlatt \times \hexlatt \longrightarrow \N$ to be graphical distance and set $B_R = \big\{ h \in \hexlatt: d( \orighex,h) \leq R \big\}$ for $R \in \N$. For $R_1,R_2 \in \N$ such that $R_1 < R_2$, write $A_{R_1,R_2} = B_{R_2} \setminus B_{R_1}$ for the annulus with inner and outer radii $R_1$ and $R_2$. The (outer) boundary of a set $A\subset \hexlatt$ is $\p A:=\{h\in\hexlatt \setminus A: d(h,A)=1\}$.

In critical percolation on $\hexlatt$, each $h \in \hexlatt$ is independently open or closed with probability one-half. The set $\{0,1\}^\hexlatt$ of percolation configurations is equipped with the usual product topology, and the events are the subsets $\cA\subseteq\{0,1\}^\hexlatt$ that are measurable with respect to the corresponding Borel sigma-algebra. For $a,b \in \hexlatt$, we write $a \lra b$ for the event that an open path of hexagons connects $a$ and $b$. For $A,B \subseteq \hexlatt$, we write $A \lra B$ if there exist $a \in A$ and $b \in B$ such that $a \lra b$. 
For $R_1,R_2 \in \N$ such that $1 \leq R_1 < R_2$, we write $R_1 \lra R_2$ to indicate that $\partial B_{R_1} \lra \partial B_{R_2}^c$. For $R \in \N$, we also write $0 \lra R$ for $\orighex \lra \partial B_R^c$.  

The open cluster of $\orighex$, $\big\{ h \in \hexlatt: \orighex \lra h \big\}$, will be denoted by $\Cl_\orighex$.

We will use the notation $\alpha_1(R_1,R_2):=\Ps{R_1 \lra R_2}$ and $\alpha_1(R):=\alpha_1(1,R)$, this being the {\bf one-arm} probability. Furthermore, $\alpha_4(R_1,R_2)$ denotes the {\bf alternating four-arm} probability: the probability that there are two open and two closed paths connecting $\p B_{R_1}$ and $\p B_{R_2}^c$, in an alternating order: open-closed-open-closed. Again, $\alpha_4(R):=\alpha_4(1,R)$. 

Given a percolation configuration $\omega\in\{0,1\}^\hexlatt$ and an event  $\cA\subseteq\{0,1\}^\hexlatt$, we call a hexagon $h$ {\bf pivotal} for $\cA$ in $\omega$ if changing the status of $h$ changes the outcome of the event. The set of pivotal hexagons will be denoted by $\Piv_\cA(\omega)$. For instance, note that $h$ is pivotal for the left-right crossing event in a rectangular region of $\hexlatt$ if and only if there are four alternating arms connecting $h$ to the corresponding sides of the rectangle. 

Let us now recall some standard tools in percolation theory \cite{WWperc}.

\medskip
\noindent{\bf The Harris-FKG inequality.} The set $\{0,1\}^\hexlatt$ of percolation configurations on the hexagonal lattice has a natural partial order $\leq$. A percolation event $\cA \subseteq \{0,1\}^\hexlatt$ is called increasing if $\omega \in \cA$ and $\omega \leq \omega'$ implies that $\omega' \in \cA$. The inequality of Harris and Fortuin-Kesteleyn-Ginibre states that if $\cA$ and $\cB$ are increasing events, then $\Ps{\cA \cap \cB} \geq \Ps{\cA} \Ps{\cB}$. 

\medskip
\noindent{\bf RSW estimates.} For any $L>0$, there exists a constant $c_L > 0$ such that the probability of an open path in critical percolation between the left and right sides of the region $\hexlatt \cap [0,Ln] \times [0,n]$ is at least $c_L$, independently of $n$.

\medskip
\noindent{\bf Quasi-muliplicativity of arm probabilities.} For $\ell\in\{1,4\}$, there exists a constant $0<c_\ell$ such that, for any radii $R_1<R_2<R_3$, we have 
\begin{equation}\label{e.qmulti}
c_\ell < \frac{\alpha_\ell(R_1,R_3)}{\alpha_\ell(R_1,R_2)\, \alpha_\ell(R_2,R_3)} \leq 1\,.
\end{equation}
The right-hand inequality is trivial; for $\ell=1$, the left-hand one is a simple consequence of FKG and RSW; for $\ell=4$, more work is needed, done in \cite{KestenScaling}; see also \cite{Nolin, SchSt}. Similarly to quasi-multiplicativity, one can show that we lose only a constant factor in probability if we require our four alternating arms to have their endpoints on nice prescribed arcs of the boundary. Together with some simple results on arm probabilities in the half-plane, this implies the following bounds on the number of pivotals: if $\cA(R)$ is the left-right crossing event in the square $[0,R]^2$, then $|\Piv_{\cA(R)}|\asymp \alpha_4(R)\,R^2$, and if $\cA(R_1,R_2)$ is the annulus crossing event $R_1\lra R_2$ with $R_1<R_2/2$, then $|\Piv_{\cA(R_1,R_2)}|\asymp \alpha_4(R_2)\,R_2^2$, with constant factors independent of $R_1$.

For critical percolation on $\hexlatt$, we also know the existence and values of critical exponents: $\alpha_1(R_1,R_2)=(R_1/R_2)^{5/48+o(1)}$ by \cite{LSW}, and  $\alpha_4(R_1,R_2)=(R_1/R_2)^{5/4+o(1)}$ by \cite{SW}, as $R_2/R_1\to\infty$. In particular, $|\Piv_{\cA(R)}|=R^{3/4+o(1)}$ as $R\to\infty$. On $\Z^2$, we have the bounds
\begin{equation}
\label{e.a4}
C^{-1}\, (r/R)^{2-\eta}\le \alpha_4(r,R) \le C\, (r/R)^{1+\eta}
\end{equation}
for some fixed constants $C > 0$, $\eta \in (0,1)$ and every $1\le r\le R$. See \cite[Appendix B]{SchSmG} and the references at \cite[Eq.~(2.6)]{GPS1}. Consequently, with some different value of the constant $C$,
\begin{equation}
\label{e.piv}
C^{-1}\, R^\eta \le |\Piv_{\cA(R)}| \le C\, R^{1-\eta}\,.
\end{equation}

\noindent{\bf The near-critical window.} One can consider monotone versions of dynamical percolation, in which dynamical updates lead always either to the closure or to the opening of hexagons. These give couplings between dynamical and off-critical percolation (and also a coupling of percolation measures at different densities), and therefore information on off-critical percolation can yield bounds on dynamical percolation questions. We will use these relations (which turn out to be sharp) several times.

Kesten found the near-critical window of percolation precisely \cite{KestenScaling} (see \cite{Nolin,WWperc} for more modern accounts): 
for a system of linear size $R$, the window is given by the reciprocal of the expected number of pivotals for the left-right crossing event $\cA(R)$ at criticality. More precisely, for the annulus crossing event $\cA(R,2R)$, as $R\to\infty$, we have
\begin{equation}\label{e.winsmall}
\frac{\Pso{p_c\pm\eps}{\cA(R,2R)}}{\Pso{p_c}{\cA(R,2R)}}\to 1 \qquad \text{if}\quad \eps \ll  \frac{1}{|\Piv_{\cA(R)}|}\,, 
\end{equation}
while 
\begin{equation}\label{e.winmed}
\delta  < \Pso{p_c\pm\eps}{\cA(R,2R)} < 1-\delta \qquad \text{if}\quad \eps \asymp  \frac{1}{|\Piv_{\cA(R)}|}\,, 
\end{equation}
with $\delta\in(0,1)$ depending only on the constant factors giving the size of $\eps$, and finally,
\begin{equation}\label{e.winlarge}
\Pso{p_c+\eps}{\cA(R,2R)} \to \begin{cases} 	1 & \qquad \text{if}\quad \eps \gg \frac{1}{|\Piv_{\cA(R)}|}\,,\\
									0 & \qquad \text{if}\quad -\eps \gg \frac{1}{|\Piv_{\cA(R)}|}\,.
									\end{cases}  
\end{equation}
Kesten also proved the stability of one- and alternating four-arm probabilities inside the window:
\begin{equation}\label{e.armstab}
\begin{aligned}
\frac{\Pso{p_c\pm\eps}{\cA_\ell(1,R)}}{\Pso{p_c}{\cA_\ell(1,R)}}
& \to 1 \qquad \text{if}\quad \eps \ll \frac{1}{|\Piv_{\cA(R)}|}\,,\\
& \asymp\ 1 \qquad \text{if}\quad \eps \asymp \frac{1}{|\Piv_{\cA(R)}|}\,,
\end{aligned}								
\end{equation}
for $\ell\in\{1,4\}$. The $\eps \ll 1/|\Piv_{\cA(R)}|$ case of (\ref{e.armstab}) and (\ref{e.winsmall}) are not stated explicitly in \cite{KestenScaling}, but they clearly follow from his proof using differential inequalities.

Using the stability of the 1-arm and 4-arm probabilities in the near-critical window, he also found the off-critical exponent, a relation usually called Kesten's scaling relation \cite[Corollary 1]{KestenScaling}: 
\begin{equation}\label{e.offcrit}
\Pso{p_c+\eps}{0\llra\infty} \asymp \alpha_1(\rho(1/\eps))\,,
\end{equation}
where $\rho(r):=\inf\{ s\in\N_+ : s^2\,\alpha_4(s)\ge r\}$ for $r\geq 1$,  the inverse function of $R\mapsto |\Piv_{\cA(R)}|$. We have $\rho(r)=r^{4/3+o(1)}$ on $\hexlatt$, and $C^{-1}r^\eta \leq \rho(r) \leq C r^{1/\eta}$ for some $0<\eta,C<\infty$ on $\Z^2$, by~(\ref{e.piv}). Note here that Kesten formulated his result in terms of critical exponents, which would not be enough for us later because of the unspecified $o(1)$ terms in the exponent, but the proof clearly gives the stronger result we stated; see \cite[Chapter 6]{WWperc}.
\medskip

\noindent{\bf Dynamical percolation and a dynamical FKG inequality.} As mentioned above, we will consider dynamical critical percolation with updates from the stationary distribution (resampling the bits) at times given by Poisson clocks of rate one, with time indexed by $\R$, and, just for the sake of definiteness, with c\`adl\`ag trajectories.

We will need the following extension of the FKG inequality to increasing events of dynamical percolation, an immediate consequence of \cite[Corollary II.2.21]{LiggettBook}. A weaker form (with a very different proof) was  given in \cite[Lemma 4.2]{HMP}.

\begin{lemma}[Dynamical FKG inequality]\label{lemdynfkg}
Let $\omega,\omega': \hexlatt \times \R\longrightarrow \{ 0,1 \}$ denote two realizations of dynamical percolation on the hexagonal lattice $\hexlatt$. We say that $\omega \leq \omega'$ if $\omega_t(x) \leq \omega'_t(x)$ for all $(x,t) \in  \hexlatt \times \R$. Let $\cA,\cB \subseteq \{0,1\}^{ \hexlatt \times \R}$ be two increasing events (i.e., if $\omega\in \cA$ and $\omega\leq \omega'$, then $\omega'\in \cA$). Then $\P(\cA \cap \cB) \geq \P(\cA) \P(\cB)$.

The same holds if the dynamics is not stationary, but started at time $0$ from an arbitrary distribution on $\{0,1\}^\hexlatt$ that satisfies the static FKG inequality.
\end{lemma}

\proof Corollary II.2.21 of \cite{LiggettBook} states this for increasing events that depend on the configuration at finitely many time instances $t_1<\dots<t_n$, proved using induction on $n$ and the infinitesimal generator of the process. Since all measurable dynamical events can be approximated by events depending on finitely many time instances, our statement follows.\qed

%

\subsection{The Fourier spectrum of critical percolation}

A key tool for the analysis of dynamical percolation is discrete Fourier analysis. Here we provide the definition of the Fourier spectrum of a percolation event, explain the basic relation between the spectrum and decorrelation for the event under dynamical percolation, and collect the results from the literature that we will use. A far more thorough overview of this theory is provided by the survey article~\cite{GS}.

Let $\cA$ denote a percolation event in $B_R$, so that $\cA$ is a subset of percolation configurations in $B_R$. Define the usual inner product on the $L^2$-space on percolation configurations on $B_R$ by $\langle f,g \rangle = \Es{fg}=2^{- | B_R |} \sum_{\omega \in \{-1,1\}^{B_R}} f(\omega)g(\omega)$, and note that the collection $\big\{ \chi_S:=\prod_{i\in S}\omega(i) : S \subseteq B_R \big\}$ is an orthonormal basis for this $L^2$-space. As such, the $\{-1,1\}$-indicator function $f_\cA$ of $\cA$ has a Fourier decomposition $f_\cA=\sum_{S \subseteq B_R}\what f_\cA(S) \chi_S$. Parseval's identity $\sum_{S \subseteq B_R} \what{f_\cA}^2(S) = 1$ allows us to define a random variable $\spec_A$, the spectral sample of $\cA$, on subsets of $B_R$ according to $\Ps{\spec_A = C} = \what{f_\cA}^2 (C)$ for $C \subseteq B_R$.

Recall that the dynamical percolation process $\{\omega_t\}_{t\in\R}$ is defined using i.i.d.~rate one Poissonian updates for each bit. 
Now, the basic relation between the spectral sample and decorrelation under this dynamics is that, for percolation events $\cA$ and $\cB$ in $B_R$,
\begin{equation}\label{e.fourierrelation}
 \Es{\omega_0 \in \cA, \omega_t \in \cB} = \sum_{S \subseteq B_R} \what{f}_\cA(S) \what{f}_\cB(S) e^{-t \vert S \vert} \, .
\end{equation}
This shows that if most of the measure for at least one of the spectral samples $\spec_\cA$, $\spec_\cB$ is supported on large sets $S$, then fast decorrelation occurs.

The spectral sample $\spec_\cA$ is a random subset of $B_R$ with some similarities to, and marked differences from, the random set $\Piv_\cA$ of hexagons in $B_R$ that are pivotal for the occurrence of $\cA$ under critical percolation. As first observed by Gil Kalai, the two random variables share their first and second moments (see \cite[Section 2.3]{GPS1}),
\begin{equation}\label{e.specmoment}
 \E \vert \Piv_\cA \vert  = \E \vert \spec_\cA \vert  \, , \qquad \E \vert \Piv_\cA \vert^2  = \E \vert \spec_\cA \vert^2  \,,
\end{equation}
but not the higher ones, and their large deviations usually differ (see \cite[Remark 4.6]{GPS1}).

Of particular import to us is the case where $\cA$ is a crossing event from one boundary arc to another in some planar domain. Let us first consider $\cA(R,2R) = \big\{R \lra 2R \}$. A standard second moment argument yields the conclusion that there exists $C > 0$ such that, for all $R$, $\E { \vert \Piv_{\cA(R,2R)} \vert^2 } \leq C\, \E { \vert \Piv_{\cA(R,2R)} \vert }^2$. In light of~(\ref{e.specmoment}) and the second moment method, we see that for some small $c > 0$, 
$$\Pb{ \vert \spec_{\cA(R,2R)} \vert \geq c \, \E \vert \spec_{\cA(R,2R)}\vert } \geq c\,.$$
Thus,~(\ref{e.specmoment}) shows that, for each $s > 0$, there exists $c(s) < 1$ (with the supremum of $c(s)$ strictly less than one over any interval of the form $(\eps,\infty)$) such that, for all $R \in \N$,
\begin{equation}\label{e.partialdecor}
 \Pb{\omega_0 \in {\cA(R,2R)}, \omega_{t} \in {\cA(R,2R)}} \leq c(s)  \quad \textrm{where }t = s\, \E\vert \Piv_{\cA(R,2R)} \vert  \, ;
\end{equation}
thus, the characteristic time-scale for at least partial decorrelation of the crossing event is determined by the mean number of pivotals. We will also need the much stronger assertion, proved in \cite{GPS1}, that, as $s\to\infty$, 
\begin{equation}\label{e.fulldecor}
 \Pb{\omega_0 \in {\cA(R,2R)}, \omega_{t} \in {\cA(R,2R)}} \to 0 \quad \textrm{where }t = s\, \E \vert \Piv_{\cA(R,2R)} \vert  \, ,
\end{equation}
uniformly in $R \in \N$; on $\hexlatt$ we have the sharp upper bound $s^{-2/3+o(1)}$. That is, the crossing event in fact decorrelates fully at large multiples of the scale determined by the mean pivotal number. The bound (\ref{e.fulldecor}) arises from a detailed examination of the lower-tail of the size $\vert \spec_{\cA(R,2R)} \vert$ of the spectral sample.

Similar sharp results are proved in \cite{GPS1}  for the decorrelation of the crossing events $\cA(0,R)=\{0\lra R\}$, which are the key for the applications to exceptional times. Namely, \cite[Equation~(9.2)]{GPS1} says that, for all $s,t\geq 0$,
\begin{align}
\frac{\PB{\1\{0\mathop{\llra}\limits^{\omega_s} R\}
\1\{0\mathop{\llra}\limits^{\omega_t} R\}}}{\Ps{0\llra R}^2}
& \leq O(1) \frac{1}{\alpha_1(\rho(1/|t-s|))} \label{e.radialrho}\\
& \leq  O(1) \, |s-t|^{-1+\delta+o(1)} \label{e.radialdelta}
\end{align}
for some $\delta>0$, uniformly in $R\in\N_+$, the $o(1)$ term being understood as $|s-t|\to 0$. On $\hexlatt$, also the sharp result $\delta=31/36$ is known. For exceptional times, the importance of these decorrelation bounds lies in the fact that the exponent $\delta$ of~(\ref{e.radialdelta}) is a lower bound on the Hausdorff dimension of the set $\exc$, using the so-called Mass Distribution Principle.

\bigskip
\noindent{\bf Acknowledgments.} We thank Yuval Peres and Alain-Sol Sznitman  for useful discussions, and Jeff Steif for pointing out an error in an earlier version.

Parts of this work were done at the Theory Group of Microsoft Research, Redmond, at New York University, at the University of Toronto, and at the Fields Institute in Toronto. AH was supported by NSF grants DMS-0806180 and OISE-0730136 at New York University, and by EPSRC grant EP/I004378/1 at the University of Oxford. GP was supported by an NSERC Discovery Grant at the University of Toronto, and an EU Marie Curie International Incoming Fellowship at the Technical University of Budapest.

\section{Construction and basic properties of the local time}\label{s.localtime}

In this section, we present the proof of Theorem~\ref{t.basic}, and collect some basic and less basic properties of the finite and the limiting local time measures. We begin by examining the martingale property for the approximating local time measures $\muA_r[a,b]$ and $\muO_r[a,b]$, defined in~(\ref{e.muA}) and~(\ref{e.muO}).

Note that $\MA_R(\omega)$ is a martingale with respect to the filtration $\FA_R$ of the percolation space generated by the variables $\big\{ \1\{0\leftrightarrow r\} : r\leq R\big\}$; indeed, for any $r'>r$,
$$
\E\left(\frac{\1\{0\leftrightarrow r'\}}{\Ps{0\leftrightarrow r'}} \, \bigg\vert \, \FA_r\right)
=\frac{\Ps{0\leftrightarrow r' \given 0\leftrightarrow r}}{\Ps{0\leftrightarrow r'}}\1\{0\leftrightarrow r\}
=\frac{\1\{0\leftrightarrow r\}}{\Ps{0\leftrightarrow r}}.
$$
Similarly, it is clear from~(\ref{e.MO}) that $\MO_r(\omega)$ is a martingale with respect to the full filtration $\FO_r$ generated by $\omega^{B_r}$. Being a martingale w.r.t.~this larger sigma-algebra is more useful: 
\begin{equation}\label{e.muOMG}
\begin{aligned}
\EB{\muO_{R}[a,b] \Bgiven \FO_r[a,b]} &= 
\int_a^b \EB{\MO_{R}(\omega_s)\Bgiven \FO_r[a,b]}\, ds\\
&= \int_a^b \EB{\MO_{R}(\omega_s)\Bgiven \FO_r(\omega_s)}\, ds = 
\int_a^b \MO_r(\omega_s) \, ds = \muO_r[a,b]\,;
\end{aligned}
\end{equation}
that is,  $\muO_r[a,b]$ is a martingale w.r.t.~$\FO_r[a,b]$. On the other hand, $\muA_r[a,b]$ does not seem to be a martingale w.r.t.~$\FA_r[a,b]$, since 
$$\EB{\MA_R(\omega_s) \md \FA_r[a,b]} \not= \EB{\MA_R(\omega_s) \md \FA_r(s)}$$
in general, because of the extra information provided by $\FA_r(t)$, $t\in [a,b] \setminus\{s\}$. 

Consequently, it is much simpler to prove the convergence of $\muO_r$ to some limit $\muO$ than the convergence of $\muA_r$, though we expect that the latter also holds: as we will see in the forthcoming proof, the local time densities $\MA_r$ and $\MO_r$ are closely related to each other.

\begin{figure}[htbp]
\SetLabels
(1*.8)$\MA_R(\omega_s)$\\
(1*.4)$\MA_r(\omega_s)$\\
(.75*.15)\textcolor{red} {$\MO_r(\omega_s)$}\\
(1.02*.07)time\\
\endSetLabels
\centerline{
\AffixLabels{
\epsfysize=1.5in \epsffile{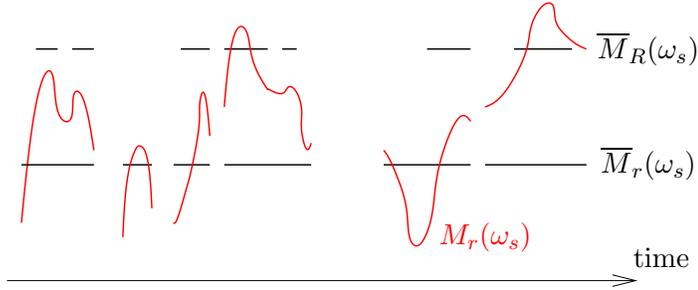}
}
}
\caption{Schematic pictures of the approximate local time densities for $\muA_r$ and $\muO_r$.}
\label{f.muOmuA}
\end{figure}

\proofof{Theorem~\ref{t.basic}}
The main task is to prove the statements for any fixed interval $[a,b]$. This implies the claims for all $(a,b) \in \Q^2$ simultaneously, and then, since rational intervals generate the Borel sigma-algebra,  
we have the statement simultaneously for all $(a,b) \in\R^2$.

First recall the quasi-multiplicativity relation (\ref{e.qmulti}), which implies, for $R>r>0$,
$$
\frac{\Ps{0\leftrightarrow R \given \omega^{B_r}}}{\Ps{0\leftrightarrow R}} \asymp \frac{\Ps{0\leftrightarrow R \given \omega^{B_r}}}{\Ps{0\leftrightarrow r}\Ps{r\leftrightarrow R}} \leq \frac{\Ps{r\leftrightarrow R \given \omega^{B_r}}\1\{0\leftrightarrow r\}}{\Ps{0\leftrightarrow r}\Ps{r\leftrightarrow R}} = \frac{\1\{0\leftrightarrow r\}}{\Ps{0\leftrightarrow r}}.
$$
Therefore, with an absolute constant $C_1<\infty$,
\begin{equation}\label{e.MOMA}
\MO_r(\omega) \leq C_1 \MA_r(\omega)\quad\text{and}\quad \muO_r[a,b] \leq C_1 \muA_r[a,b].
\end{equation}
Second, recall from (\ref{e.radialdelta}) the bound $O(1)|s-t|^{-1+\delta+o(1)}$, with $\delta>0$.  
Integrating over $s$ and $t$, this gives the second moment estimate
\begin{equation}
\label{e.mu2nd}
\Eb{\muA_r[a,b]^2} \leq 
\begin{cases} 
|b-a|^{1+\delta+o(1)} & \text{ as }|b-a|\to 0\,,\\
C_2 \, |b-a| & \text{ for all } a<b\,, 
\end{cases}
\end{equation}
uniformly in $r$, with an absolute constant $C_2<\infty$. Therefore, by (\ref{e.MOMA}), the sequence $\muO_r[a,b]$ is an $L^2$-bounded martingale w.r.t.~$\FO_r[a,b]$, and the $L^2$ martingale convergence theorem implies the existence of the limit
\begin{equation}\label{e.MGlimit}
\muO_r[a,b]\xrightarrow[L^2]{\;\text{a.s.}\;} \muO[a,b]\,.
\end{equation}

Now, we turn to the sequence $\muA_r[a,b]$. If we fix $r>0$, and take $R\to\infty$, then
$$
\Eb{\MA_{R}\bgiven \FO_r}=\frac{\Ps{0\lra R\given \FO_r}}{\Ps{0\lra R}} \xrightarrow[L^\infty]{\;\text{a.s.}\;} \MO_r\,,
$$
by the very definition of $\MO_r$, a random variable on the finite space $B_r$. Thus, for fixed $r$, the random variables $\Eb{\MA_{R} \bgiven \FO_r}$ are uniformly bounded in $R$, and
$$
\int_a^b \EB{\MA_{R}(\omega_s)\Bgiven \FO_r(\omega_s)}\, ds \xrightarrow[L^\infty]{\;\text{a.s.}\;} \int_a^b \MO_r(\omega_s) \, ds = \muO_r[a,b].
$$
On the other hand, for random variables, convergence in $L^2$ is stronger than in $L^1$, hence the hypothetical $L^2$-convergence of the unconditional $\muA_{R}[a,b]$ implies
\begin{align*}
\int_a^b \EB{\MA_{R}(\omega_s)\Bgiven \FO_r(\omega_s)}\, ds &= 
\int_a^b \EB{\MA_{R}(\omega_s)\Bgiven \FO_r[a,b]}\, ds \\
&= \EB{\muA_{R}[a,b] \Bgiven \FO_r[a,b]} \xrightarrow[L^2]{} \EB{\muA[a,b] \Bgiven \FO_r[a,b]}.
\end{align*}
One sequence can have only one $L^2$-limit, and convergence in $L^\infty$ is stronger than in $L^2$, thus
\begin{equation}\label{e.MOMAcond}
\EB{\muA[a,b] \Bgiven \FO_r[a,b]}= \muO_r[a,b]\quad\text{ in }L^2,\text{ hence almost surely}.
\end{equation}
As $r\to\infty$, $\FO_r[a,b]$ converges to the full sigma-algebra, hence the left-hand side of (\ref{e.MOMAcond}) converges a.s.~to $\muA[a,b]$ by L\'evy's zero-one law, while the right-hand side converges to $\muO[a,b]$, by (\ref{e.MGlimit}). The two limits coincide a.s., thus the proof of Theorem~\ref{t.basic} is complete.\qed
\medskip

\begin{conjecture}\label{c.muA}
The $L^2$-limit $\muA_[a,b]=\lim_{r\to\infty} \muA_r[a,b]$ exists, and then, by Theorem~\ref{t.basic}, $\muA=\muO$ almost surely.
\end{conjecture}

We collect now some basic properties of the dynamical percolation process, the exceptional set, and the associated local time. 

\begin{lemma}[Ergodicity]\label{l.erg}
The dynamical percolation process $\omega$ on the infinite lattice (in particular, the local time $\mu=\mu(\omega)$) is ergodic with respect to~time shifts.
\end{lemma}

\proof For any dynamic event $\cA$ and any $\eps>0$, there exists a radius $r\in\N$, a time $T>0$, and an event $\cA_{r,T}$ measurable with respect to~$\omega^{B_r}(-T,T)$ such that $\Pb{\cA \,\triangle\, \cA_{r,T}}<\eps$. Now, by the ergodicity of dynamical percolation in $B_r$ (a Markov chain on a finite state space), there exists $t=t(r,T)$ such that $\big|\Ps{\cA_{r,T} \cap (\cA_{r,T}+t)}-\Ps{\cA_{r,T}}^2\big|<\eps$, where $\cA_{r,T}+t$ represents the event $\cA_{r,T}$ evaluated for the dynamical configuration shifted back by time $t$. Now, if $\cA$ is invariant under time shifts, then $\big|\Ps{\cA \cap (\cA+t)} - \Ps{\cA_{r,T} \cap (\cA_{r,T}+t)}\big| \leq \Pb{\cA \, \triangle\, (\cA_{r,T} \cap (\cA_{r,T}+t))} < 2\eps$. Altogether, $\big|\Ps{\cA}-\Ps{\cA}^2\big| < 2\eps + \eps + \eps^2$. This holds for any $\eps>0$, hence $\Ps{\cA}\in\{0,1\}$.
\qed

\begin{lemma}[Perfectness]\label{l.exctop} Almost surely, the set $\Exc$ of exceptional times\\
{\bf (i)} is disjoint from the set of times at which the status of a hexagon is updated;\\
{\bf (ii)} is topologically closed;\\
{\bf (iii)} has no isolated points. 
\end{lemma}

\proof Parts (i) and (ii) are proved in \cite[Lemma 3.2]{HaPeSt}. Part (iii) is proved in \cite[Lemma 3.4]{HaPeSt} and the remark following it.\qed

\begin{lemma}[No atoms]\label{l.noatom} There are almost surely no atoms in any of the measures $\muA_r$, $\muO_r$ or $\muO$.
\end{lemma}

\proof Fix any large $n\in\N$, and cover the interval $[0,1]$ by the intervals $I_i:=[\tfrac{i}{2n},\tfrac{i}{2n}+\tfrac{1}{n}]$, $i=0,1,\dots,2n-2$. By (\ref{e.mu2nd}) and Chebyshev's inequality, for any $c>0$ and any index $i$, we have $\Ps{\muA_r(I_i) > c} \leq n^{-1-\delta+o(1)}$ as $n\to\infty$, uniformly in $r$. By a union bound, $\Ps{\exists i : \muA_r(I_i)>c} \leq n^{-\delta+o(1)} \to 0$, which implies the claim for $\muA_r$. Then (\ref{e.MOMA}) implies it for $\muO_r$, and (\ref{e.MGlimit}) implies it for $\muO$. 
\qed
\medskip

It is clear that $\muO$ is supported inside $\exc$. The following statement is very natural, but it seems hard to prove:

\begin{conjecture}
\label{c.support}
The support of the local time measure $\mu$ is almost surely the entire exceptional time set $\Exc$.
\end{conjecture}

This conjecture cannot fail by much: the Hausdorff dimension of $\supp\mu$ is the same as the dimension of $\Exc$, namely $31/36$. The reason is that the proof of the lower bound in \cite{GPS1} (just like in \cite{SchSt}) uses the approximate local time measures $\muA_r$ and a version of the Mass Distribution Principle, and, via~(\ref{e.MOMA}), it could also have used the measures $\muO_r$, hence it actually yields a lower bound on $\dim_H(\supp\mu)$. The next lemma, which will be of use later, provides a little further evidence for the conjecture.

\begin{lemma}\label{l.cond2ndM}
For any $\eps>0$, let $\mu_\eps$ denote $\muA_r[0,\eps]$ or $\muO_r[0,\eps]$ or the limit
$\muO[0,\eps]$. Then there is an absolute constant $C<\infty$ such that
\begin{equation}\label{e.musecmom}
\Eb{\mu_\eps^2 \md \mu_\eps>0} \leq C\, \Eb{\mu_\eps \md \mu_\eps>0}^2\, ,
\end{equation}
and another such constant $c>0$ such that
\begin{equation}\label{e.cc}
\Pb{\mu_\eps > c\,\E\mu_\eps \md \Exc\cap[0,\eps]\not=\emptyset} > c\,.
\end{equation}
\end{lemma}

\proof 
The left-hand side of~(\ref{e.musecmom}) equals $\Es{\mu_\eps^2}/\Ps{\mu_\eps>0}$, while the right-hand side equals
$\Es{\mu_\eps}^2/\Ps{\mu_\eps>0}^2=\eps^2/\Ps{\mu_\eps>0}^2$. Hence, we
need to show that
$$
\Eb{\mu_\eps^2} \leq C\, \frac{\eps^2}{\Ps{\mu_\eps>0}}\,.
$$
By a usual coupling between dynamical and near-critical percolation, in which dynamical updates lead always to the opening of hexagons in the latter case, we have 
\begin{equation}\label{e.nceup}
\begin{aligned}
\Ps{\mu_\eps>0}\leq \Ps{\Exc\cap[0,\eps]\not=\emptyset} &\leq \Pso{p_c+O(\eps)}{0\llra\infty}\\
&=O(1)\,\alpha_1(\rho(1/\eps))\,,
\end{aligned}
\end{equation} 
by Kesten's scaling relation (\ref{e.offcrit}). On the
other hand, taking the double integral of~(\ref{e.radialrho}) over $s,t\in[0,\eps]$, we claim that
\begin{equation}\label{e.muepstwo}
\Eb{\mu_\eps^2} \leq O(1)\, \frac{\eps^2}{\alpha_1(\rho(1/\eps))}\,,
\end{equation}
which will finish the proof of (\ref{e.musecmom}). 

By (\ref{e.MOMA}) and (\ref{e.MGlimit}), it is enough to verify~(\ref{e.muepstwo}) for $\mu_\eps=\muA_r[0,\eps]$.
Set $R = \rho(1/\eps)$ and $A_i = [C^i R, C^{i+1}R]$, $i 
\in \N$, where $C > 0$ is a large constant to be specified shortly. For $i \in \N$, write
$$
B_i = \Big\{ (s,t) \in [0,\eps]^2 : \rho \big( \vert s - t \vert^{-1} \big) \in A_i  \Big\} \, ,
$$
so that
$$
\Eb{\mu_\eps^2}  =   \int_{[0,\eps]^2}
\frac{\PB{\1\{0\mathop{\llra}\limits^{\omega_s} r\}
\1\{0\mathop{\llra}\limits^{\omega_t} r\}}}{\Ps{0\llra r}^2} \, {\rm d} s \, {\rm d} t 
 \leq   O(1) \sum_{i \geq 0} \phi_i \, ,
$$
with $\phi_i = \int_{B_i} \frac{1}{\alpha_1(\rho(1/|t-s|))} \, {\rm d} s \, {\rm d} t$;
the latter inequality is due to (\ref{e.radialrho}).

Note that $\rho(\cdot)$ is a non-strictly increasing function. By~(\ref{e.a4}), there exists an absolute constant $K>0$ such that $s^2\alpha_4(s) < (Ks)^2 \alpha_4(Ks)$ for all $s\in\Z^+$, hence $\rho(s^2\alpha_4(s))\in (s/K, s]$ for all $s\in\Z^+$.
This implies that
\begin{equation}\label{e.K}
\rho^{-1}(A_i)\subseteq \big[ (C^i R)^2 \alpha_4(C^iR),\, (C^{i+1}KR)^2\alpha_4(C^{i+1}KR)\big] \, .
\end{equation}

If $C$ is large enough, so that $(CR)^2\alpha_4(CR) >  2(KR)^2\alpha_4(KR)$, then $(1/\eps,2/\eps)\subseteq \rho^{-1}(A_0)$, hence the Lebesgue measure of $B_0$ is at least $\eps^2/2$. Therefore, 
$$
\phi_0 =  \int_{B_0} \frac{1}{\alpha_1(\rho(1/|t-s|))} \,  {\rm d} s \, {\rm d} t \geq \frac{\eps^2}{2} \, \alpha_1 \big( \rho(1/\eps) \big)^{-1}\,.
$$
On the other hand, for $i \geq 1$, using (\ref{e.K}), 
$$
\phi_i =  \int_{B_i} \frac{1}{\alpha_1(\rho(1/|t-s|))} \,  {\rm d} s \, {\rm d} t
\leq 2 \eps C^{-2i} R^{-2} \alpha_4 \big( C^i R \big)^{-1}  \alpha_1 \big( C^{i+1} R \big)^{-1} \, .
$$
Thus,
\begin{eqnarray*}
 \frac{\phi_i}{\phi_0} & \leq &  4 \eps^{-1} C^{-2i} R^{-2} \frac{\alpha_1(R)}{\alpha_4(C^i R) \alpha_1(C^{i+1}R)} \\
  & \leq & \frac{4 C^{-2i}}{\alpha_4(R,C^i R) \alpha_1(R,C^{i+1}R)}\,,
\end{eqnarray*}
where in the second inequality we used that $\eps^{-1}\leq R^{2}\alpha_4(R)$. Now, \cite[Appendix]{GPS1} says that the sum of the 1-arm and 4-arm exponents is strictly less than 2 --- properly interpreted in the case of $\Z^2$ where these exponents are not known to exist. That is, there exists some $c \in (0,1)$ such that $\phi_i/\phi_0 \leq O(1) \, c^i$ for all $i \geq 1$. Thus, 
$$\int_{[0,\eps]^2}  \frac{1}{\alpha_1(\rho(1/|t-s|))} \,  {\rm d} s \, {\rm d} t \leq O(1)\, \int_{A_0}  \frac{1}{\alpha_1(\rho(1/|t-s|))} \,  {\rm d} s \, {\rm d} t \leq O(1)\, \eps^2 \alpha_1 \big( \rho(1/\eps) \big)^{-1}\,,$$
and we have confirmed~(\ref{e.muepstwo}). 

By the Paley-Zygmund second moment inequality (a simple consequence of Cauchy-Schwarz; see, e.g., \cite[Section 5.5]{LPbook}), the above computations show that 
$$
\Ps{\mu_\eps>0}\geq \frac{(\E \mu_\eps)^2}{\Es{\mu_\eps^2}} \geq c_1\, \alpha_1(\rho(1/\eps))\,,
$$
matching the upper bound (\ref{e.nceup}) up to a constant factor. Therefore,
$$
\Pb{\mu_\eps>0 \md \Exc\cap[0,\eps]\not=\emptyset} > c_2 > 0\,.
$$
On the other hand, again by the Paley-Zygmund inequality, (\ref{e.musecmom}) implies that 
$$
\Pb{\mu_\eps > c_3\,\E\mu_\eps \md \mu_\eps>0 }  > c_3 > 0\,,
$$
for some $c_3>0$. Combining the last two displayed inequalities proves (\ref{e.cc}).
\qed
\medskip

We conclude this section with a natural question:

\begin{question}
Is the local time $\mu$ the 31/36-dimensional Minkowski content of the set $\exc$? 
Is $\mu$ the Hausdorff measure of $\exc$ for some Hausdorff gauge function? 
\end{question}

\section{Finding the Incipient Infinite Cluster}\label{s.IIC}

Given the description of the local time measure using~(\ref{e.MOIIC}), it is natural to guess that the infinite cluster at a ``typical'' exceptional time (typical with respect to~$\mu$) has the law of $\IIC$. The first exceptional time having been discredited as a candidate for the $\IIC$ by Theorem~\ref{t.FETIC}, we now prove Theorems~\ref{t.IICq} and \ref{t.IICa}, thereby verifying what may be the simplest relationship between exceptional times and the $\IIC$.

Unsurprisingly, the proofs go through the finite approximations, about which we provide a further definition.
\begin{definition}
Let $\iic_r$ denote the law on percolation configurations in $B_r$
given by $\Pbo{p_c}{ \, \cdot \md \orighex \leftrightarrow r }$. 
\end{definition}

Note that $\MA_r(\omega)$ is the Radon-Nikodym derivative ${\d\IIC_r}/{\d\P}$, while $\MO_r(\omega)$ is the Radon-Nikodym derivative ${\d\IIC^{B_r}}/{\d\P}$, where $\P=\P_{p_c}$ is critical percolation. 
Since both $\IIC_r$ and $\IIC^{B_r}$ converge to $\IIC$ as $r\to\infty$, both $\muA_r$ and $\muO_r$ can be useful in studying the relationship between dynamical percolation and the $\iic$. Indeed, in the forthcoming lemmata, the versions about $\muO_r$ will be used in finding the $\iic$ in dynamical percolation, while the versions for $\muA_r$ will be used in Section~\ref{s.FETIC} to prove that $\fetic\not=\iic$. The finite versions of our results will be slightly stronger than the infinite ones, in that they identify not only a moment where we get $\IIC^{B_r}$ or $\IIC_r$, but also an equality of entire processes. We will use the stronger, dynamic version for $\muA_r$ in Section~\ref{s.FETIC}.

\begin{lemma}[Finite $r$ quenched sampling]\label{l.IICrq}
Let $\{\omega(s): s\in [0,\infty)\}$ be dynamical percolation in $B_r$. Let  $\overline\chi_{r,T}\in\R$ be a random time sampled from $\muA_r/\muA_r[0,T]$, defined only when $\muA_r[0,T]>0$. Then, the finite dimensional distributions of $\big\{\omega(\overline\chi_{r,T}+s) : s\in [0,\infty)\big\}$ converge for almost all $\omega$ as $T \to \infty$ to those of standard dynamical percolation started from $\iic_r$ at time zero. Moreover, the law of the entire process in the Skorokhod topology converges in probability to the same limit process.

Similarly, if $\chi_{r,T}\in\R$ is a random time sampled from $\muO_r/\muO_r[0,T]$, then the same results hold for the process $\big\{\omega(\chi_{r,T}+s) : s\in [0,\infty)\big\}$.
\end{lemma}

\begin{lemma}[Finite $r$ annealed sampling]\label{l.IICra}\ 
\begin{itemize}
\item[{\bf (a)}] Let $\{\overline\omega^*(s) : s\in [0,\infty)\}$ be dynamical percolation in $B_r$ size-biased by $\muA_r[0,T]$, and $\overline\chi^*_{r,T}\in\R$ be a random time with law $\muA_r/\muA_r[0,T]$ for $\muA_r=\muA_r(\overline\omega^*)$. Then the process  $\big\{\overline\omega^*(\overline\chi^*_{r,T}+s) : s\in [0,\infty)\big\}$ is equal in law to standard dynamical percolation started from $\iic_r$ at time zero.

Similarly, if $\{\omega^*(s) : s\in [0,\infty)\}$ is dynamical percolation in $B_r$ size-biased by $\muO_r[0,T]$, and $\chi^*_{r,T}\in\R$ is a random time with law $\muO_r/\muO_r[0,T]$ for $\muO_r=\muO_r(\omega^*)$, then the process  $\big\{\omega^*(\chi^*_{r,T}+s) : s\in [0,\infty)\big\}$ is equal in law to standard dynamical percolation started from $\iic^{B_r}$ at time zero.

\item[{\bf (b)}] The Palm version $(\overline\omega^*,\overline\Pi^*_r)$ of the process $(\omega,\Pi_{\muA_r(\omega)})$ in $B_r$ is  standard dynamical percolation started from $\iic_r$ at time zero. A somewhat concrete way to realize the Palm version is  {\bf Liggett's extra head construction} \cite{Liggett}, see Figure~\ref{f.liggett}:

\begin{figure}[htbp]
\SetLabels
(1*.9)${\muA_r[0,\,\cdot\,]}$\\
(-0.03*.55)$p_1$\\
(-0.03*.6)$p_2$\\
(-0.03*.64)$p_3$\\
(-0.1*.7)$\Theta$\\
(0.16*-0)$q_{r,1}$\\
(0.32*-0)$q_{r,2}$\\
(0.79*-0)\textcolor{red}{$q_{r,J}$}\\
(1*-0)$\Pi_{\muA_r}$\\
\endSetLabels
\centerline{
\AffixLabels{
\epsfysize=4.5in \epsffile{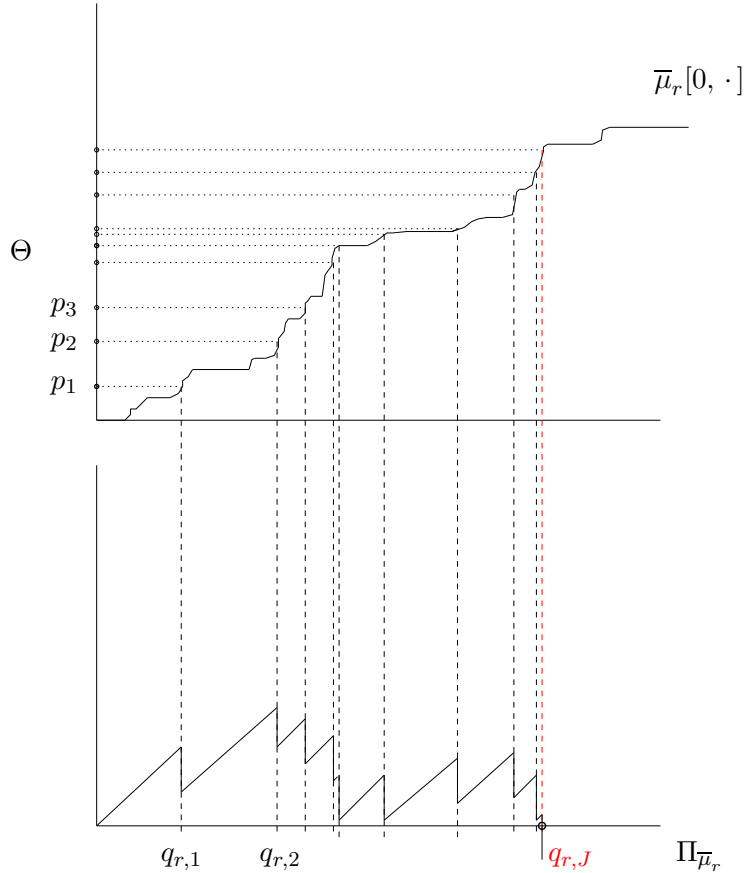}
}
}
\caption{Depicting Liggett's extra head construction.}
\label{f.liggett}
\end{figure}

Let $\big\{ p_i \in [0,\infty): i \in \N \big\}$ enumerate a Poisson point process $\Theta$ with intensity measure Lebesgue
on $[0,\infty)$, and set $q_{r,i} = \inf \big\{ t > 0: \muA_r[0,t] > p_i \big\}$. Clearly, $\Pi_{\muA_r}:=\big\{q_{r,i} : i \in \N\big\}$ is a Poisson point process with intensity $\muA_r$. Set $m_r:=\E \muA_r[0,1]$. Now let $J \in \N$ be the first integer with $\big|\Pi_{\muA_r}\cap [0,m_r^{-1}J] \big| > J$. Then shifting back time by  $q_{r,J}$  gives the Palm version of $(\omega,\Pi_{\muA_r})$.

Similarly, the Palm version $(\omega^*,\Pi^*_r)$ of the process $(\omega,\Pi_{\muO_r(\omega)})$ in $B_r$ is  standard dynamical percolation started from $\iic^{B_r}$ at time zero.
\end{itemize}
\end{lemma}

It should be intuitively quite clear why the ergodic quenched limits in Lemma~\ref{l.IICrq} lead to the size-biased finite averages in Lemma~\ref{l.IICra}: each dynamic configuration of a finite time interval appears in the ergodic quenched limit with a frequency proportional to its probability.

\proofof{Lemma~\ref{l.IICrq}} 
Note that for any percolation configuration $\zeta$ on $B_r$ satisfying $0\llra r$, by definition, $
\Pb{\omega(\overline\chi_{r,T})=\zeta}= \Eb{\int_0^T \1\{\omega_t=\zeta\} \, \d t  / \muA_r[0,T]}$, where the event on the left-hand side is taken to be unsatisfied and the ratio on the right-hand side is taken to be zero on the event that $\muA_r[0,T]=0$. Similarly and more generally, for any time instances $0=s_0\leq s_1 \leq \dots \leq s_k$ and configurations $\zeta_0,\dots,\zeta_k$,
\begin{equation}\label{e.zetaint}
\PB{\omega(\chi_{r,T}+s_i)=\zeta_i,\  i=0,\dots, k}= \EB{ \frac{1}{ \muO_r[0,T]} \int_0^T \MO_r(\zeta_0) \prod_{i=0}^k \1\{\omega_{t+s_i}=\zeta_i\} \, \d t }\,,
\end{equation}
where the random variables on both sides are again interpreted appropriately if $0 \llra r$ at no time in $[0,T]$. There is a very similar multi-point formula in the case of $\overline\chi_{r,T}$; in fact, the entire argument for the first part of the lemma runs in parallel to that for the second, and we omit it.

Dynamical percolation in $B_r$ is a tail trivial process, hence not only is it ergodic, but also the process of the entire configuration along a finite time interval is ergodic. Thus, by the ergodic theorem and the Markov property, the integral in~(\ref{e.zetaint}), divided by $T$, converges almost surely as $T\to\infty$ to 
\begin{equation}\label{e.zetaprob}
\IIC^{B_r}(\zeta_0)\prod_{i=0}^{k-1} \Pb{\omega_{s_{i+1}}=\zeta_{i+1} \md \omega_{s_i}=\zeta_i}\,,
\end{equation}
while 
$\muO_r[0,T]/T\to \E\muO_r[0,1]=1$, almost surely. Therefore, in~(\ref{e.zetaint}), we are taking the expectation of a random variable that converges almost surely to the formula in~(\ref{e.zetaprob}). This random variable is bounded, and hence convergence in expectation also follows. We have thus shown that, for almost all $\omega$, the finite dimensional distributions of $\big\{\omega(\chi_{r,T}+s) : s\in [0,\infty)\big\}$ converge as $T \to \infty$ to those of standard dynamical percolation started from $\iic^{B_r}$. 

To ameliorate this conclusion to hold for the Skorokhod topology (but only in probability, not almost surely), note that, alongside finite-dimensional distributional convergence and the c\`adl\`ag nature of all the sample paths concerned, it is enough to argue that, for any given $K > 0$ and $\eps > 0$, the probability that the process $[0,K] \to \R: t \to \omega (\chi_{r,T} + t)$
has two hexagon switches at times differing by less than $\eps$ vanishes in the high $T$ then low $\eps$ limit. To see this, note that the Lebesgue measure of the set $\cA_T$ of times $t \in [0,T]$ such that $[t,t+K]$ contains two such switch times behaves like $a_\eps T (1+o(1))$ as $T \to \infty$, where $\lim_{\eps \to 0} a_\eps = 0$; on the other hand, the Lebesgue measure of the set $\cB_T$ of times $t \in [0,T]$ such that $\omega \big\vert_{B_r}$ is the completely open configuration behaves almost surely like $b T (1+o(1))$ as $T \to \infty$, where $b > 0$. Since the Radon-Nikodym derivative of $\chi_{r,T}$ is maximized by each point in $\cB_T$, we see that 
$\PB{\chi_{r,T} \in \cA_T } \leq {|\cA_T|}/{|\cB_T|} \leq 2{a_\eps}/{b}$ almost surely for $T$ sufficiently high, where $|\cdot|$ denotes Lebesgue measure. Since $a_\eps \searrow 0$ as $\eps \searrow 0$, we verify the claim needed for convergence in the Skorokhod topology, and complete the proof.    
\qed

\proofof{Lemma~\ref{l.IICra}} The Palm version of a stationary process $(\omega,\xi)$ on $\R$, where $\xi$ is a random measure, is defined in \cite[Chapter 11]{Kallen} as follows. For any Borel set $B\subset \R$ of positive 
Lebesgue measure, and any nonnegative  measurable function $f$ on configurations $(\omega,\xi)$, consider 
$\xi_f(B) := \int_B f(\theta_s(\omega,\xi)) \, \xi(ds)$, 
where $\theta_s$ is the shift by $-s$. Then the Palm version is the law defined by 
$Q_{\omega,\xi}[f] := {\E\xi_f(B)}/{\E\xi(B)}$. It is not hard to show that this does not depend on $B$. 

If we take $\xi=\muA_r$ or $\muO_r$ and $B=[0,T]$, then this construction specializes to the processes defined in part (a). Since we know from Lemma~\ref{l.erg} that $(\omega,\muA_r, \muO_r)$ is ergodic, we can apply \cite[Theorem 11.6]{Kallen}, saying that these Palm versions equal the limit processes defined in Lemma~\ref{l.IICrq}, hence the claim of part (a) follows from that lemma. 

For part (b), there will be no difference between the proofs for $\muA_r$ and $\muO_r$, so let us just work with $\muO_r$. Take $\xi=\Pi_{\muO_r}$, and the Borel sets $B_\eps:=(-\eps,\eps)$. \cite[Theorem 11.5]{Kallen} says that the Palm version of $(\omega,\Pi_{\muO_r})$ is the same as conditioning on $|B_\eps\cap\Pi_{\muO_r}|\geq 1$ or on $|B_\eps\cap\Pi_{\muO_r}|=1$, then taking the limit $\eps\to 0$. This is the most common form of taking the ``Palm version of a point process''. Note that for the equivalence of definitions here, we need that $\muO_r$ does not have atoms (by Lemma~\ref{l.noatom}), hence $\Pi_{\muO_r}$ is a simple point process. 

(Let us give a two-sentence intuitive explanation of why the quoted theorem on the equality between the Palm process and the $\eps$-conditioning holds, at least for the time-zero configuration. Since $\muO_r$ has a density, $\MO_r$, for any static percolation configuration $\zeta$ in $B_r$, we have 
\begin{align*}
\Pb{B_\eps\cap\Pi_{\muO_r}\not=\emptyset \md \omega(0)^{B_r}=\zeta,\,\muO_r} 
&=  1 - \exp\Big(-\int_{-\eps}^\eps \MO_r(\omega_t) \, \d t\Big)\\
& \sim 2\eps \, \MO_r(\zeta)\qquad\text{a.s. as }\eps\to 0\,,
\end{align*}
by the Lebesgue integration theorem and Fubini. Therefore, $\MO_r(\zeta)$ being the Radon-Nikodym derivative $\d \IIC^{B_r}/\d\P$, the $\eps$-conditioning gives 
$$\lim_{\eps\to 0}\Pb{\omega(0)^{B_r}=\zeta \md B_\eps\cap\Pi_{\muO_r}\not=\emptyset } = \IIC^{B_r}(\zeta)\,,$$
as desired.)

Since $\Pi_{\muO_r}$ is obtained from $\muO_r$ using independent stationary randomness (the Lebesgue Poisson point process $\Theta$), the $\omega^*$ marginal in the Palm version of $(\omega,\Pi_{\muO_r})$ is the same as in the Palm version of $(\omega,\muO_r)$, which we already described in part~(a).

\begin{figure}[htbp]
\SetLabels
(0.5*0.85)\Green{$T_{x_{1,1}}$}\\
(0.5*0.6)\Green{$T_{x_{1,2}}$}\\
(0.5*0.2)\Green{$T_{x_{1,3}}$}\\
(0.31*1.04)\Green{$x_{1,1}$}\\
(0.34*0.69)\Green{$x_{1,2}$}\\
(0.29*0.35)\Green{$x_{1,3}$}\\
(0.7*0.4)\textcolor{blue}{$\Pi^*_1$}\\
(0.91*0.23)\textcolor{red}{$\Omega^*_1$}\\
(1.03*0.65)\textcolor{red}{$\omega^*_{1,1}$}\\
(1.03*0.55)\textcolor{red}{$\omega^*_{1,2}$}\\
(1.03*0.45)\textcolor{red}{$\omega^*_{1,3}$}\\
(1.03*0.35)\textcolor{red}{$\omega^*_{1,4}$}\\
(0.3*0.0)\textcolor{blue}{$T^{-1}(\Pi^*_1)$}\\
(0.2*-0.15)\textcolor{red}{$T^{-1}\big((\Pi^*_1,\Omega^*_1)\big)$}\\
\endSetLabels
\vskip 0.2 in
\centerline{
\AffixLabels{
\epsfxsize=5in \epsffile{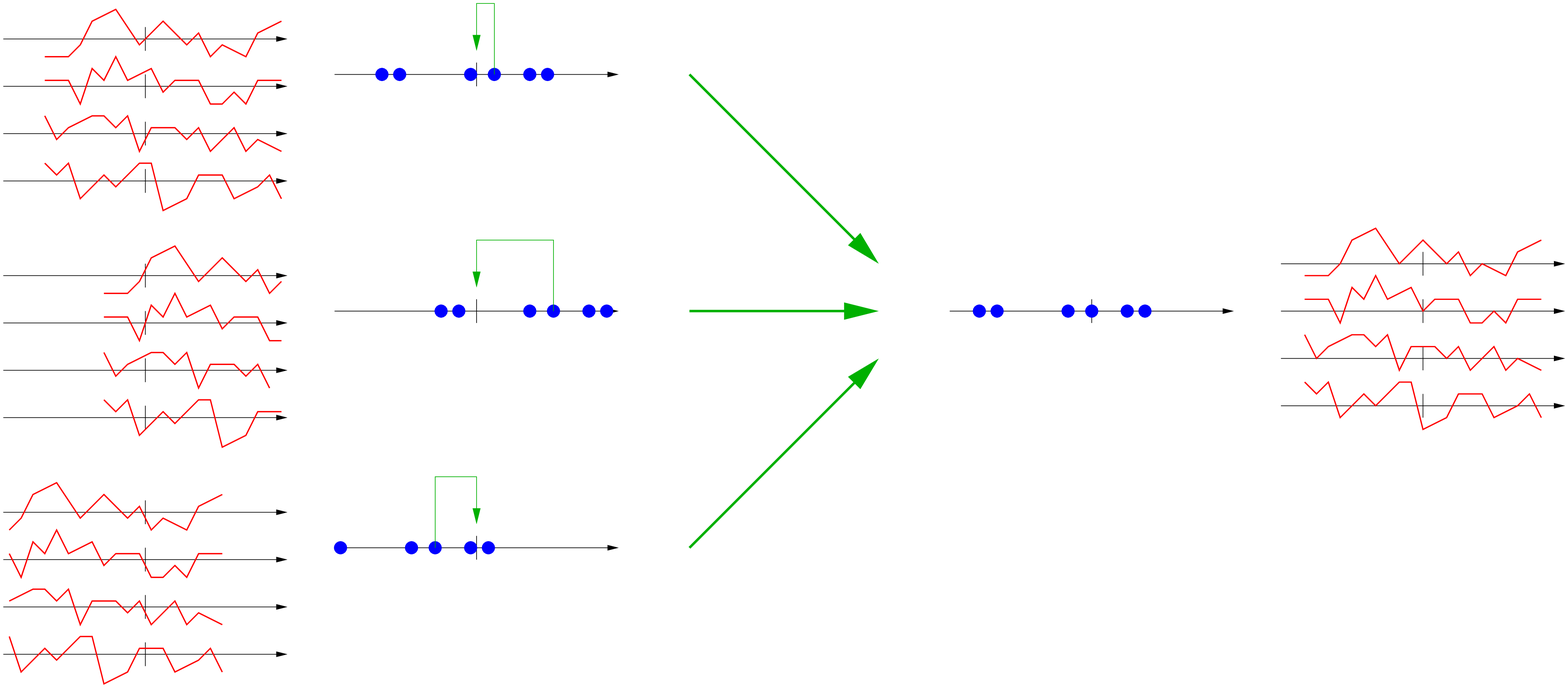}
}
}
\vskip 0.8 in
\SetLabels
(0.5*0.8)\Green{$T_{x_{2,1}}$}\\
(0.5*0.1)\Green{$T_{x_{2,2}}$}\\
(0.3*1.05)\Green{$x_{2,1}$}\\
(0.32*0.43)\Green{$x_{2,2}$}\\
(0.7*0.2)\textcolor{blue}{$\Pi^*_2$}\\
(0.91*0.06)\textcolor{red}{$\Omega^*_2$}\\
(1.03*0.6)\textcolor{red}{$\omega^*_{2,1}$}\\
(1.03*0.4)\textcolor{red}{$\omega^*_{2,2}$}\\
(0.3*-0.15)\textcolor{blue}{$T^{-1}(\Pi^*_2)$}\\
(0.2*-0.4)\textcolor{red}{$T^{-1}\big((\Pi^*_2,\Omega^*_2)\big)$}\\
\endSetLabels
\centerline{
\AffixLabels{
\epsfxsize=5in \epsffile{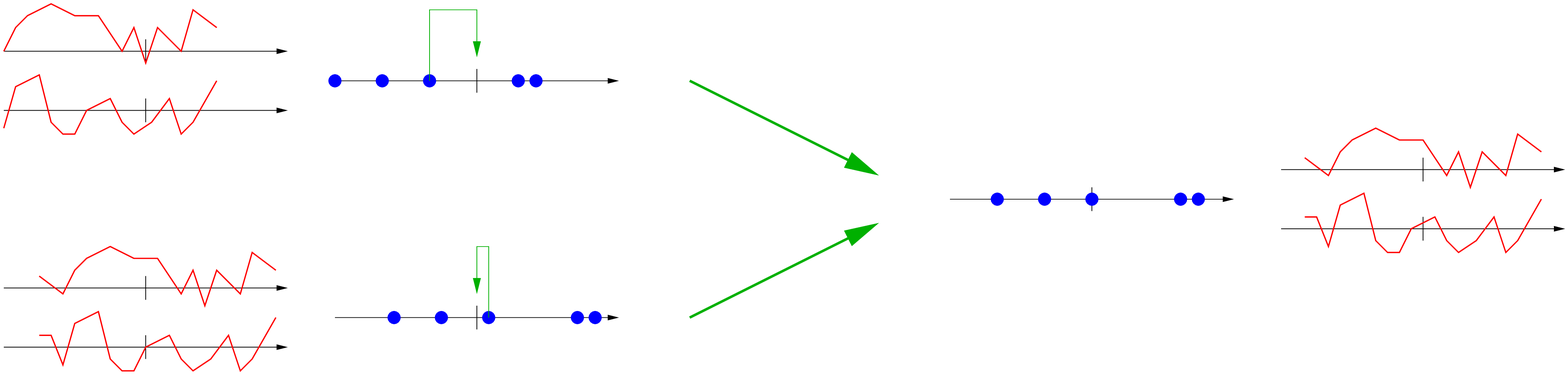}
}
}
\vskip 0.5 in
\caption{A schematic picture of the effect of Liggett's extra head time shift $T$ on $(\Pi_{\muO_r},\omega)$. (For simplicity, the figure pretends that $\Pi$ is a measurable function of $\omega$, instead of $(\omega,\Theta)$.) Different Palm point process realizations $\Pi^*_1$ and $\Pi^*_2$ may arise from a different ``amount'' of Palm dynamical percolation realizations $\Omega^*_1=\{\omega^*_{1,i}:i\in I_1\}$ and $\Omega^*_2=\{\omega^*_{2,i}:i\in I_2\}$ (a ratio 4:2 on the right side of the picture), and the preimages $T^{-1}(\Pi_1^*)$ and $T^{-1}(\Pi_2^*)$ might also have different sizes (which gives the reweighting of the Palm measure compared to the ordinary measure, a ratio 3:2 in the middle of the picture). The product of these ratios is the same as the ratio for the preimages $T^{-1}\big((\Pi^*_1,\Omega_1^*)\big)$ and $T^{-1}\big((\Pi^*_2,\Omega_2^*)\big)$.}
\label{f.ExtraPalm}
\end{figure}

Finally, regarding Liggett's extra head construction, \cite[Corollary 4.18]{Liggett} says that shifting back by $q_{r,J}$ as defined in the statement of part (b) produces the Palm version of $\Pi_{\muO_r}$. Now we need to extend this result from the marginal $\Pi_{\muO_r}$ to $(\omega,\Pi_{\muO_r})$; we will certainly need to use that Liggett's shift coupling acts nicely also on the level of $\omega$ and $\Theta$, since the result clearly would not hold for an arbitrary measurable map $(\omega,\Theta) \mapsto f(\omega,\Theta)$ with the property that $\Pi_{\muO_r(f(\omega,\Theta))} \stackrel{d}{=} \Pi^*_{\muO_r(\omega)}$. The niceness of Liggett's construction lies  in the fact that it gives a random time shift $T_{J_{q,r}}$ that is measurable with respect to~$\Pi_{\muO_r}$, where each time shift $T_x$ is a measure-preserving transformation on the space of configurations $(\omega,\Theta)$. Therefore, if $\cA$ and $\cB$ are arbitrary events for the Palm version $\Pi_{\muO_r}^*$, and $\tilde\cA$ and $\tilde\cB$ are the events for $(\omega^*,\Pi_{\muO_r}^*)$ that project to $\cA$ and $\cB$ in the second coordinate, then 
$$\frac{\P^*(\cA)}{\P^*(\cB)}=\frac{\Pb{T^{-1}(\cA)}}{\Pb{T^{-1}(\cB)}}=\frac{\Pb{T^{-1}(\tilde\cA)}}{\Pb{T^{-1}(\tilde\cB)}}\,,$$
whenever the denominator on either side of this equation is positive. See Figure~\ref{f.ExtraPalm}. Since $\frac{\P^*(\cA)}{\P^*(\cB)}=\frac{\P^*(\tilde\cA)}{\P^*(\tilde\cB)}$ by definition, we get that the effect of $T$ is the same as conditioning on $\{0\in\Pi_{\muO_r}\}$ not only on $\Pi_{\muO_r}$ but also on $(\omega,\Pi_{\muO_r})$, and we are done.
\qed

\medskip
We can now turn to sampling from the limit measure $\muO[0,T]$.

\proofof{Theorem~\ref{t.IICq}} We must argue that $\mu[0,T]>0$ for all $T$ sufficiently high, and also that, for each $r \in \N$, $\omega(\chi_T)^{B_r}$ converges weakly, as $T \to \infty$, to $\iic^{B_r}$. 

For $n \in \N$, set $I_i^n  = [i/n,(i+1/n))$. For $R \in \N$, define $f_R^n: [0,\infty) \to [0,\infty)$ according to $f_R^n(x) = n\mu_R(I_i^n)$ if $x \in I_i^n$ for $i \in \N$. Similarly define $f_\infty^n: [0,\infty) \to [0,\infty)$ according to $f_\infty^n(x) = n\mu(I_i^n)$ if $x \in I_i^n$ for $i \in \N$. 
We now argue that,  for each $\eps > 0$ and $n \in \N$, there exists $R \in \N$ such that 
\begin{equation}\label{e.mumur}
  \limsup_T   \int_0^T \Big\vert    \tfrac{f_R^n(t)}{\int_0^T f_R^n(s) ds} -  \tfrac{f_\infty^n(t)}{\int_0^T f_\infty^n(s) ds} \Big\vert {\rm d} t  \leq \eps \, .
\end{equation}
Note that this assertion allows us to  construct couplings $\coup$  of $\chi_{R,T}$ and $\chi_T$ with the following property:  for each $\eps > 0$ and $n \in \N$, there exists $R \in \N$ such that, for all large enough $T$ simultaneously, $\chi_{R,T}$ and $\chi_T$  are coupled under $\coup=\coup_{R}$ so that 
\begin{equation}\label{e.coup}
\limsup_T \coup \big( \vert \chi_{R,T} - \chi_T \vert \geq 1/n \big) \leq \eps \, .
\end{equation}

Note that~(\ref{e.mumur}) is a consequence of the next three assertions.
First, for each $\eps > 0$, there exists $R \in \N$ such that 
\begin{equation}\label{e.mumurone}
  \limsup_T  T^{-1} \int_0^T \big\vert    f_R^n(t) - f_\infty^n(t) \big\vert  {\rm d} t  \leq \eps \, , \, \, \textrm{$\P$-almost surely} \, .
\end{equation}

Second, for each $\eps > 0$, and for this same value of $R \in \N$, 
\begin{equation}\label{e.mumurtwo}
\limsup_T  T^{-1} \Big\vert \int_0^T f_R^n(s)  {\rm d} s - \int_0^T f_\infty^n(s)  {\rm d} s \Big\vert \leq \eps  \, , \, \, \textrm{$\P$-almost surely} \, .
\end{equation}

Third,
\begin{equation}\label{e.mumurthree}
\lim_T  T^{-1} \int_0^T f_\infty^n(s)  {\rm d} s = 1 \, , \, \, \textrm{$\P$-almost surely} \, .
\end{equation}

We now justify (\ref{e.mumurone}),  (\ref{e.mumurtwo}) and (\ref{e.mumurthree}).

To confirm (\ref{e.mumurone}), note that   $\int_0^{1/n} \vert    f_R^n(t) - f_\infty^n(t) \vert  {\rm d} t = \int \vert \mu(0,1/n) - \mu_R(0,1/n) \vert  {\rm d} \P$.  We fix $R \in \N$ by Theorem~\ref{t.basic} so that $\E \int \vert \mu(0,1/n) - \mu_R(0,1/n) \vert  {\rm d}  \P \leq \eps/n$. Lemma~\ref{l.erg} then provides~(\ref{e.mumurone}).

Note that (\ref{e.mumurtwo}) is a trivial consequence of (\ref{e.mumurone}).

To show (\ref{e.mumurthree}), note that, 
by definition, $\Es{\muO_r[0,1]} = 1$ for each $r \in \N$. Thus Theorem~\ref{t.basic} implies that $\Es{\mu[0,1]} = 1$. Lemma~\ref{l.erg} then implies that 
$\lim_T  T^{-1} \mu(0,T) = 1$, $\P$-almost surely. This limit coincides with that in~(\ref{e.mumurthree}), which establishes this claim. Note that in this derivation we have confirmed that indeed $\mu(0,T) > 0$ for $T$ sufficiently high, $\P$-almost surely.

We conclude the proof by arguing that, for each $\eps > 0$ and $r \in \N$, there exists $R \in \N$ such that  
\begin{equation}\label{e.coupprop}
\liminf_T \coup \Big( \omega(\chi_{R,T})^{B_r} = \omega(\chi_T)^{B_r} \Big) \geq 1 - \eps \, .
\end{equation} 
This indeed suffices for Theorem~\ref{t.IICq}, by the following argument. Recall that we must argue that, for each $r \in \N$, $\omega(\chi_T)^{B_r}$ converges weakly as $T \to \infty$ to $\iic^{B_r}$. We know by Lemma~\ref{l.IICrq} that the weak limit as $T \to \infty$ of $\omega(\chi_{R,T})$ equals $\IIC^{B_R}$. Thus, fixing any $\eps > 0$ and any $R\geq r$, for large enough $T$, the total variation distance between $\omega(\chi_{R,T})^{B_r}$ and $\IIC^{B_r}$ is at most $\eps$. (Note here that on the discrete topological space $\{0,1\}^{B_r}$, convergence in law is the same as in total variation distance.) On the other hand, by~(\ref{e.coupprop}), $\omega(\chi_T)$ coincides with $\omega(\chi_{R,T})$ on $B_r$ with $\coup$-probability at least $1-\eps$ for all high enough $T$. Thus the total variation distance between $\omega(\chi_T)^{B_r}$ and $\IIC^{B_r}$ becomes less than $2\eps$, and we are done.

It remains only to verify (\ref{e.coupprop}). In light of~(\ref{e.coup}), it is enough to argue that, for given $\eps > 0$ and $r \in \N$, there exists $n \in \N$ such that, for all $R\geq r$ and all $T$ sufficiently high, the $\coup$-probability that a hexagon in $B_r$ flips during $\big[\chi_{R,T} - 1/n, \chi_{R,T} + 1/n \big]$ is at most $\eps$. However, by Lemma~\ref{l.IICrq}, the times of hexagon flips in $B_r$ during $\big[\chi_{R,T} - 1/n, \chi_{R,T} + 1/n \big]$, shiftward backwards in time by $\chi_{R,T}$, converges weakly as $T \to \infty$ to a Poisson process of rate $\vert B_r \vert/2$ on $[-1/n,1/n]$. Choosing $n \geq C_\eps r^2$ thus gives the desired statement. \qed

\proofof{Theorem~\ref{t.IICa}} Part (a) follows from Theorem~\ref{t.IICq} --- by Lemma~\ref{l.erg} and  \cite[Theorem 11.6]{Kallen} ---  just as Lemma~\ref{l.IICra} followed from Lemma~\ref{l.IICrq}.

Part~(b) follows from Lemma~\ref{l.IICra}(b) and the next two lemmas. \qed

\begin{lemma}\label{l.exclim}
If $\tau\in\Exc$ is an exceptional time, and $\tau_n\to\tau$, then, for any $r>0$, we have
$\omega(\tau_n)^{B_r} =    
\omega(\tau)^{B_r}
$
for all sufficiently high $n$. 
\end{lemma}

\proof
By Lemma~\ref{l.exctop}(i),  there is  an open interval $I$ which contains
the exceptional time $\tau$ such that, for $t \in I$, $\omega(t)^{B_r}
= \omega(\tau)^{B_r}$; hence $\tau_n\to\tau$ implies the lemma.
\qed
\medskip 

\begin{lemma}\label{lemqjlim}
For the times $q_{r,J}$ defined in Lemma~\ref{l.IICra}~(b), the limit $q_J:=\lim_{r \to \infty} q_{r,J}$ exists almost surely and is an exceptional time.
\end{lemma}

\proof If $f_n: [0,\infty) \to [0,\infty)$ is a sequence of non-decreasing
functions converging pointwise on $[0,\infty)$ to a function
$f:[0,\infty) \to [0,\infty)$, and we write
$$
f_n^{-1}(x) = \inf \Big\{ t > 0 : f_n(t) > x \Big\},
$$
then, whenever $x \in [0,\infty)$ is a point of increase of $f:[0,\infty) \to
[0,\infty)$, we have that $\lim_{n \to \infty} f_n^{-1}(x)  = f^{-1}(x)$. 
The following thus suffices for Lemma \ref{lemqjlim}:

\begin{lemma}
For any $\rho \in \Theta$, 
$t_\rho : =  \inf \big\{ s > 0: \mu(0,s) > \rho \big\}$ is almost surely a point
of increase of $\mu(0,\cdot)$; in particular, it is contained in the support of $\mu$.
\end{lemma}

\proof
Note that the set of $\rho \in (0,\infty)$ for which $t_\rho$ is a point of
increase of $\mu(0,\cdot)$ is given by $\R \setminus \mu(0,Q)$, with $\mu(0,Q) = \big\{ \mu(0,q): q \in Q \big\}$, where $Q$ is the
collection of left-hand endpoints of intervals comprising ${\rm
  supp}(\mu)^c$. Note that $Q$ is countable, and, thus, is so $\mu(0,Q)$. 
Thus, $\Theta \cap \mu(0,Q) = \emptyset$ a.s., because $\Theta$ is
independent of $\mu(0,Q)$.  \qed


\section{$\FETIC$ is not $\IIC$}\label{s.FETIC}

In this section, we prove Theorem \ref{t.FETIC}.

\subsection{The skeleton of the argument}\label{ss.skeleton}

\begin{definition}
Let $\omega$ be a sample of dynamical percolation in the $R$-ball $B_R$. Write $\exc_R$ for the set of times such that $0 \leftrightarrow R$. Let $\fet_R = \inf \big\{ t \geq 0:  0 \build\longleftrightarrow_{}^{\omega_t}  R  \big\}$, and let 
$\FETIC_R$ be the law of $\omega_{\FET_R}$ conditioned on $\FET_R>0$. (Since $\Cl_0(\omega_0)$ is almost surely finite, it takes positive time for the first bit on its boundary to change, and hence the event $\FET_R>0$ is the same as $0\hskip 0.2 cm\not\hskip -0.2 cm\llra R$ in $\omega_0$, which is almost surely satisfied for large enough $R$.) 
\end{definition}

These finite approximations will be very useful. On the one hand, $\iic_R$ is the law of the configuration at a typical point of $\exc_R$, as we saw in Lemmas~\ref{l.IICrq} and~\ref{l.IICra}. On the other hand, by \cite[Lemma 4.5]{HMP}, we have that $\fet_R\to\fet$ almost surely as $R\to\infty$; hence $\fetic_R$ converges to $\fetic$ in law (by Lemma~\ref{l.exclim}).




There is a natural line of attack if we want to distinguish $\iic_R$ from $\fetic_R$. Let us call the left-isolated points of $\exc_R$ {\bf arrivals}, and the right-isolated points {\bf departures}. As we will see, the law of a {\em typical} arrival configuration can be easily obtained from $\iic_R$ (and will be denoted by $\iic_R'$): we get it by size-biasing with respect to the number of pivotal hexagons for the event $\{0 \lra R \}$. This is different from $\iic_R$, but not by much: it can be shown (though we will not do so) 
that its weak limit as $R\to\infty$ coincides with $\iic$. However, $\fetic_R$ is not given by a typical arrival: as usual when waiting for the first arrival of a stationary point process, the time between $\fet_R$ and the last departure before it (somewhere in the negative half-line) is a size-biased sample of the typical reconnection time between departures and arrivals, and if an arrival configuration typically occurs at the end of longer disconnection intervals, then it is more 
likely to appear in $\FETIC_R$. Since it is harder to think about dynamical percolation ending at a certain configuration than about starting it at such a configuration, our strategy to understand $\fetic_R$ will be to reverse time, start dynamical percolation from certain typical $\iic_R'$ configurations, condition on immediate termination of $\{0 \lra R \}$, and then estimate the expected time of reconnection. If we can exhibit two events at time zero that have the same positive probability under the limit measure $\iic$, but for which the expected reconnection times differ, then these events will turn out to have different probabilities under $\fetic$, and we will be done.

Roughly, of these two events under $\iic_R'$,  the first will be that the configuration looks ``normal'' in a bounded neighbourhood of $0$, while the second will be that the configuration is ``thinner'' in the same neighbourhood. (We will in fact define a thinning procedure on normal static configurations satisfying $0 \lra R$, changing the configuration in a bounded neighbourhood of $0$.) A thinner configuration falls apart more easily, and hence reconnects to distance $R$ with more difficulty; and so, one may expect that such a configuration is more probable under $\fetic_R$ than is a normal configuration, which is to say, $\fetic_R$ is thinner than $\iic_R$. This is certainly the case if the thin configuration is, say, given by a single straight line segment of open hexagons from $0$ to $\p B_R$, with all other hexagons in $B_R$ being closed. However, this $R$-dependent configuration has a vanishing probability in the limit measure $\iic$; therefore, while the imbalance in probability of this 
configuration 
distinguishes $\fetic_R$ from $\iic_R$, a distinction between $\fetic$ and $\iic$ cannot be deduced. This is why we want to require the configuration to be thin only in a bounded neighbourhood of $0$. However, the main difficulty now is that normal reconnection times are very short if $R$ is large, and that, with high probability, the configuration is entirely static in a bounded neighbourhood of the origin; hence it is not clear that our thinning will have a noticeable effect on the reconnection time. The solution will be that the expected reconnection time, though tiny, turns out to be dominated by times that are macroscopically large (independently of $R$): large enough that if the configuration close to $0$ is thin then it does indeed start falling apart, making expected reconnection time noticeably larger when the thinning procedure has been applied. To argue this, we will need the result from \cite{HMP} that $\fet$ has finite expectation (in fact, an exponential tail): this will tell us that the normal 
reconnection time is well behaved, making it possible to 
prove that, in expectation, it is strictly dominated by the reconnection time of thinned configurations.
\medskip

In this introductory subsection, we first explain the time-reversal and the size-biasing effects determining the relationship between $\iic_R$ and $\fetic_R$, then define the thinning procedure, and will finally show that a noticeable difference between expected reconnection times indeed implies that $\fetic$ and $\iic$ are different. In the subsequent subsections, we will prove that there is such a difference.

Recall from the above discussion that, in standard c\`adl\`ag dynamical percolation, a time $t\in\exc_R$ for which there exists $\epsilon > 0$ such that $[t - \epsilon,t)\cap\exc_R=\emptyset$ is called an arrival. Write $\arrival_R$ for the set of arrivals. Furthermore, for a static percolation configuration $\zeta$ in $B_R$ that satisfies $0\lra R$, denote by $\Piv=\Piv_{0 \lra R}(\zeta)$ the set of hexagons in $B_R$ that are pivotal in the configuration $\zeta$ for $0 \leftrightarrow R$, and recall that $\iic'_R$ denotes the law on configurations in $B_R$ whose Radon-Nikodym derivative with respect to $\iic_R$ is given by $\vert \Piv \vert$ up to normalization. 

\begin{lemma}\label{l.IIC'arrival}
The following three definitions for the process $\Ps{\,\cdot \md 0\in \arrival_R}$ are equivalent:
\begin{itemize}
\item[{\bf (i)}]  consider dynamical percolation in $B_R$ conditionally on the event $\{ 0 \leftrightarrow R \}$ occurring at time $0$ but not at time $-\epsilon$, and take the weak limit as $\epsilon \downarrow 0$;
\item[{\bf (ii)}] for large $T>0$, pick uniformly an element $\tau \in \arrival_R\cap [0,T]$, consider the shifted dynamical percolation configuration $\{\omega_{t-\tau} : t\in\R\}$, and take the weak limit as $T\to\infty$ (conditionally on $\omega$, or averaged);
\item[{\bf (iii)}] let $\omega_0$ be distributed according to $\iic_R'$, choose uniformly an element $S \in \Piv(\omega_0)$, obtain the configuration $\omega_{0^-}$ by closing the hexagon $S$, and let the rest of the evolution $\{ \omega_t: t \in \R \}$ be given by c\`adl\`ag dynamical percolation updates independently of the values of $\omega_0$ and $S$.
\end{itemize}
\end{lemma}

\proof
While the weak limits in (i) and (ii) might not exist a priori, the definition of (iii) is clearly well formulated. We first prove the equivalence of  (i) and (iii), implying the existence of the weak limit in (i), in particular. It is enough to show that, for all configurations $\zeta$ in $B_R$ such that $0 \leftrightarrow R$,
\begin{equation}\label{e.compwmiic}
 \lim_{\eps\to 0}\frac{{\rm d } \Pb{\, \cdot \md 0 \in \exc_R,\, -\eps\not\in\exc_R}}{{\rm d}\IIC_R}(\zeta) =  Z_1^{-1} \vert \Piv_{0 \leftrightarrow R}(\zeta) \vert \, ,
\end{equation}
where  $Z_1 \in (0,\infty)$ is a normalization.

Given a configuration $\zeta$ such that $0 \leftrightarrow R$, let $p_\epsilon(\zeta)$  be the probability that dynamical percolation given $\omega_0 = \zeta$ satisfies $0 \, \not \!\leftrightarrow R$ at time $-\epsilon$. If $\eps$ is tiny (depending on $R$), then the probability of having at least two hexagons flipping in the time interval $(-\eps,0)$ is much less than the probability of any specific hexagon flip. Therefore, $\lim_{\epsilon \to 0} p_\epsilon(\zeta) / \eps =  \vert  \Piv_{0 \leftrightarrow R}(\zeta) \vert$, which implies (\ref{e.compwmiic}).

To prove the equivalence of (ii) and (iii), let us reformulate the $T$-dependent law defined in (ii) as taking uniformly one from all pairs of configurations $(\omega_{t-},\omega_t) \in \exc_R^c \times \exc_R$, with $t\in [0,T]$, and then running dynamical percolation in the two directions from here. By the ergodicity of $\{\omega_t: t\in\R\}$ (Lemma~\ref{l.erg}), the weak limit of this law is the same as taking a pair of static configurations $(\zeta_1,\zeta_2)$ that differ only in one hexagon such that $0\lra R$ in $\zeta_2$ but not in $\zeta_1$ to start the dynamics. This is clearly the same as the law defined in (iii).

The equivalence of (i) and (ii) follows from the above two equivalences; or, just like in Lemma~\ref{l.IICra}, we can also quote \cite[Theorem 11.6]{Kallen} on the equivalent definitions of the Palm version of the process $(\omega,\arrival_R)$. 
\qed
\medskip

Now, as we promised, in order to understand the effect of waiting for the first exceptional time on the distribution of the configuration at that time, we time-reverse the dynamics, started from typical arrival times:

\begin{definition}
Let $\pnorm$ denote the time-reversal of $\Ps{\, \cdot \md 0\in\arrival_R}$ (i.e., $t\mapsto -t$ for all $t\in\R$). More explicitly, it is the c\`a\textbf{g}l\`a\textbf{d} (left-continuous with right limits) Markov process given as follows. Under $\pnorm$, the distribution of $\omega_0$ is $\iic'_R$. Given $\omega_0$, a uniform element $S \in \Piv$ is selected, with the configuration $\omega_{0^+}$ being set equal to $\omega_0$ modified by closing the hexagon $S$. The rest of the evolution of $\{ \omega_t: t \in\R\}$ is given by c\`agl\`ad dynamical percolation updates independently of the values of $\omega_0$ and $S$. 
\end{definition}


\begin{lemma}\label{l.weights}
Under the law $\pnorm$, recall that $0 \lra R$ is satisfied by $\omega_0$ but not by $\omega_{0^+}$; let the reconnection time $N \in (0,\infty)$ be given by  $N = \inf \big\{ t >0: 0 \build\longleftrightarrow_{}^{\omega_t}  R 
\big\}$. For each static $B_R$ configuration $\zeta$, we have that
$$
 \frac{{\rm d} \FETIC_R}{{\rm d} \iic_R}(\zeta) = Z^{-1} \Eso{\rm norm}{N \md \omega_0 = \zeta}  \, \vert  \Piv_{0 \leftrightarrow R} \vert \, ,
$$  
where $Z \in (0,\infty)$ is a normalization. 
\end{lemma}

\proof 
We claim that
\begin{equation}\label{e.compfetwm}
 \frac{{\rm d} \FETIC_R}{{\rm d } \Pb{\, \cdot \md 0 \in \arrival_R}}(\zeta) =  Z_2^{-1}  \Eso{\rm norm}{N \md \omega_0 = \zeta} \, ,
\end{equation}
where  $Z_2 \in (0,\infty)$ is another normalization. From (\ref{e.compwmiic}) and (\ref{e.compfetwm}) follows the statement of the lemma.

To prove (\ref{e.compfetwm}), let $\phi:\exc_R^c \longrightarrow \arrival_R$ associate to each moment of disconnection $0 \, \not \!\leftrightarrow  R$ in  c\`adl\`ag dynamical percolation the first connection time to its right (which is necessarily an arrival). Condition the process on $\exc_R^c \cap [-n,0] \not= \emptyset$ and pick a random time $\chi$ whose conditional law is given by normalized Lebesgue measure on $\exc_R^c \cap [-n,0]$; note that $\fetic_R$ is the weak limit as $n \to \infty$ of $\omega_{\phi(\chi)}$. Note that, in this weak  limit, 
the probability that $\omega_{\phi(\chi)}$ is a given static configuration $\zeta$ (for which $0 \leftrightarrow R$) is proportional to the mean length of an interval in $\exc_R^c$ at whose right-hand endpoint the configuration is $\zeta$. Thus we obtain~(\ref{e.compfetwm}). \qed
\medskip

Here is a straightforward variant of (\ref{e.compfetwm}). For any non-negative random variable $X$ of finite mean, $\what{X}$ will denote the size-biased version; i.e., $\Pb{\what{X} \geq t} = \Eb{X}^{-1} \Eb{X \1_{X \geq t}}$.

\begin{lemma}\label{lemsizebiasing}
Let $\what{N}$ be the size-biased version of the reconnection time $N$ under the law $\pnorm$, and let $U$ be an independent $\mathsf{Unif}[0,1]$ random variable. Then $\what{N} \, U$ has the distribution of $\fet_R$.
\end{lemma}

The following useful fact was proved in \cite{HMP}. 

\begin{lemma}\label{l.expdecay}
In dynamical percolation we have
$$
\Pb{ \fet_R > t } \leq \exp \big\{ - ct \big\}
$$
for all $t > 0$,
where $c > 0$ may be chosen uniformly in $R \in \N$.
\end{lemma}

Note that the preceding two lemmas imply that 
$\Ps{\what N > t}\leq \exp \{ -c\,t \}$, uniformly in $R$. In particular, this random variable has finite moments: for each $k \in \N$, $\Es{{\what N}^k} = \Es{N^{k+1}}/\E N < \infty$, again uniformly in $R$.
\medskip

We now introduce the thinning procedure which is central to our technique for showing that $\fetic$ differs from $\iic$. 

\begin{definition}
A circuit $\Gamma$ is a finite self-avoiding path of hexagons such that for no vertex in the hexagonal lattice are all three of the neighbouring hexagons visited by $\Gamma$ and such that $\hexlatt \setminus \Gamma$ has exactly two connected components: a finite one, denoted by $\intg(\Gamma)$, and an infinite one.
Note that a partial order on circuits $\Gamma$ is provided by containment of the enclosed regions $\intg(\Gamma)$.

Let $\zeta$ be a percolation configuration in $B_R$ such that $0 \leftrightarrow R$.
Note that if some $\zeta$-open circuit $\Gamma$ satisfies $B_r \subseteq \intg(\Gamma)$, then there is a unique $\zeta$-open circuit which encloses $B_r$ and is minimal in the partial order among such circuits. If $\zeta$ is such that this circuit exists, we label the circuit by $\Gamma_r$.
\end{definition}

\begin{definition}\label{d.Fine}
Recall the exponent $\eta \in (0,1)$ from~(\ref{e.a4}), and fix $\eps>0$ small enough that $(1+2\eps)(1-\eta)<1$.
Now assume that $r$ satisfies $r^{2(1+2\eps)}\alpha_4(r^{1+2\eps})<r/2$, which holds for all large enough $r$, by~(\ref{e.piv}). Let $R \in \N$ satisfy $R \geq r^{1 + 2\eps}$. A configuration $\zeta$ in $B_R$ is said to satisfy $\zeta\in\thin$ if the following conditions hold:
\begin{itemize}
\item $0 \llra R$;
\item the circuit $\Gamma_r$ exists and satisfies $\Gamma_r \subseteq B_{r^{1+\eps}}$;
\item the pivotal set $\Piv_{0 \lra \Gamma_r} = \Piv_{0\lra R} \cap \intg(\Gamma_r)$ satisfies $|\Piv_{0 \lra \Gamma_r}(\zeta)| \leq r^{2(1+2\eps)}\alpha_4(r^{1+2\eps})$.
\end{itemize}
Finally, a dynamical configuration $\{\omega_t : t\in\R\}$ is said to satisfy $\thin$ if $\omega_0 \in \thin$.
\end{definition}

\begin{figure}[htbp]
\SetLabels
(-0.03*.7)$B_{r}$\\
(0*.52)$B_{r/2}$\\
(0*.25)$\Gamma$\\
(1.08*.95)pivotals\\
(1.09*.88) for $0\lra\Gamma$\\
(1.08*.82) in $\chi_{r,\Gamma,a}$\\
\endSetLabels
\centerline{
\AffixLabels{
\includegraphics[width=3in]{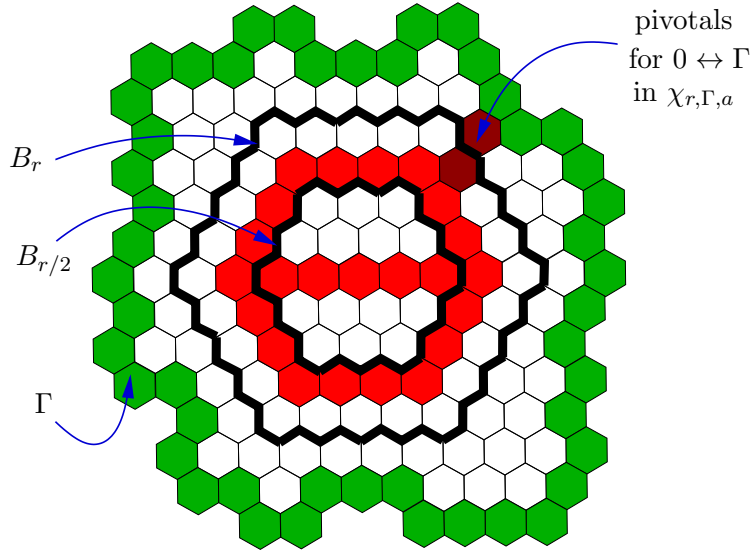}
}}
\caption{The boundary paths delimiting $B_2$ and $B_4$ are black, and the circuit $\Gamma$ is green. The red and dark red hexagons are the open hexagons of $\chi_{4,\Gamma,2}$: the red hexagonal circuit is set in such a way that its distance from $\Gamma$ is $a=2$, and the dark red path is chosen in some arbitrary but fixed way so that it realizes this distance $a$. Note that this dark red path is the set of pivotals for  $0 \lra \Gamma$.}
\label{f.thin}
\end{figure}

\begin{definition}\label{defthin}
Let $r \in \N$ be even. Let $\Gamma$ denote a circuit such that $B_r \subseteq \intg(\Gamma)$. Let $a \in \{0,\ldots,r/2 \}$.
The $(r,\Gamma,a)$-slim configuration $\chi_{r,\Gamma,a}$ is a particular percolation configuration in $\intg(\Gamma)$, as shown in  Figure~\ref{f.thin}, whose set of open hexagons in $\intg(\Gamma) \cap B_{r/2}$ consists of the hexagons in $B_{r/2}$ that intersect the $x$-axis, and for which $\vert \Piv_{0 \lra \Gamma} \vert = a$.
\end{definition}

\begin{definition}\label{defthinning}
The thinning procedure $\thinning=\thinning^R_r$ maps the set of configurations in $B_R$ to itself. Let $\zeta$ be such a configuration. If $\zeta \not\in \thin$, then set $\thinning(\zeta) = \zeta$.
If $\zeta \in \thin$, let $\thinning(\zeta)$ be the configuration in $B_R$ of the following form:
\begin{equation*}
  \thinning(\zeta)(x)  = \left\{ \begin{array}{lr} \zeta(x) &  \textrm{if $x \in B_R \setminus \intg(\Gamma_r)$,}   \\
                  \chi_{r,\Gamma_r,\vert \Piv_{0 \lra \Gamma_r}  \vert}(x) &  \textrm{if $x \in \intg(\Gamma_r)$.}  
  \end{array} \right.  
\end{equation*}  
\end{definition}

We define a coupling of $\pnorm$ with another dynamical process begun by pairing the initial condition with its thinned counterpart. We denote by $\omega'$ the process under $\pnorm$, and write $\omega''$ for the process under the measure $\pthin$ which we now introduce by coupling with $\pnorm$. We set $\pthin$ by choosing its initial condition $\omega''_0 = \thinning(\omega'_0)$; if the hexagon $S$ selected for initial closure in the definition of $\pnorm$ lies in the unbounded component of the complement of $\Gamma_r(\omega'_0)$, we set $S'' = S$; otherwise, we choose $S''$ uniformly among $\Piv_{0 \lra R}(\omega_0'') \cap \intg(\Gamma_r)$. We define $\omega''_{0^+}$ by modifying $\omega''_0$ by closing $S''$. The subsequent evolution of $\omega''$ is made in accordance with the c\`agl\`ad dynamical updates used in defining $\omega'$. Note that there might be updates that do not have an effect on $\omega'$ (the new status coinciding with the old one), and hence are not visible if we see only $\omega'$, 
while do have an effect on $\omega''$; thus $\omega''$ is not entirely measurable w.r.t.~$\omega'$, even though the extra randomness in $\omega''$ is quite simple.

We denote by $\pnorm$ and $\pthin$ the above dynamics and its thinned counterpart, and write $N$ and $\thinrv$ for the reconnection time $\inf\big\{ t > 0:  0 \build\longleftrightarrow_{}^{\omega_t} R \big\}$ under $\pnorm$ and $\pthin$. We will often use the above coupling of the two c\`agl\`ad processes, but will not need a separate notation to denote it. 
The principal result we need is now stated.

\begin{proposition}[Thinned versus Normal]\label{lemvnthin}
As in Definition~\ref{d.Fine}, fix $\eps>0$ small, and consider all large enough $r\in\N$. Then, uniformly in $R \geq r^{1 + 2\eps}$, we have $\tfrac{\Eso{\th}{\thinrv \1_{\thin}}}{\Eso{\no}{N \1_\thin}} \to \infty$ 
 as $r\to\infty$. 
\end{proposition}

\proofof{Theorem~\ref{t.FETIC}, assuming Proposition~\ref{lemvnthin}} 
We want to show that there exists a circuit $\Gamma$ in the annulus $A_{r,r^{1 + \epsilon}}$ and two configurations $\zeta'$ and $\zeta''$ on $\Gamma \cup \intg(\Gamma)$ with $\Gamma = \Gamma_r(\zeta')=\Gamma_r(\zeta'')$, such that $\iic_R(\zeta') = \iic_R(\zeta'')$ for each integer $R \geq r^{1+2\eps}$, with the common value having a positive limit as $R\to\infty$, while $\inf_{R\geq r^{1+2\eps}}\tfrac{\FETIC_R(\zeta'')}{\FETIC_R(\zeta')} > 1$. 

By Proposition \ref{lemvnthin}, we may choose $r \in \N$ so that $\Eso{\th}{\thinrv \1_{\thin}} > 2 \, \Eso{\no}{N \1_{\thin}}$ for all $R$ sufficiently high. Hence, 
 there exists a choice of circuit $\Gamma$ in $A_{r,r^{1+\epsilon}}$, and a configuration $\zeta'$ in $\intg(\Gamma) \cup \Gamma$, such that $\Gamma = \Gamma_r(\zeta')$, the second and third conditions for $\thin$ occur,  
 and, setting $\zeta''$ equal to the restriction of $\thinning(\zeta')$ to $\intg(\Gamma) \cup \Gamma$,
$$ 
\EBo{\no}{N \md \omega_0\big\vert_{\intg(\Gamma) \cup \Gamma} = \zeta'' } > 2 \, \EBo{\no}{N \md \omega_0\big\vert_{\intg(\Gamma) \cup \Gamma} = \zeta'} \, . 
$$ 
It is clear that $\iic_R(\zeta'') = \iic_R(\zeta')$; moreover, $\iic_R'(\zeta'') = \iic_R'(\zeta')$, since the number of pivotals for $\{0 \lra R \}$ is left intact by $\thinning$. Hence, by Lemma~\ref{l.weights}, we have $\fetic_R(\zeta'') > 2 \,\FETIC_R(\zeta')$. \qed
\medskip

The rest of the section will be devoted to the proof of Proposition~\ref{lemvnthin}. Let us start by collecting the main ingredients needed for the proof; these ingredients will then be proved in the remaining subsections.

Thinning will make a difference only if there is enough time before reconnection for the configuration in $\intg(\Gamma_r)$ to change significantly. To this end, as we will see, the events $\{N>1/r\}$ and $\{\thinrv>1/r\}$ will be important to us. How different are these two events? Although the set of open hexagons in $\thinning(\zeta)$ is not exactly a subset of its counterpart for $\zeta$,  we can compare the thinned and normal reconnection times in this regime under a certain event $\good$: 
\begin{equation}\label{e.goodNV}
\good \cap \big\{ N > 1/r \big\} \subseteq \big\{ \thinrv > 1/r \big\}\,,
\end{equation}
where $\good$ is defined as follows (and is applied in the above relation to the configuration before thinning):

\begin{definition}\label{d.good}
Let $R,r \in \N$ satisfy $R \geq r^{1 + 2\eps}$ where $\eps > 0$ is specified in Definition~\ref{d.Fine}. Let $\omega$ be a dynamical configuration in $B_R$. We say that $\omega \in \good$ if the following conditions are satisfied:
\begin{itemize}
 \item $\omega_0\in \thin$, as specified in Definition~\ref{d.Fine}; 
 \item for each $t \in [0,r^{-1}]$, the inner and outer boundaries of the annulus $A_{r^{1 + \epsilon},r^{1 + 2\epsilon}}$ are separated by an $\omega_t$-open circuit;
 \item for each $t \in [0,r^{-1}]$, $0  \build\longleftrightarrow_{}^{\omega_t} r^{1 + 2\epsilon}$. 
\end{itemize}
\end{definition}

Now, to see~(\ref{e.goodNV}), note that the occurrence of $\good$ implies that $0$ is connected to some open circuit $\Gamma = \Gamma(t)$ such that $B_{r^{1 + \eps}} \subseteq \Gamma$ for all $0 \leq t \leq r^{-1}$. Hence, $N > 1/r$ implies that $r^{1 + \eps} \,\,\, \not\!\!\!\llra R$ for all $t \in [0,r^{-1}]$ under $\pnorm$. Since the dynamical percolations under $\pnorm$ and $\pthin$ agree at all positive times in $A_{r^{1+\eps},R}$, we have that $r^{1 + \eps} \,\,\, \not\!\!\!\llra R$ for all $t \in [0,r^{-1}]$ also under $\pthin$. Thus, $\thinrv > 1/r$ and we obtain (\ref{e.goodNV}). 

The event $\good$ is of course useful only if it is reasonably likely to occur.  Proposition~\ref{lemgoodm}, which is the main result of the upcoming Subsection~\ref{ss.N>1/r}, will show that 
$$
\Pso{\no}{\good \md N > 1/r} \geq c_1 \, .
$$
This,~(\ref{e.goodNV}) and $\good \subseteq \thin$ imply the following ``stochastic quasi-domination'' between $\thinrv$ and $N$:
\begin{equation}\label{e.VWdom}
\Pso{\th}{\thinrv>1/r, \thin} \geq \Pso{\no}{N>1/r,\good} \geq c_1 \, \Pso{\no}{N>1/r}\,.
\end{equation}

Although the event $\{ N>1/r \}$ has minute probability when $R$ is large, a large portion of the expectation $\Es{N\1_\thin}$ is contributed by sample points realizing this event. This can be proved using the size-biasing description of the connection time discussed in Lemma~\ref{lemsizebiasing}. Indeed, by some rather general size-biasing arguments, together with the uniform boundedness of the expectation $\Es{\what N_R}<\infty$ (due to Lemmas~\ref{lemsizebiasing} and~\ref{l.expdecay} above), alongside the fact that $\Pso{\no}{\thin \md N > 1/r} \geq c_1$ (due to $\good \subseteq \thin$), it will be proved in Subsection~\ref{ss.sizebias} that
\begin{equation}\label{e.WisOK}
\Ebo{\no}{N \md N > 1/r,\thin} <  C_2 <\infty\,,
\end{equation}
and that
\begin{equation}\label{e.hatWisOK}
\Pbo{\no}{\what{N\1_\thin} > 1/r } = \frac{\Eso{\no}{N\1_{N>1/r} \1_{\thin}}}{\Eso{\no}{ N \1_{\thin}}} >  c_2 > 0\,.
\end{equation}

Finally, as we will prove in Proposition \ref{p.Visgood} of Subsection~\ref{ss.thinned}, should the dynamics begun under $\thinning$ result in at least a short reconnection time, $\thinrv>1/r$, then there is a uniformly positive probability that connection will not be reestablished until very much later: 
\begin{equation}\label{e.Visgood}
\Pbo{\th}{\thinrv > g(r) \md \thinrv > 1/r,\thin } > c_3 > 0\,,
\end{equation}
for some $g(r)\to\infty$ as $r \to \infty$. 
%

\proofof{Proposition \ref{lemvnthin}}
From the above assemblage of facts, we find that 
\begin{align*}
 \Eso{\th}{\thinrv \1_{\thin}} & \geq \Eso{\th}{\thinrv \1_{\thinrv>1/r} \1_{\thin}} \\
 & = \Ebo{\th}{\thinrv \md \thinrv >1/r,\thin} \, \Pso{\th}{\thinrv>1/r,\thin}\\
 &\geq c_3 \, c_1 \, g(r)\, \Pso{\no}{N>1/r,\thin}\,, \qquad\text{by (\ref{e.Visgood}) and  (\ref{e.VWdom})}\\
 &\geq c_3 \, c_1\, g(r)\, \frac{\Eso{\no}{N \1_{N>1/r} \1_{\thin}}}{C_2}\,, \qquad\text{by (\ref{e.WisOK})} \\
 &\geq c_3 \, c_1 \, c_2 \, g(r)\, \frac{\Eso{\no}{N \1_{\thin}}}{C_2}\,,\qquad\text{by (\ref{e.hatWisOK})}\,.
\end{align*}
Therefore, the ratio 
$\tfrac{\Eso{\th}{\thinrv \1_{\thin}}}{\Eso{\no}{N\1_{\thin}}}$ 
tends to infinity as $r\to\infty$, uniformly in $R\geq r^{1+2\eps}$, as required. \qed
\medskip
 
We will now start proving the above ingredients.

\subsection{Understanding the law $\mathbb{P}_{\mathrm{norm}}( \cdot \,|\, N > 1/r)$}\label{ss.N>1/r}

In this section, $\P$ will denote the law of c\`agl\`ad dynamical percolation with time $\R$. Recall that $\exc_R$ is 
the set of times such that $0 \leftrightarrow R$,  now a union of left-open right-closed intervals.

It is hard to understand the conditioned measure $\P':=\Pbo{\no}{ \, \cdot \md N \geq 1/r }$ directly, because the condition has a tiny probability. We will handle this issue by noticing that, for large enough $s\in\Z^+$, we have $\Pb{\Exc_R \cap (s/r,\,(s+1)/r]=\emptyset \md 0\in\Exc_R}> c>0$, uniformly in $r>0$ (see Lemma~\ref{l.temp}), and given the existence of this empty interval, $\gamma:=\sup \{\exc_R\cap [0,s/r)\}$ is a moment such that the reconnection time from it is at least $1/r$. If $s$ is bounded, then the law of dynamical percolation viewed from such a $\gamma$ (to be denoted by $\P''$, see Lemma~\ref{l.rnderiv}) turns out to be not very different from the law $\P'$ (see Lemma~\ref{l.rnd}). Therefore, once we prove that $\omega_t$ has certain good properties with high probability for all $t\in [0,(s+1)/r]$ under $\Pb{ \,\cdot\, \md 0\in\Exc_R,\ \Exc_R \cap (s/r,\,(s+1)/r)=\emptyset}$, which is already a feasible task, and hence that the dynamical configuration viewed from $\gamma$ (i.e., the measure $\P''$) 
is well behaved, we will be able to deduce almost the same for the measure $\P'$; this will be Proposition~\ref{lemgoodm}, the main goal of this subsection.

\begin{figure}[htbp]
\SetLabels
(.17*.92)\textcolor{red}{$j\in\cJ\subset \exc_R$}\\
(.67*.92)$j+\frac{s}{r}$\\
(.78*.92)$j+\frac{s+1}{r}$\\
(.6*.38)\textcolor{red}{$m:=A(j)\in\cM$}\\
(.85*.38)\textcolor{red}{$m+\ell_m$}\\
(.05*.48)$m-\frac{s}{r}$\\
(.32*.48)$m-\frac{s}{r}+\ell_m$\\
(.12*.3)\textcolor{Blue}{$L_m$}\\
(.15*.04)\textcolor{Purple}{$L_m \cap \exc_R =: \cD_m=A^{-1}(m)$}\\
\endSetLabels
\centerline{
\AffixLabels{
\includegraphics[width=\textwidth]{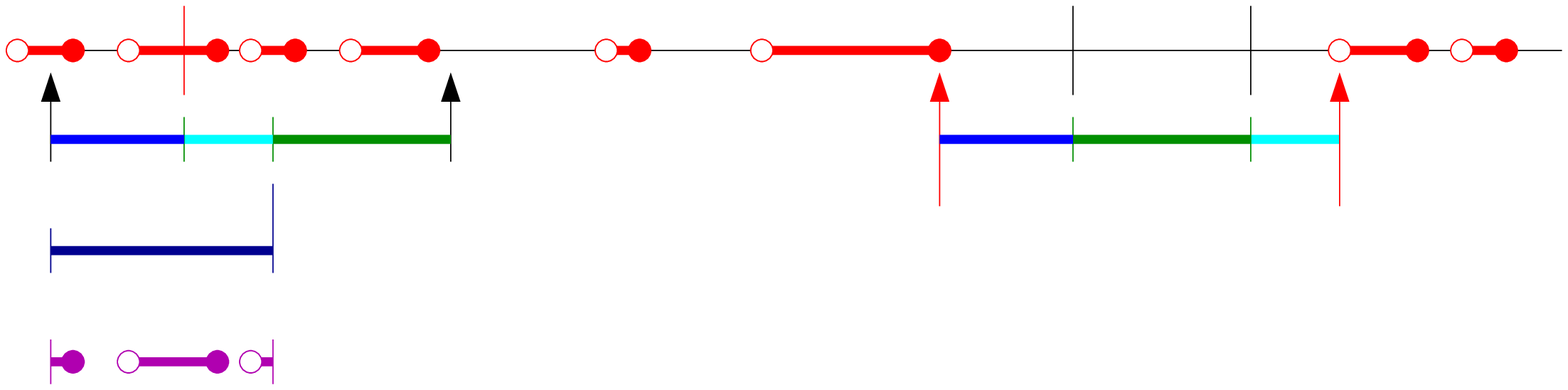}
}}
\caption{The domain of attraction $\cD_m$ appearing in the definition of $\P''$, and the map $A:\cJ\longrightarrow \cM$ appearing in the proof of $\P''=\tilde\P$ (Lemma~\ref{l.rnderiv}).}
\label{f.DOA}
\end{figure}

\begin{definition}
Call an element $x \in \exc_R$ a {\bf marker} if
$\big(x,x + r^{-1} \big] \cap \exc_R = \emptyset$. Write $\cM \subseteq \exc_R$ for the set of markers.
For $x \in \cM$, set $\ell_x \geq r^{-1}$ so that $x + \ell_x$ is the first limit point of $\exc_R$ encountered to the right of $x$.  Let $s \in\Z^+$ be a (large) integer to be determined later.  
For each $x \in \cM$, set $L_x = \big[x - s r^{-1} ,\, x - sr^{-1}+ \ell_x - r^{-1} \big]$ if $r^{-1} \leq \ell_x < (s+1)r^{-1}$; if $\ell_x \geq (s+1)r^{-1}$, take $L_x = \big[x-sr^{-1},x\big]$.
Define the {\bf domain of attraction} $\cD_x$ of $x \in \cM$ by $\cD_x = L_x \cap \exc_R$. See Figure~\ref{f.DOA}.
\end{definition}

Note that Lemma~\ref{l.IIC'arrival} has a  straightforward analogue for $\Pb{\, \cdot \md 0\in\cM}$, and we have $\P'=\Pbo{\no}{ \, \cdot \md N \geq 1/r } = \Pb{\, \cdot \md 0\in\cM}$. We now define the measure $\P''$ on dynamical configurations on $B_R$ that will be our main tool for understanding $\P'$.

\begin{definition}
Define the law $\P''$ so that, for any c\`agl\`ad dynamical percolation configuration $\omega$ satisfying $0\in \cM$,
$$
\frac{{\rm d } \P''}{{\rm d } \P'}(\omega) = Z^{-1} \vert \cD_0 \vert,
$$
where $|\cdot|$ is Lebesgue measure, and
$Z > 0$ is a normalization chosen to ensure that $\P''$ is indeed a probability measure.
\end{definition}

\begin{lemma}\label{l.rnderiv}
Let $\tilde{\P}$ denote the following dynamical process. Consider c\`agl\`ad dynamical percolation $\big\{ \omega_t: t \in \R \big\}$ in $B_R$ with $\omega_0$ distributed as $\iic_R$, and with the update decisions made independently of $\omega_0$. Condition this process on the event that $\exc_R \cap \big(sr^{-1},(s+1)r^{-1}\big) = \emptyset$. Let $\gamma \in [0,sr^{-1}]$ be given by $\gamma = \sup \{ \exc_R \cap [0,sr^{-1})\}$. Now set $\tilde{\P}$ equal to the conditional law of $\omega(\gamma + \cdot)$. Then $\tilde{\P} = \P''$.
\end{lemma}

\noindent{\bf Proof.}
Under dynamical percolation on $B_R$, let $\cJ$ denote the set of times $j \in \exc_R$ such that  $\big(j+sr^{-1},j+(s+1)r^{-1}\big) \cap \exc_R = \emptyset$.
Consider the map $A: \cJ \longrightarrow \cM$ such that, for each $j \in \cJ$, $A(j)$ is the largest element of $\cM$ preceding $j + sr^{-1}$. Note that $j \in \exc_R$ implies that $j \leq A(j) \leq j + sr^{-1}$. Note further that, for each $m \in \cM$, we have $A^{-1}(m) = \cD_m$. See Figure~\ref{f.DOA}.

Consider now an experiment in which, for $x > 0$,
dynamical percolation is sampled conditionally on $\cJ \cap [0,x] \not= \emptyset$, and an element $\chi \in \cJ \cap [0,x]$ is chosen with the conditional law of normalized Lebesgue measure on this set. Note that, by $\lim_{x\to\infty} \Ps{\cJ \cap [0,x] \not= \emptyset} = 1$, the law of $\omega_{A(\chi) + \cdot}$ (using the randomness in both $\omega$ and $\chi$ has the limit $\tilde{\P}$ as $x \to \infty$. However, from the previous paragraph we also know that $\omega_{A(\chi) + \cdot}$ has a weak limit whose Radon-Nikodym derivative with respect to dynamical percolation given $0 \in \cM$ is $\vert \cD_0 \vert$ up to normalization.\qed 

\begin{lemma}[Typical events of $\P''$ will appear in $\P'$]\label{l.rnd}
The Radon-Nikodym derivative $\frac{{\rm d } \P''}{{\rm d } \P'}$ has a second moment that is bounded above by some $B<\infty$ which might depend on the parameter $s$ but not on $R$. Consequently, $\P'(\cA)\geq \P''(\cA)^2/B$ for any event $\cA$.
\end{lemma}

 \noindent{\bf Proof.} The claim regarding the Radon-Nikodym derivative follows directly from Lemma~\ref{l.domofatt} below. The second claim then follows by Cauchy-Schwarz: 
$$
 \P''(\cA)=\int \1_\cA\, {\rm d } \P'' = \int \1_\cA \, \frac{{\rm d }\P''}{{\rm d }\P'} \, {\rm d }\P' \leq  \sqrt{\int \1_\cA^2 \, {\rm d }\P'} \, \sqrt{\int \left(\frac{{\rm d }\P''}{{\rm d }\P'}\right)^2 {\rm d }\P'} \leq \sqrt{\P'(\cA)}\, \sqrt{B}\,,
$$
as desired.
\qed

\begin{lemma}\label{l.domofatt}
Let $m_{R,r}$ denote the conditional mean under dynamical percolation of $|\exc_R \cap \big(0,r^{-1}\big)|$ given that this intersection is non-empty. 
Consider dynamical percolation $\mathbb{P}$ on $B_R$ conditionally on $0 \in \cM$.
Then the Lebesgue measure of the domain of attraction of the origin satisfies 
\begin{equation}\label{eqdoaone}
\Pb{ | \cD_0 | \geq c \, m_{R,r} \md 0 \in \cM } \geq c s^{-2} 
\end{equation}
and
\begin{equation}\label{eqdoatwo}
\Eb{ | \cD_0 |^2  \md 0 \in \cM } \leq C s^2 m_{R,r}^2\,, 
\end{equation}
for constants $C > c > 0$ which do not depend on $r,R$ or $s$.
\end{lemma}

Before starting the proof of Lemma~\ref{l.domofatt}, we need to verify a basic decorrelation result. In light of Lemma~\ref{l.rnderiv} (describing $\P''$ as $\tilde\P$), it is far from surprising that this result will be crucial in understanding the measures $\P''$ and $\P'$.

\begin{lemma}[Ensuring an empty interval]\label{l.temp}
There exists a large $s \in \Z^+$ and a small $c > 0$ such that, for each $r\in \Z^+$ and $R > R_0(r)$, the 
probability that 
dynamical percolation with initial condition $\omega_0$ distributed according to $\iic_R$
satisfies $\exc_R \cap \big( s\,r^{-1},(s+1)\,r^{-1} \big) = \emptyset$ exceeds $c$.  
\end{lemma}

An important element of the proof of Lemma \ref{l.temp} is the following claim. It is slightly more convenient to reverse time once again, just for this claim. Recall that $\rho(r)=\inf\{s: s^2\alpha_4(s)\geq r\}$, and keep in mind that its magnitude is known to be $r^{4/3+o(1)}$ for percolation on the faces of $\hexlatt$ and to lie between $C^{-1}\,r^{1+\eta}$ and  $C\,r^{1/\eta}$  for some $\eta \in (0,1)$ and  $0 < C < \infty$ for bond percolation on $\Z^2$.

\begin{lemma}\label{l.tempprep}
There exists $c > 0$ such that the following holds, independently of $r \in \N$.  
Let $\cN$ denote the event that at no time in the interval $[-r^{-1},0]$ is there an open crossing of the annulus $A_{\rho(r),\,2\rho(r)}$. For $s  > 0$, let $\cY_s$ denote the event that an open crossing of $A_{\rho(r),\,2\rho(r)}$ exists at time $s r^{-1}$. Then, for all large enough $s > 0$ (without dependence on $r$), we have $\Pb{\cN \cap \cY_s} \geq c$.
\end{lemma}

\proof By considering a coupling in which dynamical updates lead always to the closure of hexagons, we know that $\Pb{\cN} \geq c$ by~(\ref{e.winmed}), Kesten's result on the near-critical window. Let $\cN_0$ denote the time-$0$ static event that the conditional probability of $\cN$ given the time $0$ configuration is at least $c$. We have that $\Ps{\cN_0} \geq c$ by adjusting the value of $c > 0$.
Note then that, denoting by $f$ and $g$ the $\pm 1$-indicator functions of  $\cN_0$ and $\cY_0$, and by $\what{f}$ and $\what{g}$ their Fourier series, the basic relation~(\ref{e.fourierrelation}) yields
$$
 \Pb{\cN_0 \cap \cY_s} - \Pb{\cN_0} \Pb{\cY_s}  =  \sum_{S \not= \emptyset} \what{f}(S) \what{g}(S) \exp \big\{ - s r^{-1} \vert S \vert \big\} \, .
$$
We apply Cauchy-Schwarz to bound above the absolute value of the right-hand side.
Then, the  basic relation~(\ref{e.fourierrelation}) and the decorrelation estimate~(\ref{e.fulldecor}) applied to $g$ 
give the following bound on the resulting term:
\begin{multline*}
\left( \sum_{S \not= \emptyset} \what{f}^2(S) \right)^{1/2}  \left( \sum_{S \not= \emptyset}  \what{g}^2(S) \exp \big\{ - 2s r^{-1} \vert S \vert \big\} \right)^{1/2} \\
\leq \eps_1 + \left( \sum_{\vert S \vert \geq \eps_2 r} \what{g}^2(S)  \exp \big\{ -  2\eps_2 s \big\} \right)^{1/2}  
\leq \eps_1 + \exp \big\{ -  \eps_2 s \big\} \, ,                        
\end{multline*}
where $\eps_1$ depends on the choice of cutoff $\eps_2 > 0$ and may be chosen so that $\eps_1 \to 0$ as $\eps_2 \to 0$. Noting that $\Pb{\cN_0} \Pb{\cY_s} \geq c_1>0$, we see that $\Pb{\cN_0 \cap \cY_s} \geq c_1/2$ 
by making a suitable choice of $\eps_1$, $\eps_2$ and $s$. Note that $\Pb{\cN \cap \cY_s} \geq c \, \Pb{\cN_0 \cap \cY_s}$ because $\cN$ and $\cY_s$ are conditionally independent given the time-$0$ configuration. This completes the proof. \qed
\medskip

The next lemma relates the restriction of $\iic$ to a dyadic annulus to the percolation configuration in the annulus obtained by conditioning on an open crossing between the annulus' boundaries.

\begin{lemma}[Localizing the $\IIC$ conditioning]\label{l.anniic}
Let $\P_r^R$ denote the law of critical percolation in $A_{r,R}$ given that $r \llra R$, for $0\leq r <  R \leq\infty$ (where the conditional law $\Ps{ \cdot \md r \lra \infty}$ on $B_{r}^c$ is obtained as a weak limit of $\Ps{\cdot \md r \lra R}$ as $R \to \infty$, constructed by~\cite{KestenIIC}). Then, for each $\eps > 0$ there exists $\delta > 0$ such that if $\cA \in \sigma \{ A_{R,2R} \}$ (i.e., an event measurable in the annulus), then $\P_R^{2R}(\cA) \geq \eps$ implies that $\P_a^b(\cA) \geq \delta$, for all $0\leq a\leq R/2$ and $4R\leq b\leq\infty$ (in particular, for $\IIC=\P_0^\infty$).
\end{lemma}

\proof  For $\zeta$ a configuration in $A_{R,2R}$ such that $R \lra 2R$, let $W_{a,R,b}(\zeta)$ denote the conditional probability that $a \lra b$ given the occurrence of the events $\omega \big\vert_{A_{R,2R}} = \zeta$, $a \lra R$ and $2R \lra b$. We  will argue that for each $\eps > 0$ there exists $\delta > 0$ such that, for all large enough $R \in \N$ and $a,b\in\N$ with $0\leq a\leq R/2$ and $2R\leq b\leq \infty$, 
\begin{equation}\label{e.iicr}
\PB{ W_{a,R,b} \leq \delta \md  R \llra 2R } \leq \eps \, .
\end{equation}
This easily implies the lemma, as follows. Note that 
$$
\frac{\mathrm{d}\P_a^b}{\mathrm{d} \P_R^{2R}} (\zeta) = Z_{a,R,b}^{-1} \, W_{a,R,b}(\zeta) \, ,
$$
where $Z_{a,R,b} = \Ps{a \lra b \md a \lra R, R \lra 2R, 2R \lra b} \leq 1$. 
Given $\eps > 0$, choose by means of (\ref{e.iicr}) an $\eps' > 0$ such that 
 $\P_R^{2R} \big( W_{a,R,b} \leq \eps'  \big) \leq \eps/2$ for each $R \in \N$. Thus, if  
 $\cA \in \sigma \{ A_{R,2R} \}$ satisfies $\P_R^{2R}(\cA) \geq \eps$, then 
 $$\P_a^b(\cA) =  Z_{a,R,b}^{-1} \int_{\cA} W_{a,R,b}(\omega) \, \mathrm{d} \P_R^{2R}(\omega) \geq \eps' \eps/2\,,$$ 
where the inequality follows from restricting the integral to that part of $\cA$ on which $W_{a,R,b} > \eps'$. Hence the lemma holds with the choice $\delta = \eps' \eps/2$.

To prove~(\ref{e.iicr}), we introduce the function $W_{R}^\eps(\zeta)$ on configurations $\zeta$ in $A_{R,2R}$, for $R \in \N$ and $\eps \in (0,1/2)$, which is the conditional probability of $R(1-\eps) \llra 2R(1+\eps)$ under critical percolation given that  $\omega \big\vert_{A_{R,2R}}  = \zeta$. 

\begin{lemma}\label{l.mbrmbrtwo}
For each $\eps \in (0,1/2)$, there exists a constant $c = c_\eps > 0$ such that, for each $R,a,b \in \N$ as before and for all configurations $\zeta$ in $A_{R,2R}$, we have $W_{a,R,b}(\zeta) \geq c \, W_{R}^{\eps}(\zeta)$.  
\end{lemma}

\proof
Let $p_1$ denote the probability under critical percolation that there exists an open surrounding circuit in the annulus $A_{R(1-\eps),R}$, and let $p_2$ denote the corresponding probability for the annulus $A_{2R,2R(1 + \eps)}$. Note that $p_1,p_2 \geq c_\eps > 0$ for all $R$ by a simple application of RSW.
We claim that 
\begin{equation}\label{e.mrclaim}
W_{a,R,b}(\zeta) \geq p_1 p_2 W_{R}^\eps(\zeta) \, .
\end{equation}
Indeed, consider the conditioning appearing in the definition of $W_{a,R,b}(\zeta)$: under the conditional law, the configuration in $A_{R,2R}^c$ stochastically dominates critical percolation, and thus open surrounding circuits appear in the annuli $A_{R(1-\eps),R}$ and $A_{2R,2R(1+\eps)}$ with probability at least $p_1 p_2$; the presence of such circuits being an increasing event, the conditional law further conditioned on the presence of such circuits has probability at least $W_{R}^{\eps}(\zeta)$ of realizing
$R(1-\eps) \llra 2R(1+\eps)$. However, the event $R(1-\eps) \llra 2R(1+\eps)$ and the presence of the two surrounding circuits is enough, alongside the conditions met under the conditional law, to ensure that $0 \llra \infty$. In summary, we obtain~(\ref{e.mrclaim}); applying  $p_1 p_2 \geq c_\eps^2$ completes the proof. \qed




\begin{lemma}\label{l.mbreta}
For each $\delta > 0$, there exists $\eps_0 > 0$ such that, for all large enough $R \in \N$ and all $\eps \in (0,\eps_0)$,
 $$
\Pb{ W_R^{\eps} \leq 1 - \delta \md R \llra 2R} \leq \delta \, .
 $$
\end{lemma}

\proof 
Note that for any $\delta_1>0$ there is an $\eps_1>0$ such that, for all $\eps<\eps_1$, 
\begin{align*}
\Eb{W_R^\eps \md R\llra 2R} &= \Pb{R(1-\eps) \llra 2R(1+\eps) \md R\llra 2R}\\
& = 1-\frac{\Pb{R\llra 2R\,, \text{ but } R(1-\eps) \hskip 0.2 cm\not \hskip -0.2cm \llra 2R(1+\eps)}}{\Ps{R\llra 2R}}\\
&\geq 1-\delta_1\,,
\end{align*}
because the event $\big\{R\llra 2R\,, \text{ but } R(1-\eps) \hskip 0.2 cm\not \hskip -0.2cm \llra 2R(1+\eps)\big\}$ implies that there are three arms from one side of the annulus $A(R(1-\eps), 2R(1+\eps))$, from radius about $\eps R$ to radius about $R$, and this event has probability of order $\eps$, the 3-arm half-plane probability being of order $\eps^{2}$. See \cite[first exercise sheet]{WWperc}.

From this bound, applying Markov's inequality to $1-W_R^\eps$, we get that 
$$
\Pb{ 1-W_R^{\eps} \geq \sqrt{\delta_1} \md R \llra 2R} \leq \sqrt{\delta_1} \,,
$$
which implies the lemma immediately.
\qed
\medskip

Now note that (\ref{e.iicr}) follows from Lemmas~\ref{l.mbrmbrtwo} and~\ref{l.mbreta} immediately. This completes the proof of Lemma~\ref{l.anniic} (localizing the $\IIC$ conditioning) for large enough $R\in\N$; on the other hand, for $R$ bounded, the lemma is trivial.
\qed

\proofof{Lemma \ref{l.temp}} (Ensuring an empty interval.)
Let $\cC \in \sigma \big\{ \annsmall \big\}$ denote the static event consisting of configurations $\zeta$ satisfying $\rho(r) \leftrightarrow 2\rho(r)$ and such that
$$
\PB{  \rho(r) \leftrightarrow 2\rho(r) \, \, \textrm{at no time in $[sr^{-1},(s+1)r^{-1}]$} \md \omega_0 = \zeta } \geq c\,.
$$
By considering the process $\omega\big(sr^{-1} - \cdot \big)$ in Lemma \ref{l.tempprep}, we see that 
$$
\P\Big( \rho(r) \lra 2\rho(r) \textrm{ at time }0 \, ,  \ 
\rho(r) \lra 2\rho(r) \textrm{ at no time in }[sr^{-1},(s+1)r^{-1}]\Big) \geq c\,;
$$
in the notation of the statement of Lemma \ref{l.anniic}, we see that $\P_{\rho(r)}^{2\rho(r)}(\cC) \geq c$ by reducing the value of $c > 0$. By Lemma \ref{l.anniic}, we infer that for some $\delta > 0$ and for $R \geq 4 \rho(r)$, $\iic_R(\cC) > \delta$, as required for the statement of Lemma \ref{l.temp}. \qed

\proofof{Lemma~\ref{l.domofatt}} We start by a simple corollary of Lemma~\ref{l.temp} concerning the density of markers.

\begin{definition}
Let $\{I_i = \big(i/r , (i+1)/r\big): i \in \N\}$ enumerate the consecutive intervals of length $r^{-1}$ rightwards from the origin. Call any such interval {\bf active} if it has non-empty intersection with $\exc_R$. 
For any $i \in \N$, call $I_i$ {\bf promising} if $I_i$ is an active interval with the property that $\cM$ intersects $\cup_{i \leq j \leq i + s} I_j$. 
\end{definition}

\begin{lemma}\label{lemprop}
There exists $c > 0$, independent of $R$, such that the conditional probability under dynamical percolation $B_R$ given that $I_0$ is active that $I_0$ is promising is at least $c$.
\end{lemma}

\proof Let $\P_0$ denote dynamical percolation on $(0,1/r)$ weighted according to the size $\big\vert \exc_R \cap (0,1/r) \big\vert$; under $\P_0$, define $\tau$ to be an element of 
 $\exc_R \cap (0,1/r)$ with conditional law given by normalized Lebesgue measure on this set. Under $\P_0$,
 the law of dynamical percolation at times $\tau+t$, $t\geq 0$ is, by Lemma~\ref{l.IICra}, dynamical percolation started from $\iic_R$. By Lemma \ref{l.temp}, the conditional probability that $\big(\tau + sr^{-1},\tau + (s+1)r^{-1}\big) \cap \exc_R = \emptyset$ exceeds some $R$-independent constant $c > 0$. Whenever this  disjointness condition is satisfied, there exists an element of $\cM$ somewhere in the interval between $\tau$ and  $\tau + sr^{-1}$, and thus in the interval $\big(0,(s+1)r^{-1}\big)$.

We learn that the $\P_0$-probability that $I_0$ is promising exceeds an $R$-independent constant $c > 0$. Lemma \ref{lemprop} will follow once we establish this assertion for dynamical percolation conditioned on the interval $I_i$ being active, a measure we label $\P_1$.  To make this reduction, it is enough to argue that $\tfrac{{\rm d} \P_0}{{\rm d} \P_1}$ has a bounded second moment, in light of the proof of Lemma~\ref{l.rnd}, with the roles of $\P''$ and $\P'$ being played by $\P_0$ and $\P_1$.    
By Lemma~\ref{l.cond2ndM}, there exists $C > 0$ such that, for all $R > 0$, 
\begin{eqnarray*}
 \int \Big(\frac{{\rm d} \P_0}{{\rm d} \P_1}\Big)^2 {\rm d} \P_1 & = &
  \E \Big( \big\vert \exc_R \cap (0,r^{-1}) \big\vert^2 \, \Big\vert \, \exc_R \cap (0,r^{-1}) \not= \emptyset  \Big) \nonumber \\
 & \leq & C  \bigg( \E \Big( \big\vert \exc_R \cap (0,r^{-1}) \big\vert \, \Big\vert \, \exc_R \cap (0,r^{-1}) \not= \emptyset  \Big)   \bigg)^2 \, .
\end{eqnarray*}
This completes the proof of Lemma \ref{lemprop}. \qed
\medskip

We can now prove (\ref{eqdoatwo}). Let $\big\{ m_i: i \in \N^+ \big\}$ enumerate the elements of $\cM \cap (0,\infty)$ in increasing order. By ergodicity, we have almost surely that
\begin{equation}\label{eqdzero}
\Eb{ \vert \cD_0 \vert^2  \md 0 \in \cM } 
 = \lim_n n^{-1} \sum_{i=2}^n \vert \cD_{m_i} \vert^2 \, ,
\end{equation}
where the term with index $i=1$ has been harmlessly omitted for later notational convenience.
Let $\{ \lambda_i\}$ (or $\{\alpha_i\}$)  enumerate the indices $i \in \N^+$ of promising (or active) intervals $I_i$ in increasing order.  
For $i \geq 1$, consider the consecutive intervals $I_j$ beginning the interval after that containing $m_i$ and stopping at the one containing $m_{i+1}$. Among these, there are at most $s+1$ promising intervals, and $\cD_{m_{i+1}}$ is contained in the union of these promising intervals. Therefore, 
 $\sum_{i=2}^n \vert \cD_{m_i} \vert^2 \leq (s+1) \sum_{i=2}^{\lambda_{(s+1)n}} \vert I_i \cap \exc_R \vert^2$. By Lemma \ref{lemprop} and the ergodicity Lemma~\ref{l.erg}, $\lambda_n \leq 2 c^{-1} \alpha_n$ for all large enough $n$. Hence, $\sum_{i=2}^{n} \vert \cD_{m_i} \vert^2 \leq (s+1) \sum_{i=1}^{2 c^{-1} \alpha_{(s+1)n}} \vert I_i \cap \exc_R \vert^2$. By ergodicity again, this upper bound behaves like
$$
2 c^{-1} (s+1)^2 \, n \, \Eb{ \vert \exc_R \cap (0,1/r ) \vert^2  \md \exc_R \cap (0,1/r ) \not= \emptyset   }  \, \big(1+o(1)\big) 
$$ 
as $n\to\infty$.  Applying Lemma~\ref{l.cond2ndM} to  $| \exc_R \cap (0,1/r) |$ (which is just a scaled version of $\muA_R(0,1/r)$) and using (\ref{eqdzero}), we obtain (\ref{eqdoatwo}).

To prove (\ref{eqdoaone}),
in light of (\ref{eqdoatwo}), the Paley-Zygmund second moment method says that it suffices to verify that, for some $c > 0$ and all $R,r,s \in \N_+$, 
\begin{equation}\label{e.doathr}
 \Eb{ \vert \cD_0 \vert \md 0 \in \cM} \geq c \, m_{R,r}\,. 
\end{equation}
We now verify this inequality. Let $\rho = \lim_n n^{-1} \big\vert \cM \cap (0,n) \big\vert$ denote the mean number of markers in $[0,1]$, or, alternatively,  
$\rho = \Es{\vert \cM \cap (0,1) \vert}$.

We claim the following.

\begin{lemma}\label{l.dzero} 
Recall that $\cJ$ denotes the set of times $j$ such that the event $0 \leftrightarrow R$ occurs at time $j$ and at no time in the interval $\big(j+sr^{-1}, j+(s+1)r^{-1}\big)$. Then
$\rho\,  \Eb{ \vert \cD_0 \vert \md 0 \in \cM} = \Eb{| \cJ \cap [0,1] |}$.
\end{lemma}

\noindent{\bf Proof.}
Recall that the subset $\cJ$ of $\exc_R$ is partitioned
into disjoint classes given by domains of attraction $\cD_m$ and thus indexed
by the set of markers $m \in \cM$.

The quantity $\rho\, \Eb{\vert \cD_0 \vert \md 0 \in \cM}$ is thus the mean Lebesgue measure of the union of the domains of
attractions indexed by markers lying in a given unit interval. By the above partition and ergodicity, we arrive at the statement of Lemma~\ref{l.dzero}. \qed
\medskip

By translation invariance, $\Es{\vert \cJ \cap [0,1] \vert} = \lim_n n^{-1} \Es{\vert \cJ \cap [0,n] \vert}$; by Lemmas \ref{l.IICrq} and \ref{l.temp}, there exists $c > 0$ such that, for $n$ sufficiently high, 
$$
n^{-1} \Es{\vert \cJ \cap [0,n] \vert} \geq c n^{-1} \Es{\vert \exc_R \cap [0,n] \vert} \, . 
$$
By translation invariance again, $n^{-1} \Es{\vert \exc_R \cap [0,n] \vert} = \Es{\vert \exc_R \cap [0,1] \vert}$ which may be written 
$r m_{R,r} \Ps{I_0 \, \textrm{is active}}$. To summarise the derivation of (\ref{e.doathr}) thus far, the preceding inequality and Lemma \ref{l.dzero} yield
\begin{equation}\label{e.doathrsum}
\rho\,  \Eb{\vert \cD_0 \vert \md 0 \in \cM} \geq 
c \, r \, m_{R,r} \, \Ps{I_0 \, \textrm{is active}}.
\end{equation}
We will show that
\begin{equation}\label{e.doathrnew}
\rho\,  \leq r \, \Ps{I_0 \, \textrm{is active}}\,;
\end{equation}
note then that (\ref{e.doathrsum}) and (\ref{e.doathrnew}) yield (\ref{e.doathr}).

To verify (\ref{e.doathrnew}), recall that a marker is by definition an element of $\exc_R$ bordered on the right by an interval of length $r^{-1}$ having no intersection with $\exc_R$. 
Thus,
each marker lies in an active interval, and no active
interval contains more than one marker. This implies that the mean rate $\rho$ of markers is at most the mean number of active 
intervals in a given unit interval, a quantity which may be expressed as $r \, \Ps{I_0 \, \textrm{is active}}$. This
verifies (\ref{e.doathrnew}).  This completes the derivation of (\ref{e.doathr}) and thus of (\ref{eqdoaone}), which concludes the proofs of Lemmas~\ref{l.domofatt} and~\ref{l.rnd} on the Radon-Nikodym derivative $\tfrac{{\rm d}\P''}{{\rm d}\P'}$. \qed
\medskip

We are now ready to address the main goal of this subsection. Recall the notion of $\good$ from Definition~\ref{d.good}.

\begin{proposition}[$\P'$ is well behaved]\label{lemgoodm}
There exists $c > 0$ such that, for any $r>r_0$ and $R>R_0(r)$,
$$
\Pb{ \omega \in \good \md 0 \in \cM} \geq c \, .
$$
\end{proposition}

In the proof, we will use the following notion and claim. 

\begin{definition}\label{defvgood}
Fix $\eps>0$ and $r\in\N$ as in Definition~\ref{d.Fine}, and let $R \in \N$ satisfy $R \geq  r^{1 + 2\eps}$.
We say that  a dynamical configuration $\omega$ in $B_R$ is $\omega \in \vgood$ if the following conditions are satisfied:
\begin{itemize}
 \item $0 \build\longleftrightarrow_{}^{\omega_0} R$;
 \item For each $t \in  [0,(s+1)r^{-1}]$, the inner and outer boundaries of the annulus $A_{r^{1 + \eps},r^{1 + 2\eps}}$ are separated by an $\omega_t$-open circuit; 
  \item  For each $t \in  [0,(s+1)r^{-1}]$, the circuit $\Gamma_r$ exists and satisfies $\Gamma_r \subseteq B_{r^{1 + \eps}}$ in $\omega_t$;
 \item $|\Piv_{0\lra\Gamma_r}| \leq |\Piv_{0 \lra r^{1 + \eps}}(\omega_t)| \le r^{2(1+2\eps)}\alpha_4(r^{1+2\eps})$ for all such $t$.
 \item For each $t \in [0,(s+1)r^{-1}]$, $0 \build\longleftrightarrow_{}^{\omega_t}   r^{1 + 2\epsilon}$. 
\end{itemize}
\end{definition}

\begin{lemma}\label{l.verygood}
For any $\delta > 0$, there exists $r_0 \in \N$ such that, for all $r \geq r_0$ and $R \geq  r^{1 + 2\eps}$,
\begin{equation}\label{eqomdelta}
\Pb{ \omega \in \vgood \md 0 \build\longleftrightarrow_{}^{\omega_0} R  } \geq 1 - \delta \, .
\end{equation}
\end{lemma}

\noindent{\bf Proof.} Let $\P_1$ denote dynamical percolation in $B_{r^{1 + 2\epsilon}}$ with $\omega_0$ having the distribution $\Pb{ \, \cdot  \md  0 \leftrightarrow R }$, and with conditionally independent updates at rate one. Let $\P_2$ denote the asymmetric dynamical process in $B_{r^{1 + 2 \epsilon}}$ with the same initial distribution as in $\P_1$, but with the updates always leading to the closure of hexagons. As usual, we form the obvious coupling $\coup$ of $\P_1$ and $\P_2$ such that the first marginal dominates the second for all $t \geq 0$.

In this new notation, the statement of the lemma is equivalent to: 
for any $\delta > 0$, there exists $r_0 \in \N$ such that, for all $r \geq r_0$ and $R \geq  r^{1 + 2\eps}$,
\begin{equation}\label{eqomdeltanew}
\P_1( \omega \in \vgood' ) \geq 1 - \delta \, ,
\end{equation}
where $\vgood'$ is given by the second and later conditions defining $\vgood$.  We claim that, to show (\ref{eqomdeltanew}), it is enough that
\begin{equation}\label{eqnewtwo}
\P_2( \omega \in \vgood' ) \geq 1 - \delta\,.
\end{equation}

To see that (\ref{eqnewtwo}) is enough for (\ref{eqomdeltanew}), note that,
under the coupling $\coup$, it is clear that if the second ($\P_2$-distributed) marginal satisfies the second, third  and fifth conditions of Definition~\ref{defvgood},  then so does the first ($\P_1$-distributed) marginal, because these conditions are monotone. In regard to the fourth condition, write $\Piv$ for $\Piv_{0 \lra r^{1 + \eps}}$. Note that if $\omega_1$ and $\omega_2$ are two configurations in 
$B_{r^{1+\epsilon}}$ such that $\omega_1 \geq \omega_2$ and  $0 \leftrightarrow r^{1+\epsilon}$ under $\omega_2$, then $\Piv(\omega_1) \subseteq \Piv(\omega_2)$: 
indeed, were a hexagon $h$ in 
$B_{r^{1+\epsilon}}$ to satisfy $h \in \Piv(\omega_1)\setminus \Piv(\omega_2)$, then its closure would disable $0 \leftrightarrow r^{1+\epsilon}$ in $\omega_1$ but not in $\omega_2$, a circumstance which stochastic domination prevents. That is, whenever $0  \build\longleftrightarrow_{}^{\omega_t} r^{1 + \epsilon}$ occurs under $\P_2$, we have that  $\vert \Piv(\omega^1_t) \vert \leq \vert \Piv(\omega^2_t) \vert$ (where $\omega^1$ and $\omega^2$ denote the $\P_1$ and $\P_2$ marginals), and thus (\ref{eqnewtwo}) implies (\ref{eqomdeltanew}) and hence (\ref{eqomdelta}).

It remains to verify (\ref{eqnewtwo}). We start with a simple lemma.

\begin{lemma}\label{l.iiccomp}
Let $\P^\downarrow_s$ denote asymmetric dynamical percolation $\{\omega_t : t\geq 0\}$ with $\omega_0$ having the distribution $\IIC_s$, then closing hexagons at rate one. Then,
restricted to the ball $B_r$, the Radon-Nikodym derivative $\frac{\d \P^\downarrow_R}{\d \P^\downarrow_r}\big(\omega[0,t]^{B_r}\big)$ is bounded from above uniformly in $r$, $R \geq r$, $t\geq 0$, and all dynamical configurations $\omega[0,t]^{B_r}$.
\end{lemma}

\proof
We claim that 
\begin{equation}\label{e.iicrone}
\frac{\d \P^\downarrow_R}{\d \P^\downarrow_r}\big(\omega[0,t]^{B_r}\big) \leq \frac{\Ps{0 \lra r} \Ps{r+1 \lra R}}{\Ps{0 \lra R}}  \, ,
\end{equation}
with the right hand side understood simply in static critical percolation. From this, the lemma follows by quasi-multiplicativity.
For the claim, note that the Radon-Nikodym derivatives with respect to asymmetric dynamical percolation $\P^\downarrow$ started from criticality, restricted to $B_r$, can be written as 
$$\frac{\d \P^\downarrow_s}{\d \P^\downarrow}\big(\omega[0,t]^{B_r}\big) = \frac{\P^\downarrow\big( 0 \lra s \text{ in }\omega_0 \,\big|\, \omega[0,t]^{B_r}\big)}{\P^\downarrow(0 \lra s \text{ in }\omega_0 )}\,,$$ 
for any $s\geq r$; in particular, for $s\in\{r,R\}$. On the other hand, 
$$\P^\downarrow\big(0 \lra R \text{ in }\omega_0 \,\big|\, \omega[0,t]^{B_r}\big) \leq \P^\downarrow\big(0 \lra r \text{ in }\omega_0 \,\big|\,  \omega[0,t]^{B_r}\big) \, \P^\downarrow(r+1 \lra R \text{ in }\omega_0 )\,.$$
Since the distribution of $\omega_0$ under $\P^\downarrow$ is simply critical percolation, from the last two displays follows~(\ref{e.iicrone}).
\qed
\medskip

\proofof{(\ref{eqnewtwo})}
Let $\P_3$ denote asymmetric dynamical percolation $\P^\downarrow_{r^{1 + 2\eps}}$ in $B_{r^{1+2\eps}}$, with the notation of the previous Lemma~\ref{l.iiccomp}. By that lemma, it is enough to verify~(\ref{eqnewtwo}) with $\P_3$ in place of $\P_2$. We are going to show that each of the four conditions defining $\vgood'$ happens with probability close to 1 if $r$ is large enough.

\begin{figure}[htbp]
\SetLabels
(.05*.48)$0$\\
(.57*.38)$0$\\
(.17*.54)$x$\\
(.82*.45)$x$\\
(.35*.27)$r^{1+\eps}$\\
(.86*.15)$r^{1+\eps}$\\
\endSetLabels
\centerline{
\AffixLabels{
\raise 0.5in\hbox{\epsfxsize=2.5 in \epsffile{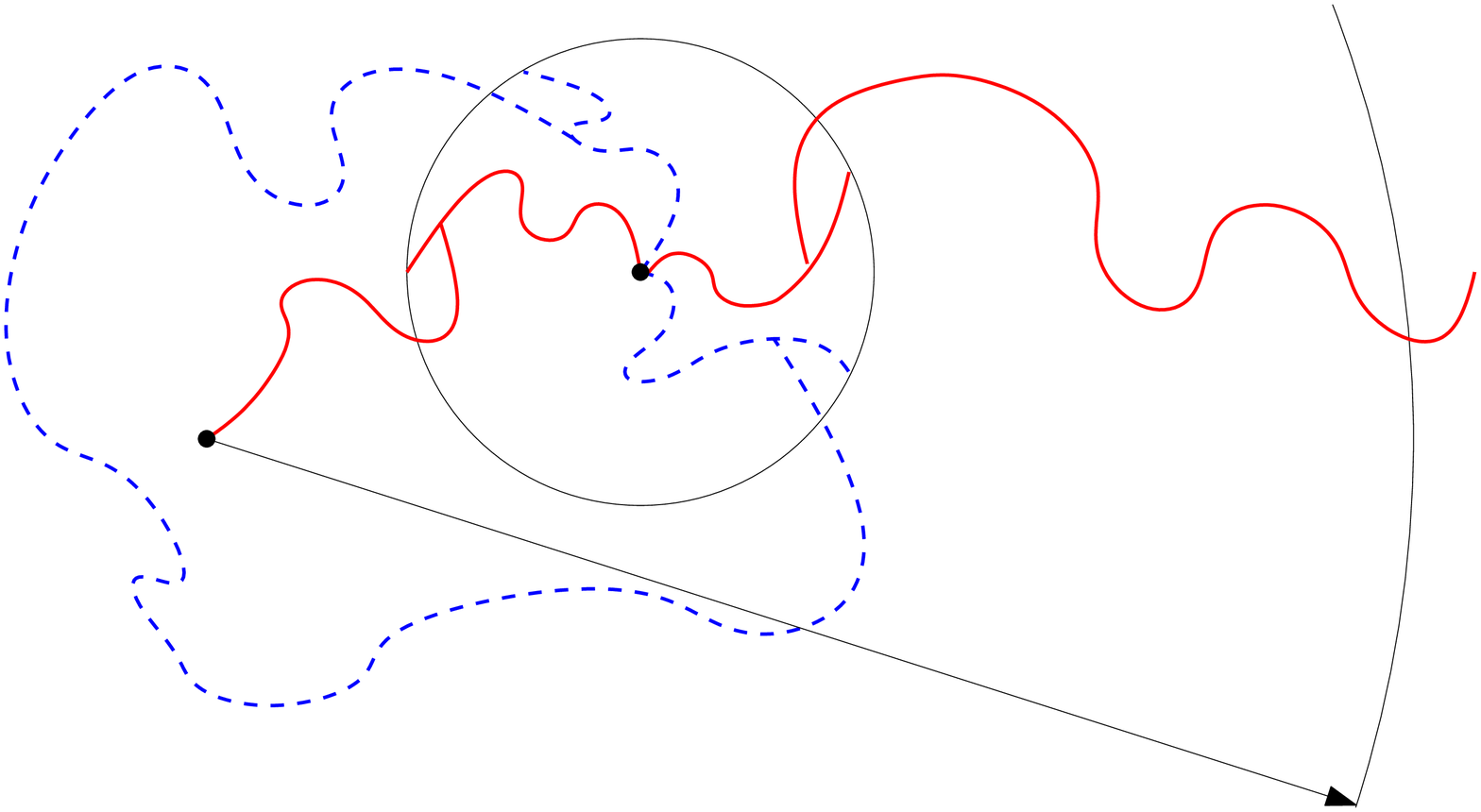}}
\hskip 0.5 in
\epsfxsize=3 in \epsffile{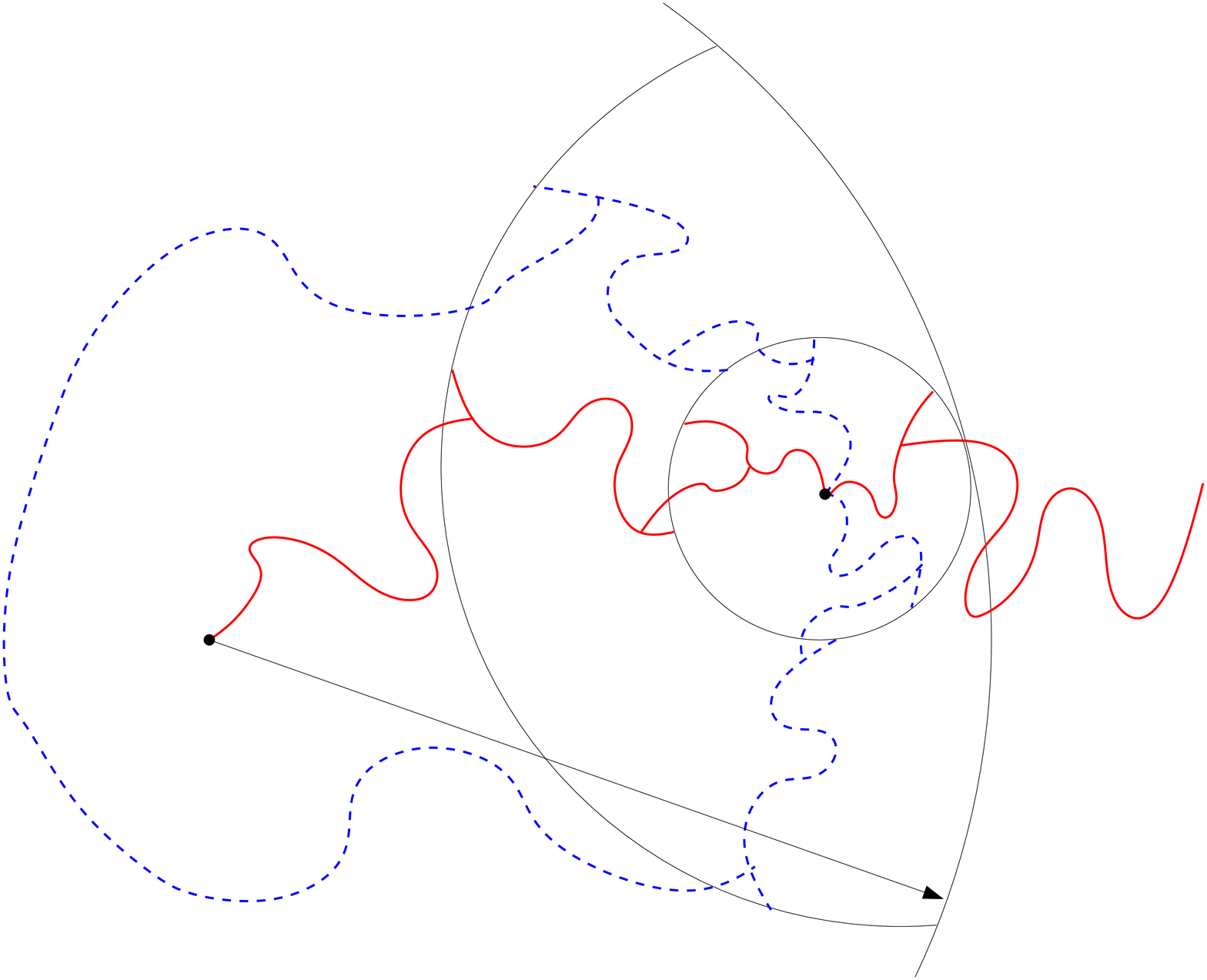}
}
}
\caption{Conditioned on $0\lra r^{1+2\eps}$, one 4-arm event (first picture) or one 4-arm and one 3-arm events (second picture) are roughly equivalent to $x$ being pivotal for $0\lra r^{1+\eps}$.}
\label{f.1armPiv}
\end{figure}

Let us first look at the four conditions at time zero. The fifth condition (that $0\lra r^{1+2\eps}$) is automatically satisfied under $\P_3$. The second and third conditions (open circuits in $A_{r^{1 + \eps},r^{1 + 2\eps}}$ and in $A_{r^{1+\eps},r}$) are satisfied with high probability in critical percolation by RSW along several scales, and also under the conditioning $0\lra r^{1 + 2\eps}$ by FKG. The fourth condition (there are not too many pivotals for $0\lra r^{1+\eps}$) follows from standard quasi-multiplicativity arguments. Namely, as illustrated on Figure~\ref{f.1armPiv}, we have
$$
\Pb{x\in\Piv_{0\lra r^{1+\eps}} \md 0\lra r^{1 + 2\eps} } \asymp \alpha_4\big(\dist(x,\p B_{r^{1+\eps}})\wedge \dist(0,x)\big)\, \alpha_3\big(\dist(x,\p B_{r^{1+\eps}}),r^{1+\eps}\big)\,,
$$
which can be summed up over the possible hexagons $x\in B_{r^{1+\eps}}$ to get 
$$\EB{|\Piv_{0\lra r^{1+\eps}}| \md 0\lra r^{1+2\eps}} = O(1)\,r^{2(1+\eps)}\alpha_4(r^{1+\eps})\,.$$ 
By quasi-multiplicativity and (\ref{e.a4}), we have 
$$
\frac{\alpha_4(r^{1+\eps})}{\alpha_4(r^{1+2\eps})} < C\, (r^\eps)^{2-\eta} \ll r^{2\eps}\,,$$
hence Markov's inequality yields 
$$
\PB{|\Piv_{0\lra r^{1+\eps}}| >  r^{2(1+2\eps)}\alpha_4(r^{1+2\eps}) \md 0\lra r^{1+2\eps}} \to 0\,,
$$ 
as $r\to \infty$, as desired.

We now have to prove that the four conditions are also satisfied with high probability at time $t=(s+1)r^{-1}$; then, by the earlier monotonicity argument, we have the result for all $t\in  [0,(s+1)r^{-1}]$, as well. 

By the exponent bound~(\ref{e.piv}) and the choice $(1+2\eps)(1-\eta)<1$ made in Definition~\ref{d.Fine} and onwards, we have that $r^{-1} \ll 1/\big(r^{2(1+2\eps)} \alpha_4(r^{1+2\eps})\big)$, as $r\to\infty$. Thus, the constant closing of hexagons for time $(s+1)r^{-1}$ keeps the system $B_{r^{1+2\eps}}$ well inside the critical window of percolation, established by Kesten, as described in~(\ref{e.winsmall}) and~(\ref{e.armstab}). Therefore, the above arguments for the second to fourth conditions of Definition~\ref{defvgood} apply verbatim. The fifth condition can be verified in a similar manner:  using (\ref{e.armstab}), 
we have 
$$
\P^\downarrow\Big(0 \build\longleftrightarrow_{}^{\omega_{(s+1)r^{-1}}} r^{1+2\eps}  \,\Big|\, 0  \build\longleftrightarrow_{}^{\omega_0} r^{1+2\eps} \Big) = 1-o(1)\,,
$$ 
as $r\to\infty$. This finishes the proof of~(\ref{eqnewtwo}) and Lemma~\ref{l.verygood}.
\qed
\medskip

\proofof{Proposition~\ref{lemgoodm}}
Whenever (\ref{eqomdelta}) holds, by Lemma \ref{l.temp} we also have that
\begin{align}
\Pb{ \omega \in \vgood \md 0  \build\longleftrightarrow_{}^{\omega_0} R ,\ \exc_R \cap (sr^{-1},(s+1)r^{-1}) = \emptyset}& \nonumber\\
&\hskip -9 cm \geq~ \frac{ \Pb{ \exc_R \cap (sr^{-1},(s+1)r^{-1}) = \emptyset \md 0 \in \exc_R } - \Pb{\omega \not\in \vgood \md 0\in\exc_R}}{\Pb{ \exc_R \cap (sr^{-1},(s+1)r^{-1}) = \emptyset \md 0 \in \exc_R }}\nonumber \\
&\hskip -9 cm \geq~ 1 - \tfrac{\delta}{c}\,.\label{eqomdeltatwo}
\end{align}
Note that if a realization of dynamical percolation in $B_R$ realizes $\vgood$, then the process $\omega(\gamma + \cdot)$ identified in Lemma \ref{l.rnderiv} realizes $\good$. By Lemma \ref{l.rnderiv} and (\ref{eqomdeltatwo}), we find then that $\P''(\omega \in \good) \ge 1 - \delta/c$. Now Lemma~\ref{l.rnd} implies that an appropriate small choice of $\delta$ in~(\ref{eqomdelta}) forces $\P'(\omega \in \good) \geq c'$ for some absolute constant $c'>0$, concluding the proof of Proposition~\ref{lemgoodm}.\qed

\subsection{Size-biasing arguments}\label{ss.sizebias}

In this subsection, we will prove the bounds (\ref{e.WisOK}) and (\ref{e.hatWisOK}), used in the proof of Proposition~\ref{lemvnthin} at the end of Subsection~\ref{ss.skeleton}. To start with, note that 
\begin{equation}\label{e.hatWisOKold}
\Pb{\what N > 1/r } = \frac{\Es{N\,\1_{N>1/r}}}{\E N} >  c > 0\,.
\end{equation}
Indeed, by Lemma~\ref{lemsizebiasing}, the distribution of $\what{N} = \what{N}_R$ stochastically dominates that of $\fet_R$, which implies (\ref{e.hatWisOKold}) trivially.
\medskip

\noindent{\bf Deriving (\ref{e.WisOK}).} 
Our goal is to show that
\begin{equation}\label{e.WisOKold}
\Eb{N  \md N > 1/r} <  C <\infty\,,
\end{equation}
uniformly in $r$ and $R$ for which $R \geq r^{1 + 2\eps}$, since Proposition~\ref{lemgoodm} and $\good \subseteq \thin$ then imply  that, 
for such values of $R$ and $r$, 
$$
\Eb{N \md N > 1/r,\thin}  \leq \frac{\Eb{N \md N > 1/r}}{\Pb{\thin \md N > 1/r}} \leq c^{-1} \, \Eb{N \md N > 1/r} < c^{-1}C < \infty\,,
$$
which was the statement of (\ref{e.WisOK}).

Lemmas~\ref{lemsizebiasing} and~\ref{l.expdecay} imply that $\Eb{\what{N}}<C$ for some constant $C<\infty$ that is independent of $R$. This, together with the lower bound (\ref{e.hatWisOKold}), plugged into the next lemma with $X:=N$ and $t:=1/r$, implies (\ref{e.WisOKold}).

\begin{lemma}[Rough size-biasing]\label{l.RSB} 
If $X$ is a non-negative random variable, and $0<t<1$ is such that  $\Ps{\what{X} > t} > c > 0$ and $\Es{\what{X}} < C < \infty$, then $\Eb{X \md X>t} < C'<\infty$, where $C'$ depends only on $c$ and $C$, and not on $t$.
\end{lemma}

\proof 
Note that $\Eb{X \md X>t} = \Ps{\what{X} > t} \frac{ \Es{X} }{ \Ps{X>t} }$. Hence, we need to show that $\Es{X} \leq C' \,\Ps{X>t}$. We will need two ingredients for this:
\begin{itemize}
\hitem{(A)}{(A)} There exists an absolute constant $A<\infty$ such that 
$$\Eb{X \1_{t \geq X}} < A\, \Ps{X > t}\,.$$
\hitem{(B)}{(B)} For all $b>0$ there is some $K<\infty$ such that
$$\Eb{X \1_{X > K}} < b\, \Eb{X \1_{ X > t}}\,,$$ 
and therefore $\Eb{X \1_{X > K}} < b' \, \Eb{X \1_{K \geq  X > t}}$ with $b'=b/(1-b)$.
\end{itemize}

How would we conclude from here? 
\begin{align*}
\Es{X} &= \Eb{X\1_{t \geq X}} +  \Eb{X\1_{K \geq X >t}} + \Eb{X\1_{X > K}} \\ 
&  < A \, \Ps{X>t} + K (1+b') \, \Ps{K \geq X >t} \\ 
& < (A+K(1+b'))\, \Ps{X>t}\,,
\end{align*}
and we are done.
\medskip

Now, for the proof of \iref{(A)}, let us look at
\begin{align*}
C \geq \Es{\what{X}} = \frac{ \Es{X^2} }{ \Es{X} } 
&  = \frac{ \Eb{X^2 \1_{X >K}} + \Eb{X^2 \1_{K \geq X > t}} + \Eb{X^2 \1_{t \geq X}} }
            { \Eb{X \1_{X >K}} + \Eb{X \1_{K \geq X > t}} + \Eb{X \1_{t \geq X}} }\\
&   \geq \frac{ \Eb{X^2 \1_{X>K}} }{ \Eb{X^2 \1_{X >K}}/K + K \, \Pb{K \geq X > t} + \Eb{X \1_{t \geq X}} }\,,
\end{align*}
hence 
$$
CK\, \Pb{K \geq X > t} + C\, \Eb{X \1_{t \geq X}} \geq \left(1-\frac{C}{K}\right)\,  \Eb{X^2 \1_{X>K}}\,,
$$
for $K>t$ to be fixed later. Assuming the opposite of \iref{(A)}, we have that $\Eb{X \1_{t \geq X}} \geq A \, \Pb{K \geq X > t}$, and the last displayed inequality implies that 
$$
\left(\frac{CK}{A}+C\right)\,\Eb{X \1_{t \geq X}} \geq   \left(1-\frac{C}{K}\right) \, \Eb{X^2 \1_{X>K}} \geq \frac{K}{2}\, \Eb{X \1_{X>K}}\,,
$$
whenever $K\geq 2C$. Therefore,
\begin{align*}
c < \frac{ \Eb{X \1_{X>t}} }{\Es{X}} &\leq \frac{\Eb{X \1_{K\geq X>t}} +\Eb{X \1_{X>K}} }{\Eb{X \1_{t\geq X}} }\\
& \leq  \frac{K \,\Pb{K \geq X > t}}{\Eb{X \1_{t\geq X}}}+ \frac{2 \left(\frac{CK}{A}+C\right)}{K}
 \leq \frac{K}{A} + \frac{4C}{K}\,,
\end{align*}
whenever $A \geq K$. The first inequality is due to $\Ps{\what{X} > t} > c$. By choosing $K$ then $A$ large enough (depending only on $c$ and $C$), this gives a contradiction, proving \iref{(A)}.
\medskip

Now, to prove \iref{(B)}, assume that it is not satisfied for some $b>0$ and an arbitrarily large $K>0$. 
Then
\begin{align*}
C\, \Es{X} \geq \Es{X^2} \geq \Eb{X^2 \1_{X > K}} \geq K\, \Eb{X \1_{X >K}} \geq bK\, \Eb{X\1_{X>t}}\,.
\end{align*}
For large enough $K$, this contradicts the bound $\Ps{\what{X}>t}>c>0$, and we are done.
\qed
\medskip

\noindent{\bf Deriving (\ref{e.hatWisOK}).}
Recall that we want to show that $\Pb{\what{N\1_\thin} > 1/r } > c_2 > 0$, uniformly in $r$ and $R$. Because of the monotonicity in $r$, it is enough to prove this for some fixed $r=r_0$ (say, $r_0=2$). We obviously have
$$
\frac{\Eb{N \1_{N > 1/r_0} \1_\thin}}{\Eb{N \1_\thin}} \geq \frac{r_0^{-1}\,\Pb{N>1/r_0, \thin}}{\E N}.
$$
We have already noted that Proposition~\ref{lemgoodm} implies that $\Ps{\thin \md N>1/r_0} > c >0$, hence the numerator is at least $c\,r_0^{-1} \, \Ps{N>1/r_0}$. For the denominator, in Lemma~\ref{l.RSB} we have proved that $\E N < C' \,\Ps{N>1/r_0}$. Thus we get that the ratio is at least $c\, r_0^{-1} / C'$, and we are done.\qed

\subsection{Reconnection from thinned configurations}\label{ss.thinned}

The missing ingredient in the proof of Proposition~\ref{lemvnthin} at the end of Subsection~\ref{ss.skeleton} is (\ref{e.Visgood}), namely: 
\begin{proposition}[Things fall apart]\label{p.Visgood}
For some  $g(r)\to\infty$ as $r \to \infty$, we have that
$$
\Pbo{\th}{\thinrv > g(r) \md \thinrv > 1/r \, , \thin } > c_3 > 0 \, .
$$
\end{proposition}

The main step in proving this proposition is:

\begin{proposition}[The centre cannot hold]\label{p.fallone}
Consider dynamical percolation in $B_n$ with an initial condition in which only the hexagons intersecting the $x$-axis are open. Then, for some function $g:\N \longrightarrow \R^+$ satisfying $g(r) \to \infty$ as $r \to \infty$, the probability that at some time between $1/(2n)$ and $g(2n)$ there exists an open path realizing
$\orighex \lra n$ is bounded away from one, uniformly in $n$.  
\end{proposition}

\proofof{Proposition \ref{p.Visgood} assuming Proposition~\ref{p.fallone}} Recall the dynamics $\pthin$ specified after Definition~\ref{defthinning}, and note that 
\begin{eqnarray*}
\Pbo{\th}{\thinrv \leq g(r) \md \thinrv > 1/r \, , \thin } & = & \Pbo{\mathrm{thin}}{\, \exists \, t\in[1/r,g(r)] : 0 \build\longleftrightarrow_{}^{\omega_t}   r  \md \thinrv>1/r \, , \thin } \\
& \leq & \Pbo{\mathrm{thin}}{\, \exists \, t\in[1/r,g(r)] : 0 \build\longleftrightarrow_{}^{\omega_t}  r/2  \md \thinrv > 1/r  \, , \thin }\,.
\end{eqnarray*}

Under $\P_{\mathrm{thin}}\big( \cdot \, \vert \, \thin \big)$, the starting configuration $\omega_0$ is specified in Definition~\ref{defthin}; inside $B_{r/2}$, this is a deterministic configuration with only the hexagons intersecting the $x$-axis being open. Since a point mass trivially satisfies the static FKG inequality, we can apply the dynamical FKG inequality Lemma~\ref{lemdynfkg} for $\pthin(\cdot\mid \thin)$ inside $B_{r/2}$. Namely, for any $s\in[0,1]$, consider the dynamical event 
$$
\cA_s:=\Big\{ [0,\infty) \mathop{\longrightarrow}\limits^\omega \{0,1\}^{B_{r/2}}\text{ c\`agl\`ad} : \Pb{\thinrv > 1/r \md \omega , \thin } \geq s \Big\}
$$
in $B_{r/2}$. This event is decreasing, so that Lemma~\ref{lemdynfkg} tells us that it is negatively correlated with the increasing event $\{\exists \, t\in[1/r,g(r)] : 0 \build\longleftrightarrow_{}^{\omega_t} r/2 \}$, that is, 
\begin{eqnarray*}
 & & \PBo{\mathrm{thin}}{\cA_s \cap \big\{ \exists \, t\in[1/r,g(r)] :  0 \build\longleftrightarrow_{}^{\omega_t} r/2 \big\} \md \thin } \\
 & \leq & \Pbo{\mathrm{thin}}{\cA_s \md \thin } \, \Pbo{\mathrm{thin}}{ \exists \, t\in[1/r,g(r)] :  0 \build\longleftrightarrow_{}^{\omega_t} r/2 \md \thin }\,.
\end{eqnarray*}
Integrating over $s\in[0,1]$ gives
\begin{eqnarray*}
 & & \PBo{\mathrm{thin}}{\{T>1/r\} \cap \big\{ \exists \, t\in[1/r,g(r)] :  0 \build\longleftrightarrow_{}^{\omega_t} r/2 \big\} \md \thin }  \\
 & \leq & \Pbo{\mathrm{thin}}{ T>1/r \md \thin } \, \Pbo{\mathrm{thin}}{ \exists \, t\in[1/r,g(r)] :  0 \build\longleftrightarrow_{}^{\omega_t} r/2 \md \thin }\,.
\end{eqnarray*}
Summarising,
$$
\Pb{\thinrv \leq g(r) \md \thinrv > 1/r \, , \thin } 
\leq \Pbo{\mathrm{thin}}{\, \exists \, t\in[1/r,g(r)] :  0 \build\longleftrightarrow_{}^{\omega_t} r/2 \md \thin} \, .
$$
By Proposition~\ref{p.fallone}, the right hand side is bounded away from one, uniformly in $r$. This completes the proof. \qed
\medskip


\proofof{Proposition~\ref{p.fallone}} 
Let $\horiz_n$ denote the set of hexagons in $\boxx_n$ intersecting the $x$-axis. The elements of $\horiz_n$ will be labelled $\big\{ h_i: i \in \{-n,\ldots,n\} \big\}$ by the $x$-coordinate of the triangular lattice point at the centre of the hexagon.
We let $\boxx_n \setminus \horiz_n = \upp_n \cup \low_n$ decompose $\boxx_n \setminus \horiz_n$ into its two components above and below the $x$-axis. 
The domain $U_n$ has the shape of  a half-hexagon, whose inner boundary naturally decomposes into four paths of hexagons, each along a straight line segment: $H_n\cup\ell^1_n\cup\ell^2_n\cup\ell^3_n$, where $\ell^2_n$ denotes the horizontal path of hexagons on the top side of $U_n$ (so that the ``corner'' hexagons containing the points given in complex coordinates by $n e^{i\pi/6}$ and $n e^{i \pi/3}$ belong to $\ell_2^n$).

We will denote by $\P_{\horiz_n}$ the dynamical percolation process of the proposition, under which only elements of $\horiz_n$ are open at time $0$.

\begin{figure}[htbp]
\SetLabels
(1.02*.5){$x$}\\
\endSetLabels
\centerline{
\AffixLabels{
\epsfysize=2 in \epsffile{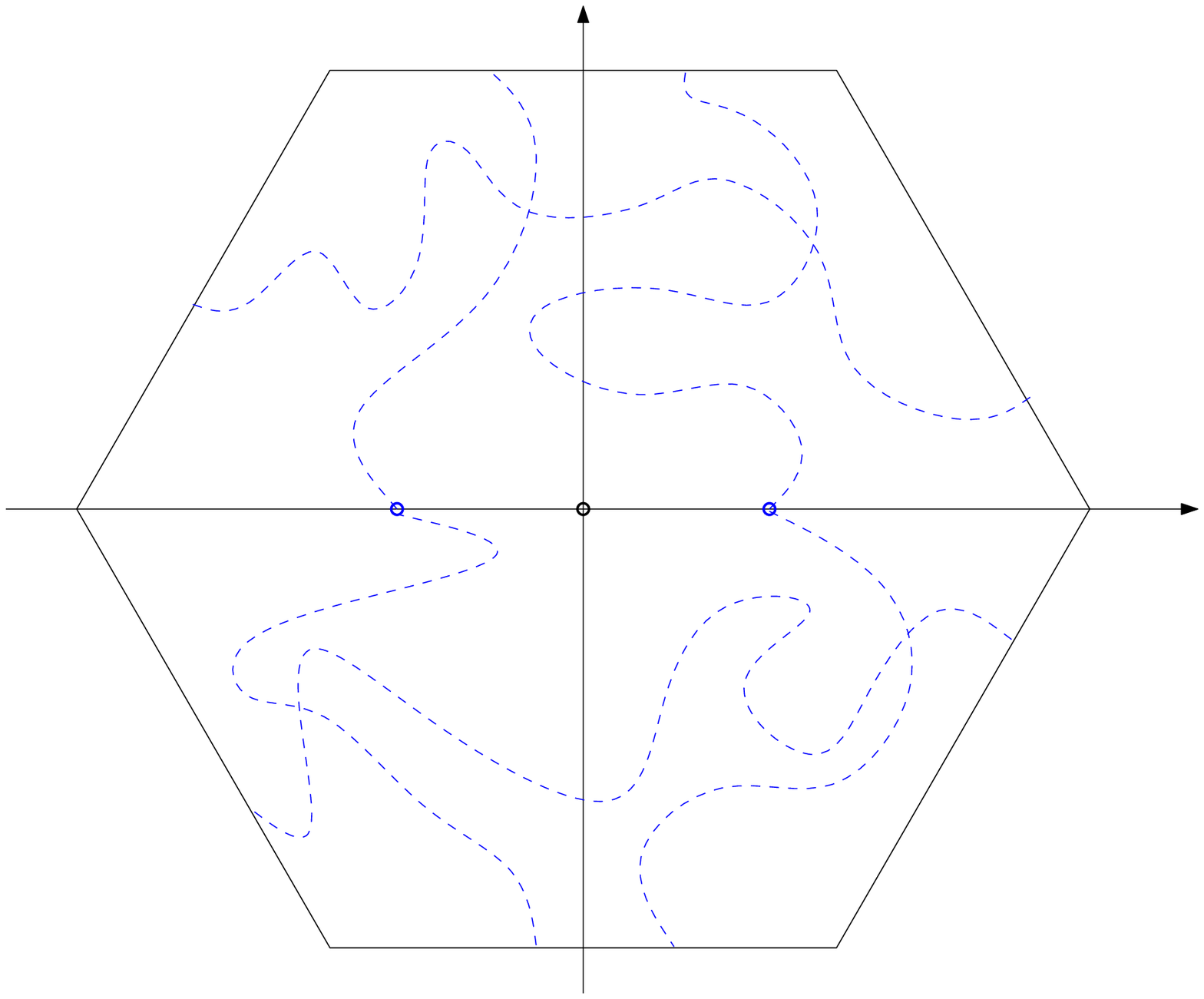}
}
}
\caption{The events $\clup$, $\cllow$, $\tmall^n_+$, $\tmall^n_-$.}
\label{f.CSULnew}
\end{figure}

Let $\clup$ denote the event that there is a closed path in $\upp_n$ from $\ell^1_n$ to $\ell^3_n$. For each $i \in \{-n/2,\ldots,n/2 \}$, let $\smallup_{h_i}$ denote the event that there is a closed path $\gamma$ in $\upp_n$ from a hexagon bordering $h_i$ to $\ell^2_n$. The events $\clup$ and  $\smallup_{h_i}$ have counterparts $\cllow$ and  $\smalllow_{h_i}$ defined verbatim after reflection in the $x$-axis. 
Finally, define 
$$
\tmall^n_+:= \big\{ \exists \, i \in \{0,\ldots,n/2\}:   h_i \textrm{ is closed}, \ \smallup_{h_i},\  \smalllow_{h_i} \big\}
$$
and
$$
\tmall^n_-:= \big\{ \exists \, i  \in \{ -n/2,\ldots,0 \} : h_i \textrm{ is closed}, \ \smallup_{h_i},\  \smalllow_{h_i} \big\} \, .
$$
Figure~\ref{f.CSULnew} illustrates that, for any $t \in (0,\infty)$,
\begin{equation}\label{e.CSUL}
 \big\{ \omega_t \in \clup \cap \cllow \cap \tmall^n_+ \cap \tmall^n_- \big\} \subseteq  \big\{ 0 \,\,\, \not\!\!\!\llra n \textrm{ in }\omega_t \big\}.
\end{equation}


Given (\ref{e.CSUL}), Proposition \ref{p.fallone} will easily follow from the next two lemmas. 

\begin{lemma}\label{l.subcrit}
For each $t > 0$,
$$
\Pbo{\horiz_n}{\bigcap_{0 < s < t} \big\{ \omega_s \in  \clup \cap \cllow \big\} } \to 1
$$
as $n \to \infty$.  
\end{lemma}
\proof Initially the set of open hexagons in $U_n$ is empty; thus, at time $s$, it has the law of a Bernoulli percolation $\P_{\tfrac{1}{2}(1 - e^{-s})}$. For $0<s<t<\infty$, let $\omega_{s,t}$ denote 
the configuration in which a hexagon is open if the hexagon is open under $\P_{H_n}$ at some time during $[s,t]$. Note then that the marginal law of $\omega_{s,t}$ in $U_n$ is a percolation whose parameter is at most 
$\tfrac{1}{2}(1 - e^{-s}) + \tfrac{1}{2}(1 + e^{-s})(1 - e^{-(t-s)})$. For any given $s > 0$, the percolation parameter of $\omega_{s,s + e^{-s}/2}$ is subcritical. By a standard subcritical percolation estimate, then, for each $s > 0$,
$\Pbo{\horiz_n}{\bigcap_{s < t < s +  e^{-s}/2} \big\{ \omega_t \in \clup  \big\} } \to 1$. By a union bound over at most $2s e^s$ sets, we see that 
$\Pbo{\horiz_n}{\bigcap_{0 < t < s} \big\{  \omega_t \in \clup  \big\}  } \to 1$. The statement of the lemma follows by symmetry in the $x$-axis. \qed

\begin{lemma}\label{l.uplow}
There exists $c > 0$ such that, for all $C > 0$ and for all $n$ sufficiently high,  
$$
 \PBo{\horiz_n}{ \bigcap_{1/n \leq t \leq C} \big\{ \omega_t \in \tmall^n_+ \big\} } \geq c \, .
$$
\end{lemma}
\noindent{\bf Proof.} We will argue that, for some $c > 0$, and for all $n$,
\begin{equation}\label{e.timelow}
 \PBo{\horiz_n}{ \bigcap_{1/n \leq t \leq c} \big\{ \omega_t \in \tmall^n_+  \big\} } \geq c \, ,
\end{equation}
and also that, for any $s \in (0,\infty)$,
\begin{equation}\label{e.timehigh}
 \lim_{n \to \infty} \PBo{\horiz_n}{ \bigcap_{s \leq t \leq s + e^{-s}/2} \big\{ \omega_t \in \tmall^n_+  \big\} } =1 \, .
\end{equation}
Note that (\ref{e.timelow}) and (\ref{e.timehigh}) prove the lemma.

Note that each of the percolations $\omega_t$ for $s \leq t \leq s + e^{-s}/2$ is stochastically dominated in $U_n \cup L_n$ by $\omega_{s,s+e^{-s}/2}$ which, as we just noted, is a subcritical percolation, of parameter $p_s < 1/2$.

Let $Q_n$ denote the set of hexagons  in $\hexlatt$ that lie in the upper-half plane and that intersect the rectangle with vertices $-n^{1/4} e_1$, $n^{1/4} e_1$, $-n^{1/4} e_1 + \tfrac{n}{2} e_2$ and $n^{1/4} e_1 + \tfrac{n}{2} e_2$. Let $\cR_n^+$ denote the event that there exists a closed path in $Q_n$ from a hexagon on the top side of $Q_n$ to one that borders $h_0$. By \cite[Theorem 11.55]{Grimmett}, for any $p < 1/2$, $\liminf_n \P_p(\cR_n^+) > 0$.

Let $\cR_n^-$ denote the event $\cR_n^+$ defined after reflection in the $y$-axis, and let $\cR_n = \cR_n^+ \cap \cR_n^- \cap \{ h_0 \, \textrm{is closed} \}$. 
Clearly, $c = \liminf_n \P_p(\cR_n) > 0$.  By partitioning $\{0,\ldots,n/2\}$ into order $n^{1/4}$ disjoint intervals and considering the analogue of $\cR_n$ for each one, we see that 
$$
\Ps{ \exists \, i \in \{0,\ldots,n/2\}:   h_i \textrm{ is closed}, \ \smallup_{h_i} \cap \smalllow_{h_i} \, \textrm{under $\omega_{s,s+e^{-s}/2}$}} 
\geq 1 - \big( 1 - c_s \big)^{n^{1/4}} \, ,
$$
where for each $s > 0$, $c_s > 0$.
Hence, we obtain~(\ref{e.timehigh}).

It is a simple matter to verify~(\ref{e.timelow}). With a probability that is bounded away from zero uniformly in $n$, some hexagon $h_i$, $0 \leq i \leq n/2$, closes during $[0,1/n]$, and remains closed until at least time one. For some $c > 0$, the marginal of $\omega_{0,c}$ in $U_n \cup L_n$ is a subcritical percolation. Thus, $\smallup_{h_i} \cap \smalllow_{h_i}$ occurs with positive probability under all $\omega_s$ for $0 \leq s \leq c$. This verifies~(\ref{e.timelow}) and completes the proof of Lemma~\ref{l.uplow}. \qed
\medskip

\proofcont{Proposition~\ref{p.fallone}} Note that Lemma~\ref{l.uplow} has a verbatim counterpart for the event  $\tmall^n_-$. 
Combining these two lemmas with the aid of the  dynamical FKG Lemma~\ref{lemdynfkg} for the process $\P_{\horiz_n}$, and using Lemma~\ref{l.subcrit}, we find that, for any $C > 0$, the left-hand side of~(\ref{e.CSUL}) is satisfied simultaneously for $1/n \leq t \leq C$ with probability tending to one as $n \to \infty$. Hence,~(\ref{e.CSUL}) proves the result. \qed




\section{The collapse of the connection near the exceptional set}\label{s.collapse}

In this section, we address the question of how quickly the infinite cluster $\Cl_0$ in dynamical percolation disintegrates as time varies away from a typical exceptional time. In view of Theorems~\ref{t.IICq} and~\ref{t.IICa}, we may rephrase the question as how rapidly this collapse occurs at small positive times in dynamical percolation where $\omega_0$ is chosen to have the law $\iic$.
In constructing approximative local times in Section~\ref{s.localtime}, we mentioned that there are several natural measurements for how 
close a finite cluster $\Cl_0$ is to being infinite. We write $\size(\Cl_0)$ as a label for any such notion, and consider three possibilities for it: the volume $|\Cl_0|$, the radius $\sup\{\|x\|: x\in\Cl_0\}$, or the ``helpfulness''  (in providing the event $0 \lra \infty$) $\help(\Cl_0)=\MO_{\Cl_0}(\omega)$ which was defined in~(\ref{e.MO}). Using any of these notions of size, one may try to define a static percolation exponent $\sigma_\size$ that measures the robustness of the infinite cluster $\Cl_0(\omega_0)$, a dynamical percolation exponent $\delta_\size$ that measures how the size of $\Cl_0(\omega_t)$ degrades with time, and then may try to relate the two exponents, a relation that is expected to reflect the fact that the ``speed'' of the dynamical process is governed by the number of pivotals in critical percolation. We first give a rough heuristic description of such a general scaling relation; however, since the existence of classical critical exponents is known only for $\hexlatt$, our actual 
theorem will reformulate the relation in a way that does not use the existence of exponents, and is valid also for the case of $\Z^2$.


To understand the robustness of the initial infinite cluster $\Cl_0(\omega_0)$, we measure the size of its restrictions to finite balls. Thus, we define the {\bf static percolation exponents} by
\begin{equation}\label{e.sigma}
\sigma_\size:=\lim_{n\to\infty}\frac{\log \E_{\IIC}\big(\size(\Cl_0\cap B_n(0))\big)}{\log n}\,,\qquad\size\in\{\vol,\radius,\help\}.
\end{equation}

From \cite{LSW} and \cite{KestenVol} we know the existence and values of the classical critical exponents
\begin{equation*}
\frac{1}{\rho}:=\lim_{n\to\infty}\frac{-\log\Pb{\radius(\Cl_0)>n}}{\log n}=\frac{5}{48}\,,\qquad
\frac{1}{\delta}:=\lim_{n\to\infty}\frac{-\log\Pb{|\Cl_0|>n}}{\log n}=\frac{1}{2\rho-1}=\frac{5}{91}\,,
\end{equation*}
which imply, with some work, that the exponents~(\ref{e.sigma}) can be given as 
\begin{equation}\label{e.three}
\sigma_\vol=\frac{\delta}{\rho}=2-\frac{1}{\rho}\,, \qquad \sigma_\radius=\frac{\rho}{\rho}=1\,,\qquad
\sigma_\help=\frac{1}{\rho}\,.
\end{equation}
The first one was established in \cite[Theorem (8)]{KestenIIC}. The second one is a triviality. For the third one, an upper bound on $\E_{\IIC}\big(\help(\Cl_0\cap B_n)\big)$ follows from~(\ref{e.MOMA}), while a lower bound can be given by the following argument. Under $\IIC$, the smallest open circuit  $\Gamma_{n/2}$ that surrounds $B_{n/2}$ is contained in $B_n$ with a uniform probability $c>0$. When conditioning on $\omega^{B_n}$, let us restrict ourselves to the part of the probability space where $\Gamma_{n/2}\subset B_n$, condition first on $\omega^{\intg(\Gamma_{n/2})}$, and then, for $R> n$, use the bound 
\begin{align*}
\EBo{\IIC}{ \Pb{0\lra R \md \omega^{B_n}} }&\geq 
\EBo{\IIC}{  \1_{\{\Gamma_{n/2} \subset B_n\}} \, \EBo{\IIC}{ \Pb{0\lra R \md \omega^{B_n}} \md \omega^{\intg(\Gamma_{n/2})} } }\\
&\geq \EBo{\IIC}{  \1_{\{\Gamma_{n/2} \subset B_n\}} \, \Pb{n/2\lra R} }\\
& \geq c\,\Pb{n/2\lra R}
\end{align*} 
to find that 
$$
\EBo{\IIC}{\lim_{R\to\infty} \frac{\Pb{0\lra R \md \omega^{B_n}}}{\Ps{0\lra R}} } \geq 
\limsup_{R\to\infty} \frac{c\,\Ps{n/2\lra R}}{\Ps{0\lra R}} \geq c'\, \Ps{0\lra n/2}^{-1} \, .
$$
In the first inequality, we used quasi-multiplicativity to obtain the uniform boundedness ${\Pb{0\lra R \md \omega^{B_n}}}/{\Ps{0\lra R}} \leq C / \Ps{0 \lra n}$ and are thus able to apply the dominated convergence theorem; the second inequality likewise uses quasi-multiplicativity. This concludes the argument for the third equality of~(\ref{e.three}) above.

For the dynamical scaling relation, we will also need the static exponent for the number of pivotals for left-right and annulus crossings:
\begin{equation*}
\tau:=\lim_{n\to\infty}\frac{\log \E_{p_c} \big| \Piv_{\cA(n)} \big| }{\log n}=\lim_{n\to\infty}\frac{\log \E_{p_c} \big| \Piv_{\cA(n,2n)} \big| }{\log n} = \frac{3}{4}\,,
\end{equation*}
following from~\cite{SW}, as already mentioned in Subsection~\ref{ss.background}.

Now, we define the {\bf dynamic percolation exponents} by 
\begin{equation*}
\delta_\size = \inf \Big\{ y \geq 0: \liminf_{t \downarrow 0} t^y \, \size(\Cl_0(\omega_t)) = 0 \Big\}\,,
\end{equation*}
starting the process from $\omega_0$ having the law of $\IIC$. Note that this is a reasonable notion of measuring the collapse of the $\IIC$ near $\omega_0$: time 0 is a limit point of exceptional times, hence $\size(\Cl_0(\omega_t))$ is infinite along some sequence $t_n\downarrow 0$, but at typical times the cluster is finite and should indeed get smaller with time, according to the following mechanism.

As we will see, for short times $t > 0$, a fragment of the original infinite cluster $\Cl_0(\omega_0)$ survives at all times $s \in [0,t]$, with the radius of this fragment determined by the maximal scale on which a pivotal hexagon rings during $[0,t]$. As such, we expect that,
for any of the above three notions of size, 
\begin{equation}\label{t.relation}
\tau \, \delta_\size = \sigma_\size \, .
\end{equation}

In the interests of concision, we will prove this relation only when $\size = \radius$. The next theorem reformulates the relation in this case, in a way that is valid even for the case of $\Z^2$. The rest of the section is devoted to the theorem's proof.

\begin{theorem}\label{t.collapse}
Consider dynamical percolation $\piic$ with $\omega_0$ having the distribution $\iic$. For $t > 0$, set $\chi_t = \inf_{s \in [0,t]} \radius\big(\Cl_0(\omega_s)\big)$. We then have $\log\chi_t \sim \log \rho(1/t)$ $\P$-a.s.~as $t \searrow 0$, with $\rho(\cdot)$ introduced in~(\ref{e.offcrit}). In particular, on $\hexlatt$, we have $\chi_t = t^{-4/3 + o(1)}$. 
\end{theorem}

\proof We start by showing the upper bound on the radius, i.e., by proving that the cluster of the origin falls apart fast enough. The following lemma will be a key step.

\begin{lemma}\label{l.collapsezeta}
There exists $c > 0$ such that the following holds.
Let $t\in(0,1)$ and $r \geq \rho(1/t)$. Let $\zeta$ denote a configuration in the annulus $A_{r,2r}$.  Let $N_r^t(\zeta)$ denote the event that the conditional probability of the inner and outer boundaries of $A_{r,2r}$ not being connected by an open path at time $t$, given that $\omega_0$ in $A_{r,2r}$ equals $\zeta$, exceeds $c > 0$. Then
$\iic \big( \{\zeta: N_r^t(\zeta)\} \big) \geq c$. Moreover, the same conclusion holds for the measure $\iic(\,\cdot\,|\,O\text{ is open})$, where $O$ is any given circuit in $A_{r/4,r/2}$ surrounding $B_{r/4}$, and where $c > 0$ may be chosen independently of $O$.
\end{lemma}

\proof Let $\cC$ denote the event that $r \llra 2r$.
That 
\begin{equation}\label{e.decor}
\Ps{\cC(\omega_0) \cap \cC(\omega_t)^c} \geq c \, ,
\end{equation}
where $c > 0$ is uniform in $r \in \N$ and $t \geq 1/(r^2\alpha_4(r))$, is a standard and simple consequence of the discrete Fourier analysis approach to critical percolation, already stated as~(\ref{e.partialdecor}). 

Note that~(\ref{e.decor}) implies the statement of the lemma when $\omega_0$ has the law of critical percolation conditioned to have the crossing. To obtain the same statement when $\omega_0$ has the law $\iic$, we can apply Lemma~\ref{l.anniic}. For the case of  $\iic(\,\cdot\,|\,O\text{ is open})$, we can apply a direct analogue of Lemma~\ref{l.anniic}, using $\P_O^\infty=\Ps{\,\cdot\md O\lra \infty}$ in place of $\P_a^b$.\qed
\medskip

We want to argue that, for any $\eps>0$, we have $\P_\iic$-a.s., for all small enough $t>0$, that  
\begin{equation}\label{e.collapseub}
  \chi_t \leq \rho(1/t)\,t^{-\epsilon} \, .
\end{equation}
We define an iterative procedure in an effort to prove~(\ref{e.collapseub}).
Let $\ell_1 \in \N$ be minimal such that $2^{\ell_1} \geq \rho(1/t)$. Write $\omega_0^{(\ell_1)}$ for $\omega_0$ restricted to $A_1 : = A_{2^{\ell_1},2^{\ell_1 + 1}}$. If $N_{2^{\ell_1}}^t(\omega^{(\ell_1)})$ occurs, and if no open path connects the inner and outer boundaries of $A_1$ at time $t$, then the procedure terminates. If one or other of these conditions is unsatisfied, let $\ell^*_1$ be the minimal $\ell \geq \ell_1 + 1$ such that $A_{2^\ell,2^{\ell + 1}}$ contains an open circuit which encloses $B_{2^{\ell}}$. Set $\ell_2 = \ell^*_1 + 2$. Write  $A_2 = A_{2^{\ell_2},2^{\ell_2 + 1}}$ and denote by $\omega_0^{(\ell_2)}$ the configuration $\omega_0$ restricted to $A_2$. If $N_{2^{\ell_2}}^t(\omega^{(\ell_2)})$ occurs, and if no open path connects the inner and outer boundaries of $A_2$ at time $t$, then the procedure terminates. Otherwise, it 
continues to its next step. The generic step has a similar description to the second one.

\begin{lemma}\label{l.terminatingindex}
Let $J \geq 1$ denote the index of the step at which the procedure terminates. Then there exists $c > 0$ such that, for each $k \in \N$, $\P \big( \ell_J - \ell_1 \geq k \big) \leq \exp \big\{ - ck \big\}$. 
\end{lemma}

\proof
Note that, by Lemma \ref{l.collapsezeta}, there exists $c > 0$ such that $J = 1$
with $\piic$-probability at least $c^2$. Under the law $\piic$ given the event that either  $N_{2^{\ell_1}}^t(\omega^{(\ell_1)})$ does not occur, or $N_{2^{\ell_1}}^t(\omega^{(\ell_1)})$ occurs and $2^{\ell_1}     \build\longleftrightarrow_{}^{\omega_t} 2^{\ell_1 + 1}$, note that the conditional distribution of $\omega_0$ in $B_{2^{\ell_1 + 1}}^c$ stochastically dominates critical percolation. (This statement is true because it is valid for $\piic$ conditionally on an arbitrary configuration in $B_{2^{\ell_1 + 1}}$ that satisfies $0 \lra 2^{\ell_1 + 1}$ at time zero.) By RSW, FKG and independence on disjoint sets, each dyadic annulus with index at least $\ell_1 + 1$ independently has probability at least $c > 0$ to contain an open circuit disconnecting its boundaries. Thus, conditionally on the value of $\ell_1$, the random variable $\ell_1^* - \ell_1$ is stochastically dominated by a geometric random 
variable (which we call $X_1$). Let $O_1$ denote the innermost of the surrounding open circuits located in $A_{2^{\ell_1},2^{\ell_1 + 1}}$. Conditionally on $\omega_
0$ taking a given form on $O \cup \intg(O)$, the conditional distribution of $\omega_0$ in the exterior of $O$ is given by $\iic$ given that $O$ is open. Thus, we may apply the $\iic(\,\cdot\,|\,O\text{ is open})$ case of Lemma~\ref{l.collapsezeta} to learn that there is probability at least $c$ that   
$N_{\ell_2}(\omega^{(\ell_2)})$ occurs. Should this event not occur, or should this event occur alongside the event $2^{\ell_2} \build\longleftrightarrow_{}^{\omega_t} 2^{\ell_2 + 1}$, then, as previously, the conditional distribution of $\ell_2^* - \ell_2$ is stochastically dominated by a geometric random variable, which we call $X_2$. 

In this way, we see that $\ell_J - \ell_1$ is stochastically dominated by $\sum_{i=1}^{G_1-1}{X_i} \, + \, 2(G_1-1)$, where $G_1$ is a geometric random variable and $\big\{ X_i: i \in \N \big\}$ is an independent sequence of i.i.d.~geometric random variables. This completes the proof of Lemma~\ref{l.terminatingindex}. \qed
\medskip

Note that the inner and outer boundaries of $A_{2^{\ell_J},2^{\ell_J + 1}}$ are disconnected at time $t$. Therefore, by Lemma~\ref{l.terminatingindex}, $\chi_t \leq \rho(1/t)\, t^{-\epsilon}$ has probability at least $1 - c^{\log_2(t^{-\epsilon})}$. We have thus verified~(\ref{e.collapseub}).
\medskip

To complete the proof of Theorem~\ref{t.collapse}, it remains to argue that, $\P_\iic$-a.s., 
\begin{equation}\label{e.collapselb}
  \chi_t \geq \rho(1/t)\, t^{\epsilon}
\end{equation}
for all small enough $t>0$. To prove this, we need the following lemma. 

\begin{lemma}\label{l.iiciicr}
Let $R \in \N$. For each $\epsilon > 0$, there exists $\delta > 0$ such that if $\cA \in \sigma \{ B_R \}$ satisfies $\iic(\cA) \geq \epsilon$, then
$\iic_R(\cA) \geq \delta$. 
\end{lemma}

\proof Recalling the definitions made in (\ref{e.muA}, \ref{e.MO}), the Bayes' rule computation (\ref{e.MOIIC}), and the quasi-multiplicativity bound (\ref{e.MOMA}), we have that, for each configuration $\zeta$ in $B_R$ realizing $0 \llra R$,
$$
\frac{{\rm d} \iic}{{\rm d} \iic_R}(\zeta) = \frac{\MO_R(\zeta)}{\MA_R(\zeta)} \leq C_1 \,,
$$
with an absolute constant $C_1<\infty$. This readily implies the claim.\qed
\medskip

Starting dynamical percolation from $\IIC_R$, and using the coupling in which bits always turn off, Kesten's near-critical one-arm stability~(\ref{e.armstab}) shows that the probability of still having the connection $0\llra R$ at all times until $\big(R^{2(1-\eps)}\alpha_4(R^{1-\eps})\big)^{-1}$ is $1-o(1)$, as $R\to\infty$. By Lemma~\ref{l.iiciicr}, the same statement holds when the initial condition is $\iic$-distributed. From this,~(\ref{e.collapselb}) follows readily. This completes the proof of Theorem~\ref{t.collapse}. \qed

\eject

\ \\
{\bf Alan Hammond}\\
Department of Statistics, University of Oxford\\
\url{http://www.stats.ox.ac.uk/~hammond/}\\
\\
{\bf G\'abor Pete}\\
Institute of Mathematics, Technical University of Budapest\\ 
\url{http://www.math.bme.hu/~gabor/}\\
\\
{\bf Oded Schramm} (December 10, 1961 -- September 1, 2008)\\
Microsoft Research\\ 
\url{http://research.microsoft.com/en-us/um/people/schramm/}\\


\addcontentsline{toc}{section}{References}

\begin{thebibliography}{XXXXX00}

\def\arXiv#1#2{\href{http://front.math.ucdavis.edu/#1}{{\tt arXiv:#1 [#2]}}} 
\def\arXivo#1{\href{http://front.math.ucdavis.edu/#1}{{\tt [arXiv:#1]}}} 

\bibitem[Ahl11]{Ahlberg} D. Ahlberg. 
The asymptotic shape, large deviations and dynamical stability in first-passage percolation on cones.
{\it Preprint}, \arXiv{1107.2280}{math.PR}



\bibitem[AGdHS08]{invperctree} O. Angel, J. Goodman, F. den Hollander and G. Slade. Invasion percolation on regular trees. {\it Ann. Probab.} {\bf 36} (2008), 420--466. \arXivo{math.PR/0608132}

\bibitem[BHPS03]{BHPS} I. Benjamini, O. H\"aggstr\"om, Y. Peres, and J. E. Steif. 
Which properties of a random sequence are dynamically sensitive? 
{\it Ann. Probab.} {\bf 31} (2003), 1--34. 

\bibitem[BKS99]{BKS} I. Benjamini, G. Kalai, and O. Schramm. 
Noise sensitivity of Boolean functions and applications to percolation. 
{\it Inst. Hautes \`Etudes Sci. Publ. Math.} {\bf 90} (1999), 5--43.
\arXivo{math.PR/9811157}

\bibitem[BS98]{BSplanes} I. Benjamini and O. Schramm. 
Exceptional planes of percolation. 
{\it Probab. Theory Related Fields} {\bf 111} (1998), no.~4, 551--564. 

\bibitem[BrGS12]{exclusion}
E. I. Broman, C. Garban and J. E. Steif.
Exclusion sensitivity of Boolean functions.
{\it Probab. Theory Related Fields} (2012), to appear.
\arXiv{1101.1865}{math.PR} 

\bibitem[BrS06]{BrSt} E. I. Broman and J. E. Steif. 
Dynamical stability of percolation for some interacting particle systems and $\eps$-movability. 
{\it Ann. Probab.} {\bf 34} (2006), no.~2, 539--576. \arXivo{math.PR/0605641}

\bibitem[CN07]{CNconv} F. Camia and C. M. Newman. 
Critical percolation exploration path and SLE$_6$: a proof of convergence. 
{\it Probab. Theory Related Fields} {\bf 139} (2007), no.~3-4, 473--519.
\arXivo{math.PR/0604487}

\bibitem[DSV09]{DaSaVa} M. Damron, A. Sapozhnikov and B. V\'agv\"olgyi. Relations between invasion percolation and critical percolation in two dimensions. {\it Ann. Probab.} {\bf 37}  (2009), 2297--2331. \arXiv{0806.2425}{math.PR}

\bibitem[DCGP11]{DCGP} H. Duminil-Copin, C. Garban, and G. Pete.
The near-critical planar FK-Ising model. 
{\it Preprint}, \arXiv{1111.0144v3}{math.PR}.

\bibitem[Dur96]{Durrett}
R. Durrett. {\it Probability: theory and examples. Second edition.} 
Duxbury Press, 1996.


\bibitem[FNRS09]{web} L. R. G. Fontes, C. M. Newman, K. Ravishankar, and E. Schertzer.
Exceptional times for the dynamical discrete web.
{\it Stochastic Processes and their Applications} {\bf 119} (2009), no.~9, 2832--2858.
\arXiv{0808.3599}{math.PR}

\bibitem[GPS10]{GPS1}
C. Garban, G. Pete, and O. Schramm.
The Fourier spectrum of critical percolation. 
{\it Acta Math.} {\bf 205} (2010), no.~1, 19--104. 
\arXiv{0803.3750}{math.PR}



\bibitem[GS12]{GS}
C. Garban and J. E. Steif. {Lectures on noise sensitivity and percolation.}
In: {\it Probability and statistical physics in two and more dimensions}
(D. Ellwood, C. Newman, V. Sidoravicius and W. Werner, ed.).
Proceedings of the Clay Mathematical Institute Summer School and XIV
Brazilian School of Probability (Buzios, Brazil), Clay Mathematics
Proceedings 15 (2012), 49--154. \arXiv{1102.5761}{math.PR}

\bibitem[Gri99]{Grimmett}
G. Grimmett. {\it Percolation. Second edition.}
Grundlehren der Mathematischen Wissenschaften, 321. Springer, 1999.

\bibitem[H\"agPS97]{HaPeSt}
O. H{\"a}ggstr{\"o}m, Y. Peres, and J.~E. Steif.
Dynamical percolation.
{\it Ann. Inst. H. Poincar\'e Probab. Statist.} {\bf 33} (1997), no.~4, 497--528.

\bibitem[HamMP12]{HMP} A. Hammond, E. Mossel, and G. Pete.
Exit time tails from pairwise decorrelation in hidden Markov chains, with applications to dynamical percolation.
{\it Electron. J. Probab.} (2012), to appear. \arXiv{1111.6618}{math.PR}

\bibitem[HarS00a]{HaraSladeOne} T. Hara  and G. Slade. The scaling limit of the incipient infinite cluster in
              high-dimensional percolation. {I}. {C}ritical exponents. {\it J. Statist. Phys.} {\bf 99} (2000) no.~5--6, 1075--1168.

\bibitem[HarS00b]{HaraSladeTwo} T. Hara  and G. Slade. The scaling limit of the incipient infinite cluster in
              high-dimensional percolation. {II}.  {I}ntegrated
              super-{B}rownian excursion.  {\it J. Math. Phys.} {\bf 41} (2000) no.~3, 1244--1293.

\bibitem[Hof06]{Hoff} C. Hoffman. 
Recurrence of simple random walk on $\Z^2$ is dynamically sensitive. 
{\it ALEA Lat. Am. J. Probab. Math. Stat.} {\bf 1} (2006), 35--45 (electronic). \arXivo{math.PR/0503065}

\bibitem[Jar03]{Jar} A. J{\'a}rai. 
Incipient infinite percolation clusters in 2{D}.  
{\it Ann. Probab.} {\bf 31} (2003) no.~1, 444--485. 

\bibitem[Kal02]{Kallen} O. Kallenberg. 
{\it Foundations of modern probability. Second edition.} 
Probability and its Applications (New York). Springer-Verlag, New York, 2002.

\bibitem[Kes86]{KestenIIC} H. Kesten.
The incipient infinite cluster in two-dimensional percolation.
{\it Probab. Th. Rel. Fields} {\bf  73} (1986), 369--394.

\bibitem[Kes87a]{KestenScaling} H. Kesten.
Scaling relations for 2D-percolation. 
{\it Comm. Math. Phys.} {\bf 109} (1987), 109--156. 

\bibitem[Kes87b]{KestenVol} H. Kesten.
A scaling relation at criticality for 2D-percolation. 
{\it Percolation theory and ergodic theory of infinite particle systems (Minneapolis, Minn., 1984-1985)}, 203--212, IMA Vol. Math. Appl., 8, Springer, New York, 1987. 

\bibitem[Kho08]{Davar} D. Khoshnevisan. 
Dynamical percolation on general trees. 
{\it Probab. Theory Related Fields} {\bf 140} (2008), no.~1-2, 169--193.
\arXiv{0705.0140}{math.PR}

\bibitem[LSW02]{LSW} G. F. Lawler, O. Schramm, and W. Werner. 
One-arm exponent for critical 2D percolation. 
{\it Electron. J. Probab.} {\bf 7} (2002), no.~2, 13 pp. (electronic). 
\arXivo{math.PR/0108211}

\bibitem[Lig02]{Liggett} T. M. Liggett.
Tagged particle distributions or how to choose a head at random. 
{\it In and out of equilibrium (Mambucaba, 2000)}, Progr. Probab.~51, pp. 133--162. Birkh\"auser Boston, Boston, MA, 2002. 
 \url{http://www.math.ucla.edu/~tml/tagged11.ps}

\bibitem[Lig05]{LiggettBook} T. M. Liggett.
{\it Interacting particle systems.}
Reprint of the 1985 original. Classics in Mathematics. Springer-Verlag, Berlin, 2005.

\bibitem[LyP11]{LPbook} R. Lyons, with Y. Peres. {\it Probability on
trees and networks}. Book in preparation, present version is at 
\url{http://mypage.iu.edu/~rdlyons}.

\bibitem[Nol08]{Nolin} P. Nolin. 
Near-critical percolation in two dimensions. 
{\it Electron. J. Probab.} {\bf 13} (2008), paper no.~55, 1562--1623.
\arXiv{0711.4948}{math.PR}

\bibitem[PSS09]{PerSchSt} Y. Peres, O. Schramm, and J. E. Steif. 
Dynamical sensitivity of the infinite cluster in critical percolation. 
{\it Ann. Inst. Henri Poincar\'e Probab. Stat.} {\bf 45} (2009), no.~2, 491--514.
\arXiv{0708.4287}{math.PR}

\bibitem[Sap11]{SapDoesnot} A. Sapozhnikov. The incipient infinite cluster does not stochastically dominate the invasion percolation cluster in two dimensions. {\it Electron. Comm. Probab.} {\bf 16} (2011), 775--780. \arXiv{1110.5269}{math.PR}

\bibitem[Sch00]{SchSLE} O. Schramm. 
Scaling limits of loop-erased random walks and uniform spanning trees. 
{\it Israel J. Math.} {\bf 118} (2000), 221--288. 
\arXivo{math.PR/9904022}

\bibitem[SSmG11]{SchSmG}
O. Schramm and S. Smirnov, with an appendix by C. Garban. 
On the scaling limits of planar percolation. 
{\it Ann. Probab.} {\bf 39} (2011), no.~5, 1768--1814.
Memorial Issue for Oded Schramm.
\arXiv{1101.5820}{math.PR}

\bibitem[SchSt10]{SchSt}
O. Schramm and J. E. Steif. 
Quantitative noise sensitivity and exceptional times for percolation. 
{\it Ann. Math.} {\bf 171} (2010), no.~2., 619--672. 
\arXivo{math.PR/0504586}

\bibitem[Smi01]{Smi} S. Smirnov. 
Critical percolation in the plane: conformal invariance, Cardy's formula, scaling limits. 
{\it C. R. Acad. Sci. Paris S\'er. I Math.} {\bf  333} (2001), no.~3, 239--244. \arXiv{0909.4499}{math.PR}
 
\bibitem[Smi06]{SmiICM} S. Smirnov. 
Towards conformal invariance of 2D lattice models. 
In {\it International Congress of Mathematicians. Vol. II}, pages 1421--1451. 
Eur. Math. Soc., Z\"urich, 2006.
\arXiv{0708.0032}{physics.math-ph}

\bibitem[SmW01]{SW} 
S. Smirnov and W. Werner. 
Critical exponents for two-dimensional percolation. {\it Math. Res. Lett.} {\bf 8} (2001), no.~5-6, 729--744.
\arXivo{math.PR/0109120} 

\bibitem[Ste09]{Steif} J. E. Steif. 
A survey of dynamical percolation. 
In {\it Fractal geometry and stochastics IV}, Birkh\"auser, pp.~145--174, 2009. \arXiv{0901.4760}{math.PR}. 

\bibitem[Szn11]{Sznitman} A-S. Sznitman. {Topics in occupation times and Gaussian free fields.} Notes of the course ``Special topics in probability'', Spring term 2011. Z\"urich Lectures in Advanced Mathematics, EMS, Z\"urich. \url{http://www.math.ethz.ch/u/sznitman/SpecialTopics.pdf}

\bibitem[Wer09]{WWperc} 
W. Werner. 
Lectures on two-dimensional critical percolation. In {\it Statistical Mechanics}, IAS/Park City Math. Ser., 16, pp. 297--360. Amer. Math. Soc., Providence, RI, 2009. \arXiv{0710.0856}{math.PR}

\end{thebibliography}
\end{document}